\RequirePackage[english,frenchb]{babel}

\documentclass{amsart}

\usepackage{ucs}
\usepackage[utf8x]{inputenc}

\usepackage{amsmath,txfonts}
\usepackage[dvips]{graphicx}
\usepackage{enumerate}
\usepackage{xypic}
\xyoption{curve}

\newcommand{\C}{\mathbb{C}}

\newcommand{\N}{\mathbb{Z}_{>0}}
\newcommand{\Q}{\mathbb{Q}}
\newcommand{\R}{\mathbb{R}}
\newcommand{\Z}{\mathbb{Z}}

\newtheorem{thm}{Theorem}[section]
\newtheorem{cor}[thm]{Corollary}
\newtheorem{lem}[thm]{Lemma}
\newtheorem{keylem}[thm]{Key Lemma}
\newtheorem{prop}[thm]{Proposition}
\newtheorem*{conj}{Conjecture}

\newtheorem*{prop-gen-exist}{Proposition~\ref{prop-generators-exist}}
\newtheorem*{introcor}{Corollary}

\newtheorem{introtheorem}{Theorem}

\theoremstyle{definition}
\newtheorem{defn}[thm]{Definition}
\newtheorem{exam}[thm]{Example}

\newtheorem*{ackn}{Acknowledgment}
\newtheorem*{claim*}{Claim}

\DeclareMathOperator{\poly}{Poly}
\DeclareMathOperator{\ang}{Angle}
\DeclareMathOperator{\supp}{supp}
\DeclareMathOperator{\crit}{Crit}

\DeclareMathOperator{\per}{per}
\DeclareMathOperator{\Int}{int}

\DeclareMathOperator{\comp}{Comp}
\DeclareMathOperator{\area}{Area}

\DeclareMathOperator{\sector}{Sector}

\theoremstyle{remark}
\newtheorem{rem}[thm]{Remark}

\renewcommand{\theenumi}{\rm (\roman{enumi})}

\renewcommand{\theenumii}{\rm (\alph{enumii})}

\def\cA{{\mathcal{A}}}
\def\cB{{\mathcal{B}}}
\def\cC{{\mathcal{C}}}

\def\cF{{\mathcal{F}}}

\def\cH{{\mathcal{H}}}
\def\cK{{\mathcal{K}}}
\def\cL{{\mathcal{L}}}
\def\cM{{\mathcal{M}}}
\def\cN{{\mathcal{N}}}
\def\cO{{\mathcal{O}}}
\def\cP{{\mathcal{P}}}

\def\cR{{\mathcal{R}}}
\def\cS{{\mathcal{S}}}
\def\cT{{\mathcal{T}}}
\def\cU{{\mathcal{U}}}
\def\cV{{\mathcal{V}}}
\def\cW{{\mathcal{W}}}

\def\cMK{{\mathcal{MK}}}
\def\cQG{{\mathcal{QG}}}
\def\cI{{\mathcal{I}}}

\def\CC{\mathbb{C}}
\def\C{\mathbb{C}}

\def\R{\mathbb{R}}

\def\Q{\mathbb{Q}}

\def\Z{\mathbb{Z}}

\def\S{{\mathbb{R}/\mathbb{Z}}}

\def\QS{\Q / \Z}
\def\CDC{\CC \setminus \overline{\Delta}}

\def\thetap{{\theta^\prime}}

\def\fnot{{f_0}}
\def\lanot{{\lambda_0}}
\def\redschema{T}
\def\redschemanot{\redschema(\lanot)}
\def\schema{T^U}

\def\deltanot{\delta_\lanot}
\def\sigmanot{\sigma_\lanot}
\def\conecnot{\cC (\lanot)}
\def\renornot{\cR (\lanot)}
\DeclareMathOperator{\hyp}{Hyp}

\newcommand{\imunit}{i}

\newcommand{\straightening}{\chi}
\newcommand{\schemamap}{{\sigma_\lambda}}
\newcommand{\schemadeg}{{\delta_\lambda}}
\newcommand{\schemamapnot}{{\sigma_{\lambda_0}}}
\newcommand{\schemadegnot}{{\delta_{\lambda_0}}}

\def\change{\marginpar{}}
\def\changetwo{\marginpar{}}

\title{Combinatorics and topology of straightening maps I: \\
compactness and bijectivity \\
}
\author{Hiroyuki Inou}
\author{Jan Kiwi}
\address{Department of Mathematics, Kyoto University, Kyoto 606-8502, Japan.}
\urladdr{www.math.kyoto-u.ac.jp/\~{}inou/}
\address{Facultad de Matem\'aticas, Pontificia Universidad Cat\'olica de Chile.}
\thanks{The second author was partially supported by Research Network on Low Dimensional Systems, PBCT/CONICYT, Chile.}
\urladdr{www.mat.puc.cl/\~{}jkiwi}

\begin{document}
\setcounter{tocdepth}{3}


\selectlanguage{english}

\begin{abstract}
  We study the parameter space structure of degree $d \ge 3$  one
 complex variable polynomials as dynamical systems acting on $\C$.  We
introduce and study {\it straightening maps}. These maps are a natural
higher degree generalization of the ones introduced by Douady and
Hubbard to prove the existence of small copies of the Mandelbrot set
inside itself.  We establish that straightening maps are
always injective and that their image contains all the corresponding
hyperbolic systems. Also, we characterize straightening maps with
compact domain.  Moreover, we give two classes of bijective straightening maps.
The first produces an infinite collection of embedded copies of
the $(d-1)$-fold product of the Mandelbrot set in the connectedness
locus of degree $d \ge 3$.  The second produces an infinite
collection of full families of quadratic connected filled Julia sets
in the cubic connectedness locus, such that each filled Julia set is
  quasiconformally embedded.
\end{abstract}


\selectlanguage{english}
\maketitle
\tableofcontents

\section*{Introduction}
The Mandelbrot set consists of all parameters $c \in \C$ for which the filled Julia set $K(z^2+c)$ of the quadratic polynomial $z^2 + c$
is connected. 
In any reasonable picture of the Mandelbrot set $\cM$ one immediately
observes the presence of small copies of it contained in itself.
Douady and Hubbard~\cite{Douady-Hubbard-poly-like} gave a mathematical proof of this fact.
They established the existence of homeomorphisms $\chi$ from carefully chosen subsets
$M'$ of the Mandelbrot set onto the whole Mandelbrot set $\cM$. 
These homeomorphisms $\chi : M' \rightarrow \cM$ are fundamental examples of what are usually called ``straightening maps''.
 The degree $ d \ge 3$ analogue of the Mandelbrot set is called the connectedness locus $\cC(d)$ of degree $d$.
In this paper we propose a definition of straightening maps in the context of higher degree polynomial dynamics.
This involves introducing appropriate  domains  as well as 
the corresponding target sets and the maps themselves.
In contrast with quadratic polynomials, one encounters a large diversity of target sets.
Therefore, a greater variety of structures replicate all over the corresponding connectedness loci.
In particular, for cubic polynomials, the target sets include the cubic connectedness locus $\cC(3)$, 
the connectedness locus of biquadratic polynomials $\{(a,b) \in \C^2;~ K((z^2 + a)^2 + b)) \mbox{ is connected} \}$, 
the full family of quadratic connected filled Julia sets
$\cMK = \{ (c,z) \in \C^2 ;~ c \in \cM, z \in K(z^2+c) \} $, and
the cross product of the Mandelbrot set with itself $\cM \times \cM$. 
In this paper we also establish the existence of a large collection of embeddings of $\cM \times \cM$ and 
inclusions of $\cMK$ in the cubic connectedness locus. The latter restrict to quasiconformal embeddings of every quadratic filled 
Julia set.

\medskip
We proceed with a brief overview of our main results. The next section (Section~\ref{sec-outline}) contains
 a more detailed account of them, including the relevant definitions and statements.

\medskip
The first objective of this paper is to propose a definition of straightening maps (i.e., to introduce the domains, the target sets, and the maps).
We will work in the parameter space $\poly (d)$ of monic centered polynomials of degree $d \ge 2$ with complex coefficients (i.e., polynomials
of the form $z^d + a_{d-2}z^{d-2} + \cdots + a_0$). 
The {\it connectedness locus $\cC(d)$ of degree $d$} is the set formed by all the polynomials $f \in \poly (d)$ which have
 connected filled Julia set $K(f)$. 
We consider a post-critically finite polynomial  $f_0 \in \cC(d)$, with at least one Fatou critical point, and define a straightening map whose domain is, in a certain sense, centered
at $f_0$ and whose target space is introduced via {\it mapping schemata}.
Following Milnor~\cite{Milnor-hyp}, one may associate with $f_0$ a combinatorial object, called the {\it reduced mapping schema $T$ of  $f_0$},
 that encodes the dynamics along the Fatou critical orbits of $f_0$. 
{\it Polynomial maps over $T$} are certain polynomial dynamical systems which  act on a disjoint union of finitely many copies of $\C$.
The parameter space formed by all monic centered polynomial maps over $T$ is called  {\it the universal polynomial model space $\poly (T)$ of $T$}.
The notions of Fatou, Julia and filled Julia sets extend to maps in $\poly (T)$.
In $\poly(T)$, the analogue of the connectedness loci is the {\it fiberwise connectedness locus $\cC(T)$}.
The fiberwise connectedness locus $\cC(T)$ will be the target space for a straightening map
defined around $f_0$.
The domain  $\cR(\lanot) \subset \cC(d)$ will be prescribed with the aid of the rational lamination $\lanot$ of $f_0$.
The definition of the map involves extracting a {\it polynomial-like map $g$ over $T$}
from each polynomial $f \in \cR(\lanot)$ (i.e., a renormalization procedure). Polynomial-like maps over mapping schemata 
are a (straightforward) generalization of the classical ones introduced by Douady and Hubbard. 
The classical notion of hybrid conjugacy between polynomial-like maps also extends easily to our context. 
Straightening will assign to each polynomial $f \in \cR(\lanot)$ a map in the fiberwise connectedness locus 
$\cC(T)$ which is hybrid equivalent to the polynomial-like map
extracted from $f$.
However, in order to avoid ambiguities, we have to introduce ``external markings'' for polynomial-like maps over $T$.
After carefully introducing markings, we obtain the desired definition of  the corresponding 
straightening map $\chi_\lanot: \cR(\lanot) \rightarrow \cC(T)$.

Our second objective is to study basic properties of straightening maps.
In contrast with the quadratic case, higher degree straightening maps are often discontinuous (see~\cite{Inou-discont}). 
Nevertheless, straightening maps are still a natural and useful tool in higher degree polynomial dynamics.
Basic questions that one may ask are the following. Are straightening maps injective? What is the image of a given straightening map? When is the domain compact? When is the domain connected?
We will show that the straightening map $\chi_\lanot: \cR(\lanot) \rightarrow \cC(T)$, briefly described above, is injective and its image contains all
hyperbolic dynamical systems in the fiberwise connectedness locus $\cC(T)$. Moreover, we give a 
characterization of  straightenings with compact domains. These straightenings include those that arise from
``primitive'' hyperbolic post-critically finite polynomials.

Straightening maps that arise from ``primitive'' hyperbolic post-critically finite polynomials are natural candidates to
be {\it onto} the connectedness locus of the corresponding model space.
We say that a polynomial  $f_0 \in \cC(d)$ is {\it primitive} if, for all distinct and bounded Fatou components $U, V$ of $f_0$, we have that
$\overline{U} \cap \overline{V} = \emptyset$.
Our third objective is to provide supporting evidence for the following conjecture:

\begin{conj}
If $f_0$ is a primitive hyperbolic post-critically finite polynomial with rational lamination $\lanot$ and reduced mapping schema $T$, 
then $\chi_\lanot : \cR(\lanot) \to \cC(T)$ is a bijection with compact and connected domain.
\end{conj}

We establish the above conjecture for two large classes of cubic straightening maps: 
``primitive disjoint type'' and ``primitive capture type''.
Moreover, we show that the primitive disjoint type gives rise, via the inverse of straightening,
to homeomorphically embedded copies of $\cM \times \cM$ in the cubic connected locus.
Similarly, each primitive capture type, leads to an injective map from the full family of quadratic filled Julia sets $\cMK$ into the cubic 
connectedness locus.
The inclusions of $\cMK$ in $\cC(3)$ have the nice extra property that restrict to a quasiconformal
 embeddings of every fiber $\{ c \} \times K(z^2+c)$.


Our results also show that 
primitive disjoint type straightenings are bijections in any degree
(which establishes the existence
of a large collection of dynamically 
defined subset of  $\cC(d)$ homeomorphic to $\cM^{d-1}$) and,
 our techniques might generalize
to degree $d \ge 4$ versions of primitive capture type straightenings.
However, already for cubic polynomials,  the above conjecture is open  since 
our techniques do not fully apply for the two remaining types of cubic straightening maps (adjacent and bitransitive).

\bigskip
In~\cite{Buff-Henriksen}, Buff and Henriksen showed that $K(z^2 +c)$ is quasiconformally 
embedded in the cubic connectedness locus provided that $z^2+c$ has a non-repelling 
fixed point. The Buff and Henriksen embeddings are a particular case of the ones obtained via a straightening map  of primitive capture type.
In~\cite{Inou-intertwine}, the first author showed that given a degree $d \ge 2$ polynomial $f \in \cC(d)$ 
with filled Julia $K(f)$, then there is a natural embedding
of $K(f)$ into the connectedness locus of degree $d'$, for all $d' > d$. 
Epstein and Yampolsky~\cite{Epstein-Yampolsky} obtained embedded copies of $\cM \times \cM$ by using the intertwining surgery technique.

As mentioned above, one of the main difficulties in  the study of higher degree straightening maps
stems from the lack of continuity (see~\cite{Inou-discont}). 
The other main difficulty stems from the insufficient understanding of  the global structure of $\cC(d)$, for $d \ge 3$.
In fact, the Douady and Hubbard techniques used in the  study of the image of quadratic straightening maps break down since they are based
on the  ``topological holomorphy'' property of straightening (a very strong form of continuity) and
on available knowledge about the global structure of the Mandelbrot set. 
As suggested by the previous paragraph, the problem of  describing the image of higher degree 
straightening maps has been already addressed in the literature, sometimes in another equivalent language.
One approach has been to construct the inverse of some straightening maps via the intertwining 
construction~\cite{Epstein-Yampolsky}.
The intertwining technique, which is certainly of intrinsic interest, may only solve the problem for a limited collection of straightening maps.
Buff and Henriksen applied holomorphic motions techniques in 
appropriately chosen one dimensional slices of parameter space.
Our approach involves a ``combinatorial tuning'' technique and obtaining suitable one dimensional restrictions
of  a given straightening map which behave as in the quadratic case. 
According to Lyubich~\cite{Lyubich-ql-space}, quadratic straightening may be also regarded as a holomorphic motion.
Nevertheless, we only use the local properties of quadratic straightening to obtain surjectivity rather
than the global properties of the holomorphic motions used  by Buff and Henriksen.  

\begin{ackn}
 The authors would like to thank Mitsuhiro Shishikura, Peter
 Ha\"issinsky and Tomoki Kawahira for valuable discussions.
\end{ackn}

\section{Outline of results}
\label{sec-outline}
The aim of this outline is to provide the reader with a self-contained exposition of the statements of our results.
These statements involve several definitions which are introduced in paragraphs~\ref{subsec-combinatorially renormalizable} through~\ref{subsec-internal}.
Paragraphs~\ref{subsec-injectivity} through~\ref{subsec-cubic} contain the statements of theorems~B, C, D, E and F, which are the main
results of this paper.

\medskip
Recall that $\poly(d)$ denotes the space of monic centered polynomials of degree $d \geq 2$.
That is, polynomials of the form $z^d + a_{d-2}z^{d-2} + \cdots + a_0$.
Therefore, $\poly(d)$ is naturally identified with $\C^{d-1}$. 
Also, recall that we denote the connectedness locus by $\cC (d)$.
According to Branner and Hubbard~\cite{Branner-Hubbard}, $\cC (d)$ is a compact subset of $\poly(d)$.
Although the results contained in the body of the paper are stated and proved in further generality, 
for the sake of simplicity of the exposition, here we restrict our attention to straightening maps that arise from
post-critically finite polynomials. 
Thus:

\medskip
 {\it For rest of this section, we let $\fnot$ be a post-critically finite  polynomial with at least one Fatou critical point.}

\medskip
As usual we denote the Julia set by $J(\fnot)$, the Fatou set by $F(\fnot)$ and the set of critical points by $\crit(\fnot) \subset \C$. 

\subsection{Combinatorially renormalizable polynomials}
\label{subsec-combinatorially renormalizable}
Given a polynomial $f \in \cC(d)$, the {\it rational lamination $\lambda_f$ of $f$} is the equivalence
relation in $\QS$ that identifies two arguments $\theta, \thetap$ if
and only if the external rays of $f$ with arguments $\theta$ and
$\thetap$ land at a common point (e.g., see~\cite[Section 6.4]{McMullen}). \change

Post-critically finite polynomials in $\cC(d)$ are completely determined by their rational laminations~\cite{Poirier}.
All the constructions in this paper only depend on a reference rational lamination $\lanot$.
For this outline of results:

\medskip
 {\it We let $\lanot$ be the rational lamination $\fnot$.}

\medskip
We say that $f \in \cC(d)$ is {\it combinatorially $\lanot$-renormalizable} if $\lambda_f \supset \lanot$.
The set formed by all {combinatorially $\lanot$-renormalizable polynomials} is denoted by $\cC(\lanot)$.

\subsection{Mapping schemata}
\label{subsec-introschemata}
In~\cite{Milnor-hyp} a ``model'' for $\conecnot$ is suggested in terms of the 
reduced mapping schema $\redschemanot$ for $\fnot$. 

In general, a {\it (resp. reduced) mapping schema $T$} is a triple $(|T|,\sigma,\delta)$
where $|T|$ is a finite set, $\sigma$ is a map from $|T|$ into itself, and
$\delta$  is a map from $|T|$ into the integers such that
(resp. $\delta(v) \geq 2$) $\delta (v) \ge 1$ for all $v \in |T|$. The map $\delta$ is called the {\it degree function}.
The sum $1+\sum_{v \in |T|}(\delta(v)-1)$ is called the {\itshape total
degree} of $T$. \changetwo

Most of the mapping schemata considered in this paper are reduced.

The {\it reduced mapping schema
$\redschemanot$ for $\lanot$ (or $\fnot$)} is   $(|\redschemanot|, \schemamapnot, \schemadegnot)$
where:
\begin{enumerate}
\item 
$|\redschemanot| = \crit(\fnot) \cap F(\fnot)$ is the set of Fatou critical points of $\fnot$.

\item 
$\schemamapnot : |\redschemanot| \rightarrow  |\redschemanot|$ where $\schemamapnot (v) = \fnot^{\ell_v}(v)$ if
  $\ell_v = \min \{ k \in \N ;~ f^k(v) \in  |\redschemanot| \}$.

\item
$\schemadegnot :  |\redschemanot| \rightarrow \N$
where $\schemadegnot (v)$ is the local degree of $\fnot$ at $v$.
\end{enumerate}

The number $\ell_v$ as above, will be called the {\it return time of $v$ to $|\redschemanot|$} and will be consistently denoted 
by $\ell_v$. 

Note that apparently $\redschemanot$ not only depends on $\lanot$ but also in $\fnot$.
However, we will define $\redschemanot$ purely in terms of $\lanot$ in Section~\ref{subsec-renormalizations}.

\subsection{Universal polynomial model space}
\label{subsec-introuniversal}
Each reduced mapping schema $T = (|T|, \sigma, \delta)$ determines a complex affine space
$\poly(T)$ called the {\it universal polynomial model space for $T$}.
More precisely, $\poly(T)$ consists of all
maps $f$ from $|T|\times \C$ to itself such that 
the restriction of $f$ to each component $\{ v \} \times \C$ is a monic centered 
polynomial map $f_v$ of degree $\delta(v)$, taking values in $\{ \sigma(v) \} \times \C$.
For short we say that $f$ is a {\it polynomial map over $T$}.
That is,
$$\begin{array}{cccl}
f:& |T|\times \C  & \to & |T|\times \C \\
  & (v,z)        & \mapsto & (\sigma(v), z^{\delta(v)} + a_{\delta(v)-2} (v) z^{\delta(v)-2} + \cdots + a_0 (v))
\end{array}
$$
where $a_j (v) \in \C$ for all $v \in |T|$ and $0 \le j \le \delta(v)-2$.
Thus, $\poly(T)$ is naturally endowed with a complex affine structure via its parameterization by the coefficients 
 $a_j (v)$.

Given $f \in \poly(T)$, the {\it filled Julia set $K(f)$} of $f$ is the set of points
in $|T|\times \C$ with precompact forward orbit. The boundary $\partial K(f)$ is called the {\it Julia set $J(f)$ of $f$}.
We say that $K(f)$ is {\it fiberwise connected} if the intersection of $K(f)$ with every fiber $\{ v \} \times \C$ is connected.
The subset of $\poly(T)$ formed by the maps $f$ with fiberwise connected filled Julia set $K(f)$ is 
called the {\it connectedness locus $\cC(T)$ for $T$}.
It follows that $\cC(T)$ is a compact subset of $\poly(T)$.
We will denote  the $v$-fiber of the filled Julia set $K(f)$ by $K(f,v)$.
Also, we will abuse of notation and sometimes regard $f \in \poly(T)$ as a collection $f=(f_v : \C \to \C)_{v \in |T|}$ of polynomials.

Polynomial maps over a schema are a very particular case of the ``fibered polynomial dynamics'' studied by O.\ Sester in~\cite{Sester}.

\subsection{Polynomial-like maps over schemata}
\label{subsec-intropolylike}
Now we generalize the notion of
polynomial-like maps introduced by Douady and Hubbard~\cite{Douady-Hubbard-poly-like}.

Consider a  mapping schema $T = (|T|, \sigma, \delta)$.
Let $U^\prime \Subset U$ (i.e. $U^\prime$ is compactly contained in $U$) 
be open subsets  of $|T| \times \C$ which are fiberwise topological disks.
That is, for all $v \in |T|$ we have that
$U^\prime \cap (\{ v \} \times \C)=\{ v \} \times U^\prime_v $ and $U \cap (\{ v \} \times \C) =\{ v \} \times U_v $ for some
topological disks $U^\prime_v, U_v \subset \C$ such that  $U^\prime_v \Subset U_v$.
A proper and holomorphic skew product over $\sigma$,
$$\begin{array}{rccl}
g : & U^\prime & \rightarrow & U \\
    & (v, z)  &  \rightarrow & (\sigma(v), g_v (z))
\end{array}$$
is called a 
{\it polynomial-like map over $T$} 
if the degree of $g_v : U^\prime_v \rightarrow U_{\sigma(v)}$ is $\delta(v)$, for every $v \in |T|$.
The {\it filled Julia set $K(g)$} is the set of all $(v,z) \in U^\prime$ such that 
$g^n ((v,z))$ is well defined, for all $n \in \N$.
Again we will abuse of notation and sometimes simply identify $g$ with the collection $(g_v :U^\prime_v \rightarrow U_{\sigma(v)})_{v \in |T|}$.
The $v$-fiber of $K(g)$ will be denoted by $K(g,v)$. 

Simple examples of polynomial-like maps over a reduced mapping schema $T$ are obtained from polynomial maps in $\cC(T)$ after
restriction to an appropriate neighborhood of their filled Julia sets.

Two polynomial-like maps $g_0$ and $g_1$  over a mapping schema $T$ are said to be {\it hybrid equivalent} if there exists 
a fiberwise quasiconformal map $\psi$ defined on a neighborhood of~$K(g_0)$, mapping the $v$-fiber of $g_0$ into $v$-fiber of $g_1$,
 such that
$\psi \circ g_0  = g_1 \circ \psi $  and 
  $$ \dfrac{\partial \psi}{\partial \bar{z}} (v,z) \equiv 0 \mbox{ a.e. on }  K(g_0).$$

\medskip
An {\itshape analytic family of polynomial-like maps over a mapping
schema $T=(|T|,\sigma,\delta)$} is a family
$((g_{\mu,v}:U_{\mu,v}' \to U_{\mu,\sigma(v)})_{v \in 
|T|})_{\mu \in M}$ of polynomial-like maps
over $T$, parameterized by a complex manifold
$M$ such that:
\begin{enumerate}
 \item $\cU_v=\{(\mu,z);\ \mu \in M, z \in
       U_{\mu,v}\}$ and
       $\cU_v'=\{(\mu,z);\ \mu \in M, z \in
       U_{\mu,v}'\}$ are homeomorphic over $M$ to $M
       \times \Delta$.
 \item The projection from the closure of $\cU'_v$ in $\cU_v$ to
       $M$ is proper.
 \item The map $g_v:\cU_v' \to \cU_{\sigma(v)}$ given by $g_v(\mu,z)
       = (\mu, g_{\mu,v}(z))$ is complex analytic and proper.
\end{enumerate}
The {\itshape fiberwise connectedness locus} of such a family is the set
of parameters $\mu \in M$ with fiberwise connected Julia sets.

This is a straightforward generalization of the notion of analytic
family of polynomial-like maps, given by Douady and Hubbard
\cite{Douady-Hubbard-poly-like}.

\subsection{External markings}
\label{subsec-intromarkings}
We will show that every polynomial-like map $g$ over a reduced mapping schema $T$ with fiberwise connected Julia set is hybrid equivalent to 
a polynomial map over $T$ (see Theorem~\ref{introstraightening} stated in Section~\ref{subsec-introstraightening} and proved in Section~\ref{sec-poly-like}). However, in order to avoid ambiguities we are forced either to  work  with the 
moduli space of polynomial maps over $T$ (affine conjugacy classes of such maps) or to introduce {\it external markings}.
For our purpose, the latter will be more convenient. Roughly speaking, an external marking of such a polynomial-like map $g$ over $T$ is a collection
of accesses, one per fiber, to an appropriate periodic or preperiodic point in $J(g)$. 

\begin{defn}
 Let $g:U' \to U$ be a polynomial-like map over a reduced mapping schema
 $T=(|T|,\sigma,\delta)$ with fiberwise connected filled Julia set. 
 A {\itshape path to $K(g)$} is a continuous map $\gamma:[0,1] \to U'$ such that
 $\gamma((0,1]) \subset U' \setminus K(g)$ and $\gamma(0) \in J(g)$.
 We say two paths $\gamma_0$ and $\gamma_1$ to $K(g)$ are
 {\itshape homotopic}
 if there exists a continuous map $\tilde{\gamma}:[0,1] \times [0,1] \to
 U'$ such that   
 \begin{itemize}
  \item $t \mapsto \tilde{\gamma}(s,t)$ is a path to $K(g)$ for all
	$s \in [0,1]$;
  \item $\tilde{\gamma}(0,t)=\gamma_0(t)$ and 
	$\tilde{\gamma}(1,t)=\gamma_1(t)$ for all $t \in [0,1]$;
  \item $\tilde{\gamma}(s,0)=\gamma_0(0)$ for all $s \in [0,1]$.
 \end{itemize}
 An {\itshape access to $K(g)$} is a homotopy class of paths to $K(g)$.
\end{defn}

\begin{defn}
 Let $g:U' \to U$ be a polynomial-like map over a reduced mapping schema
 $T=(|T|,\sigma,\delta)$ with $K(g)$ fiberwise    connected. 
 An {\itshape external marking of $g$} is 
 a collection $\Gamma=(\Gamma_v)_{v \in |T|}$ where each $\Gamma_v$ is an access to $K(g)$, contained in $\{ v \} \times \C$,
such that $\Gamma$ is forward invariant in the following sense.
For every $ v \in |T|$ and every representative $\gamma_v \subset (\{v\}\times \C) \cap U'$ of $\Gamma_v$,
the connected component of $g(\gamma_v) \cap U'$ that contains a point of $J(g)$ is a representative of  $\Gamma_{\sigma(v)}$.

 An {\itshape externally marked polynomial-like map over $T$} is a pair
 $(g,\Gamma)$ of a polynomial-like map over $T$ and an external
 marking of it.
\end{defn}

\subsection{Standard external marking of polynomial maps in $\poly(T)$}
\label{subsec-introstandard}
Polynomial maps over a reduced mapping schema $T$ are naturally endowed with
a standard marking which we introduce with the aid of appropriately defined B\"ottcher coordinates.
More precisely, given $f \in \poly(T)$,
the standard arguments for polynomials (e.g., see~\cite{Milnor3} and \cite[Proposition~2.7]{Sester})  
generalize to prove that there exists a {\itshape B\"{o}ttcher
coordinate} $\varphi$ (at infinity) for $f$. That is,  there exists a neighborhood $U$ of $|T| \times \{\infty\}$ in $|T|
\times \hat{\C}$ and a conformal map $\varphi:U \to |T| \times \hat{\C}$
of the form 
$\varphi(v,z) = (v,\varphi_v(z))$ such that $\varphi \circ f(v,z) =
(\sigma(v), (\varphi_v(z))^{\delta(v)})$ and
$\varphi$ is tangent to the identity as $z \to \infty$ on each fiber.

Let $\Delta = \{ z \in \C; |z| < 1 \}$. 
If $f \in \cC(T)$, then a B\"ottcher coordinate can be extended
uniquely to a conformal isomorphism 
\[
 \varphi: (|T| \times \C) \setminus K(f) \to |T| \times \left( \CDC \right),
\]
which we also denote by $\varphi$.
In this case, $\varphi$ is uniquely determined by $f \in \cC(T)$.

For $(v,\theta) \in |T| \times \R/\Z$, define the {\itshape external ray
for $f$} by
\[
   R_f(v,\theta) = \{\varphi^{-1}(v,r\exp(2\pi i \theta));~ 1 < r < \infty\}
\]  
We say that {\itshape $R_f(v,\theta)$ lands at $(v,z) \in J(f)$}
if 
\[
 \lim_{r \searrow 1} \varphi^{-1}(v,r\exp(2\pi i \theta)) = (v,z).
\]
The landing theorems easily generalize to this context (e.g., see~\cite{Milnor3}).
In particular,  external rays with arguments in $\QS$ always land (at eventually periodic points).

For $f \in \cC(T)$, the collection 
\[
 (R_f(v,0))_{v \in |T|}
\]
formed by  the external rays with    angle $0$ naturally induces an external marking
of any polynomial-like map over $T$ obtained as a restriction of iterates of $f$ to appropriate domains.  
We call it the {\itshape standard external marking of $f$}.
\begin{figure}[tbh]
 \begin{center}
  \fbox{\includegraphics[width=10cm,clip]{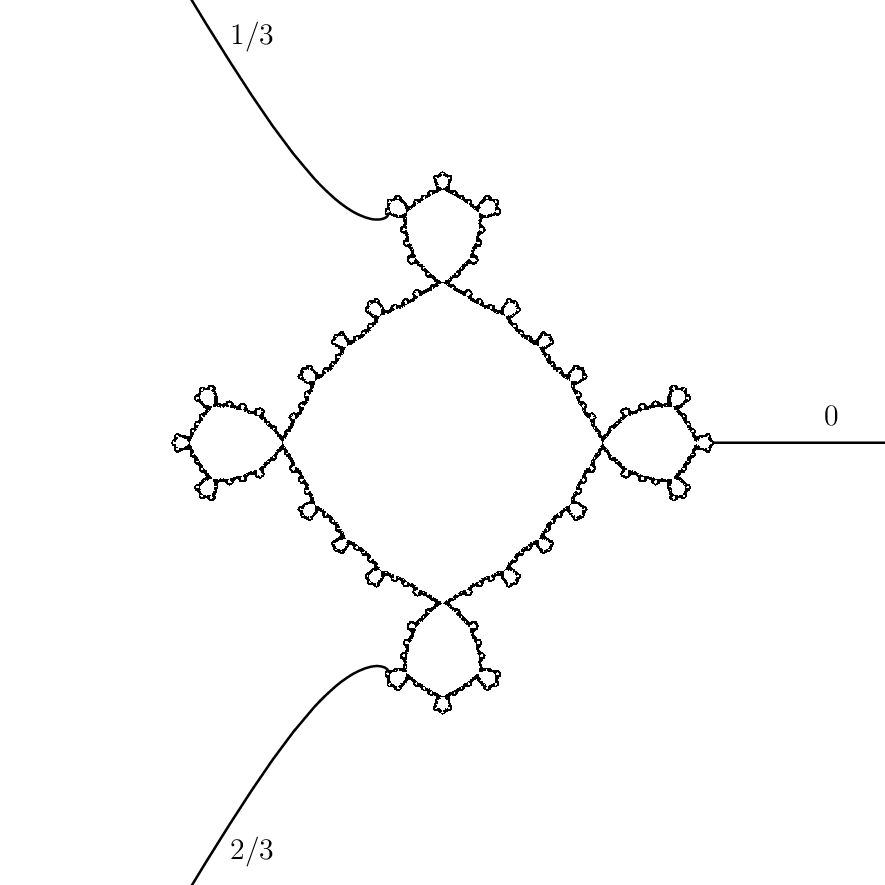}}
  \caption{The invariant rays for a quartic polynomial $z^4-1$.
  Different invariant rays represent different external markings.}
 \label{fig-ext mark}
 \end{center}
\end{figure}
Figure~\ref{fig-ext mark} illustrates three possible markings
for a 
polynomial of degree $4$, each one corresponds
to some invariant external ray. The standard marking corresponds to the
$0$-ray (i.e., the horizontal one).

\subsection{Straightening of polynomial-like maps over reduced schemata}
\label{subsec-introstraightening}
Consider two externally marked polynomial-like maps 
$(g_0:U_0' \to U_0,\Gamma_0)$, $(g_1:U_1' \to U_1,\Gamma_1)$ over $T$
with connected filled Julia sets and assume there exists a topological
conjugacy $\psi:U_0 \to U_1$ between $g_0$ and $g_1$.
We say {\itshape $\psi$ respects external markings} if 
$\Gamma_1 = ([\psi \circ \gamma_v])_{v \in |T|}$
where $(\gamma_v)_{v \in |T|}$ is a representative of $\Gamma_0$.
Note that this definition is independent of the choice of representative.

Section~\ref{sec-poly-like} contains the following generalization of Douady and Hubbard straightening theorem
to our context.

\begin{introtheorem}[Straightening]
\label{introstraightening}
Let $T$ be a reduced mapping schema and consider a polynomial-like map $g$ over $T$.
Then, there exists a polynomial map $f \in \poly (T)$ hybrid conjugate to $g$.
Moreover, 
\begin{enumerate}
\item 
If $K(g)$ is fiberwise connected, then $f \in \cC(T)$ and $f$ is unique up to affine conjugacy.

\item
If $K(g)$ is fiberwise connected and $g$ is externally marked by $\Gamma$, then for exactly one $f \in \cC(T)$ there exists a hybrid conjugacy 
that respects the external markings $\Gamma$ and $\Gamma'$, where $\Gamma'$ is the standard external marking of $f$. \change
\end{enumerate}
\end{introtheorem}

Based on the previous result which ``straightens'' a single polynomial-like map, we will define  ``straightening maps''. The domain
of a ``straightening map'' will be the set formed by ``renormalizable polynomials''.  

\subsection{Renormalizable polynomials}
Recall from Section~\ref{subsec-combinatorially renormalizable} that $\cC(\lanot)$ denotes the $\lanot$-combinatorially renormalizable locus. 
Under certain conditions it will be possible to extract a polynomial-like map over $\redschemanot$ from a
given $f \in \cC(\lanot)$. In order to be more precise, below, for every $v \in |\redschemanot|$, we identify a subset $K_f (v)$ of $K(f)$ 
which we call the {\it $v$-small filled Julia set}.

Given $\theta, \theta^\prime \in \QS$ so that the corresponding rays $R_{\fnot}(\theta), R_{\fnot}(\thetap)$ land at a common
point (i.e.,  $\theta, \thetap $ are $\lanot$-equivalent) let 
$\sector_\fnot (\theta, \theta^\prime; v)$ be the connected component of 
$\C \setminus \overline{(R_{\fnot}(\theta) \cup  R_{\fnot}(\thetap))}$
that contains $v$.
If $f \in \cC(\lanot)$, then we denote by 
$\sector_f (\theta, \theta^\prime; v)$
the connected component of
$\C \setminus \overline{(R_{f}(\theta) \cup  R_{f}(\thetap))}$
such that, for all $t \in \S$,  $$ R_{\fnot}(t) \subset \sector_\fnot (\theta, \theta^\prime; v) \mbox{ if and only if }
R_{f}(t) \subset \sector_f (\theta, \theta^\prime; v).$$
Now we define,
$$K_f (v) = K (f) \cap \bigcap_{\theta \sim_{\lanot} \thetap} \overline{ \sector_f (\theta, \theta^\prime; v)}.$$
It follows 
that $K_f (v)$ is connected and $f^{\ell_v} (K_f (v)) = K_f(\sigmanot(v)) $ (see Proposition~\ref{prop-unlinked class}).

We say that $f \in \conecnot$ is {\it $\lanot$-renormalizable} if for all $v \in \crit(\fnot) \cap F(\fnot)$ there exist
topological disks $U_v^\prime \Subset U_v$ such that $K_f (v) \subset U_v^\prime$
and $f^{\ell_v} : U^\prime_v \rightarrow  U_{\sigmanot(v)}$ is a proper map of degree $\deltanot(v)$. 
We denote the set of $\lanot$-renormalizable polynomials by $\renornot$.

Consider $f \in \renornot$ and, using the notation above,  we let 
$$U^\prime = \{ (v,z) ;~ z \in U^\prime_v \},$$
$$U = \{ (v,z) ;~ z \in U_v \}$$
and $g_v = f^{\ell_v}$. It follows that (see Definition~\ref{defn-renormalization} and Proposition~\ref{small-are-small-p}): 
$$
\begin{array}{rccl}
g: & U^\prime & \rightarrow & U \\
   & (v, z)  & \mapsto & (\sigmanot (v), g_v (z) )
\end{array}
$$
is a polynomial-like map over $\redschemanot$ with (fiberwise connected) filled Julia set
$$K(g) = \{ (v, z) ;~ z \in K_f (v) \}.$$ 
We say that $g$ is a {\it $\lanot$-renormalization} of $f$.
Note that $g$ is uniquely defined over $K(g)$, however there is a choice involved for the domain $U^\prime$.

\subsection{Internal angle systems and induced external markings}
\label{subsec-internal}
We consistently mark the polynomial-like maps over $\redschemanot$ extracted from maps 
in $\renornot$. More precisely, we introduce {\it internal angle systems} for $\fnot$ 
and describe how an internal angle system determines an external marking for a $\lanot$-renormalization of every $f \in \renornot$.

For every $v \in |\redschemanot|$ denote by $\gamma_v$ the boundary of the Fatou component of $\fnot$ that contains $v$.
Since $\fnot$ is post-critically finite, $\gamma_v$ is a Jordan curve (e.g., see~\cite{Milnor3}).
Moreover, there exists a collection $\alpha= (\alpha_v:\gamma_v \to \S)_{v \in
 |\redschema (\lanot)|}$ of homeomorphisms such that:
$$\alpha_{\schemamapnot (v)}(\fnot^{\ell_v} (z)) = \schemadegnot (v)\alpha_v(z)$$ for  all $z 
	\in \gamma_v$.

 We call $\alpha = (\alpha_v)_{v \in |\redschema (\lanot)|}$ an {\itshape internal
 angle system} for $\fnot$.

An internal angle system $\alpha = (\alpha_v)_{v \in |\redschema (\lanot)|}$   
determines  an external marking of any $\lanot$-renormalization of every $f \in \cR(\lanot)$.
  For each $v \in |\redschema (\lanot)|$ choose an argument $\theta_v$ so that the external ray of $\fnot$ with argument $\theta_v$ lands at
$\alpha_v^{-1}(0)$.
  Given a $\lanot$-renormalization $g$ of a polynomial $f \in \cR(\lanot)$,
let $\Gamma_v$ be the access with representative the connected component of 
$\overline{R_f (\theta_v)} \cap U^\prime_v$ that contains a point of $K_f (v)$. \change 
We say that $\Gamma = (\Gamma_v)_{v \in |\redschema (\lanot)|}$ is the {\itshape external marking of $g$ determined by  
the internal angle system  $\alpha$}.
We will show that this external marking  depends only on the internal angle system $\alpha$
(see~Definition~\ref{induced-marking-d} and Remark~\ref{induced-marking-r}).

\subsection{Straightening map and injectivity}
\label{subsec-injectivity}
Finally, we are ready to define the straightening maps under consideration.
More precisely, as a consequence of Theorem~\ref{introstraightening} we obtain the following result.

\begin{introcor}
Given an 
internal angle system $\alpha$ for $\fnot$ and $f \in \cR(\lanot)$, 
there exists a unique map $P \in \cC(\redschemanot)$ such that
there exists a hybrid equivalence between $P$ and a $\lanot$-renormalization $g$ of $f$,
respecting external markings, where the external marking of $g$ is the one determined by $\alpha$ and the external marking of
$P$ is the standard external marking.
\end{introcor}

With  the notation of the above corollary, 
given an internal angle system $\alpha$ for $\fnot$, 
we say that the {\it associated straightening map   $\chi_\lanot: \cR(\lanot) \to \cC(\redschemanot)$}
is the function defined by $\chi_\lanot(f) =P$.
\begin{figure}[t]
 \begin{center}
  \fbox{\includegraphics[width=6cm]{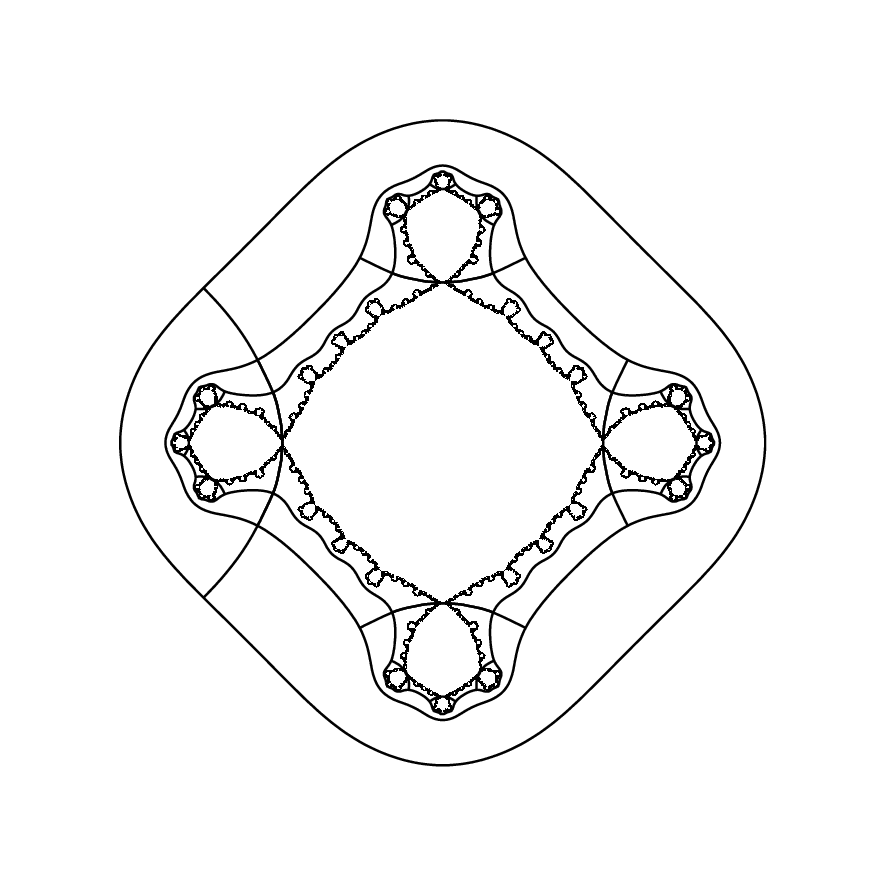}}
  \fbox{\includegraphics[width=6cm]{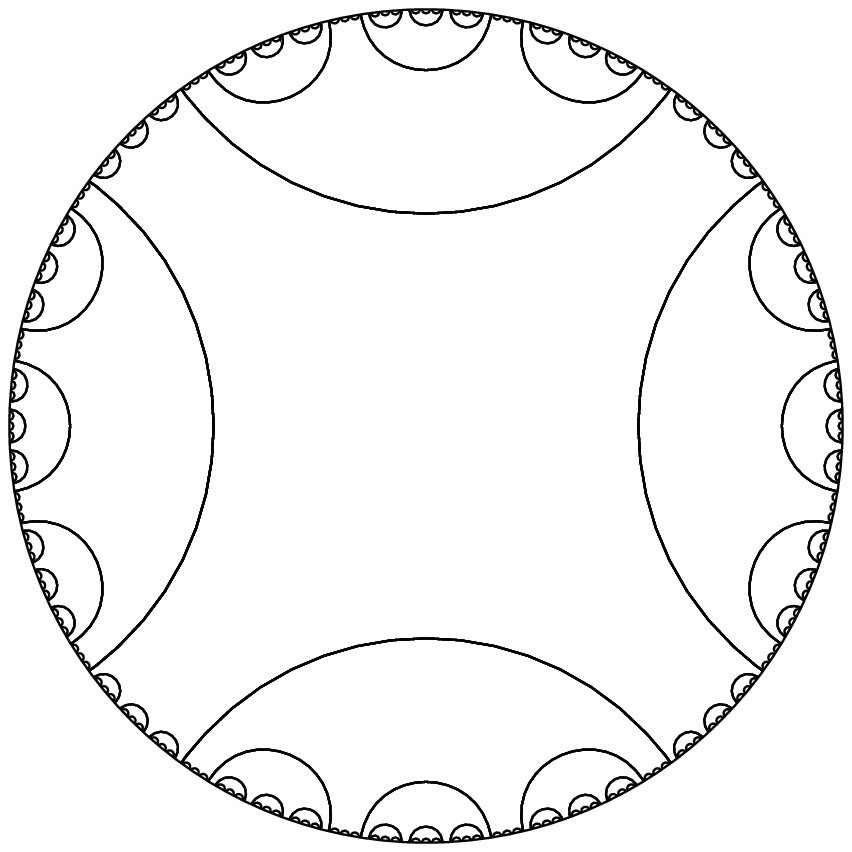}}
  \caption{The Julia set and Yoccoz puzzles, and the rational
  lamination for $z^4-1$.}
  \label{fig-quar2}
 \end{center}
\end{figure}
\begin{figure}[bth]
 \begin{center}
  \fbox{\includegraphics[width=6cm]{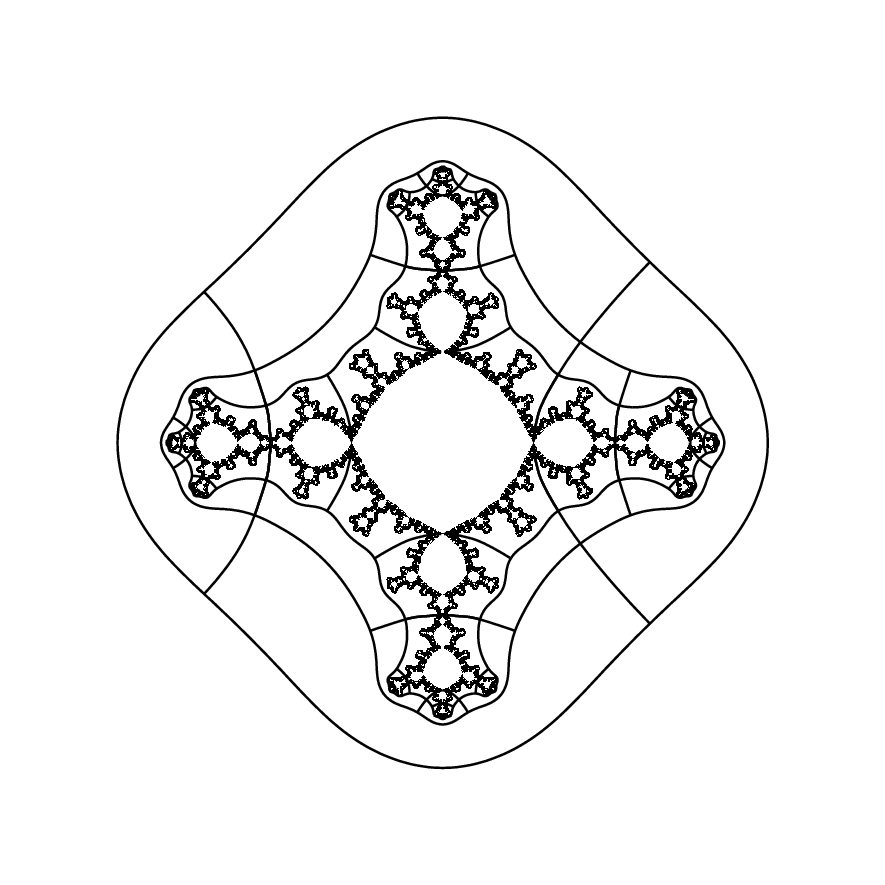}}
  \fbox{\includegraphics[width=6cm]{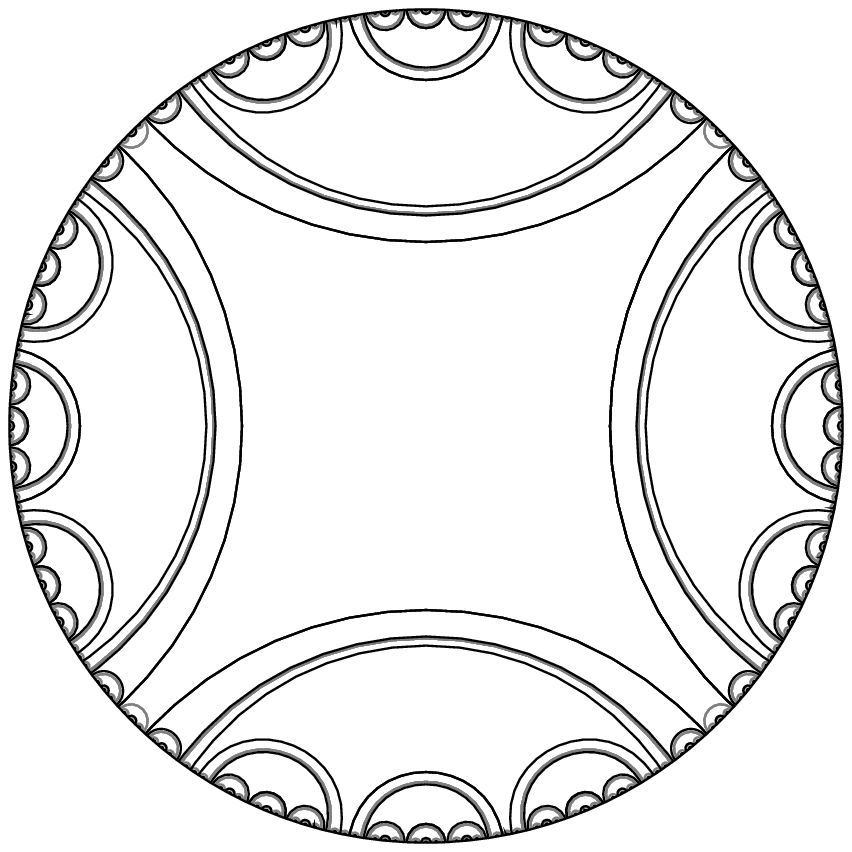}}

  \fbox{\includegraphics[width=6cm]{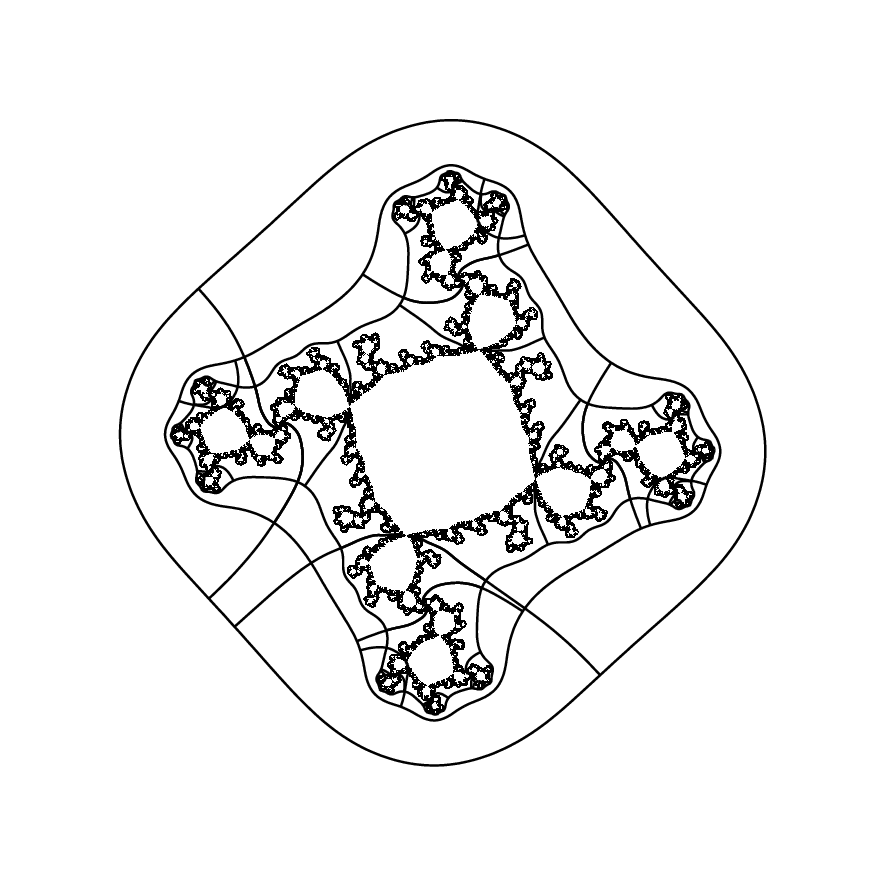}}
  \fbox{\includegraphics[width=6cm]{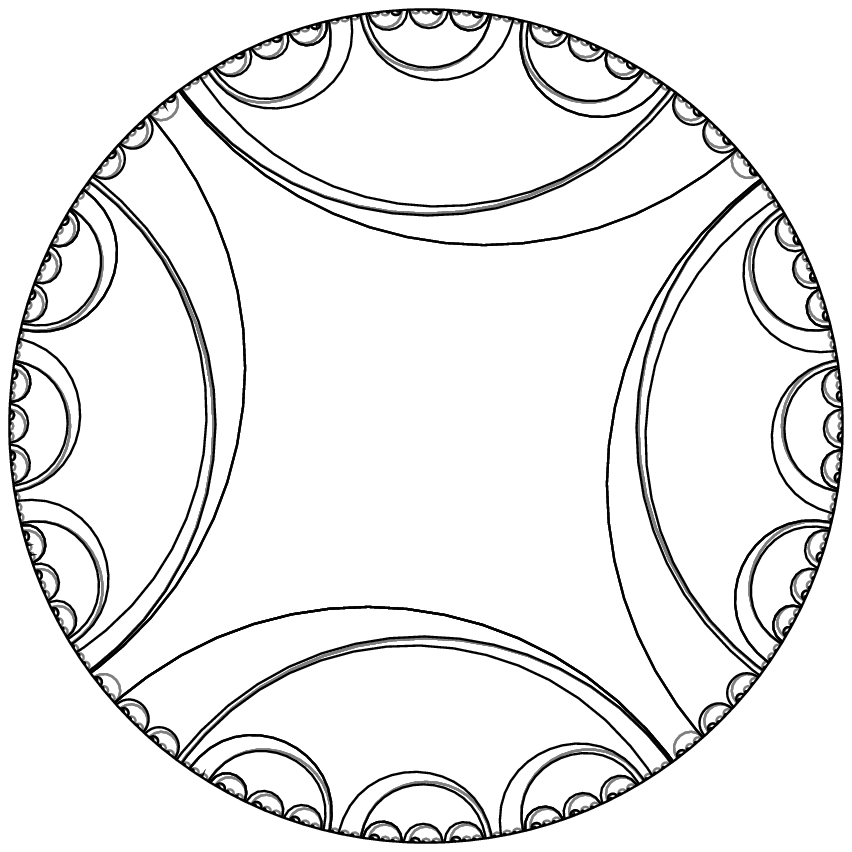}}
  \caption{Tunings (see Section~\ref{ss-onto hyperbolic}) of $z^4-1$ with itself with different external
  markings, given by the external rays of angle $0$ and $1/3$.}
  \label{fig-quar2x2}
 \end{center}
\end{figure}

\change
For example, consider $f_0(z)=z^4-1$ as in Figure~\ref{fig-quar2}.
There are three polynomials $f_1, f_2$, and $f_3 \in \renornot$ such
that $\chi_\lanot(f_i)$ is affinely conjugate to $f_0$ depending on the
choice of the external marking of $f_0$.
Figure~\ref{fig-quar2x2} shows two  $f_i$'s and the other one is just the
complex conjugate of the second one.

In Section~\ref{sec-injectivity} we prove the following result.

\begin{introtheorem}[Injectivity of Straightening]
\label{introthm-injectivity}
  Consider an internally angled  post-critically finite polynomial $\fnot$.
  Let $\chi_\lanot : \renornot \rightarrow \cC(\redschemanot)$ be the associated straightening map.
  Then, $\chi_\lanot$ is injective.
\end{introtheorem}

\subsection{Onto hyperbolic maps}
\label{ss-onto hyperbolic}
As a first step towards understanding the image of general straightening maps, 
in Section~\ref{sec-tuning} we study the action of straightening on rational laminations and the inverse of this action. That is, we study the ``combinatorial straightening'' and   ``combinatorial tuning'' procedures.
In Theorem~\ref{thm-combinatorial onto}, we establish that
straightening maps are, from a combinatorial viewpoint, (almost) surjective.
  
Following~\cite{Milnor-hyp} (also compare with~\cite{Sester}) a polynomial map $f$
over a mapping schema is called {\itshape hyperbolic} if and only if every bounded critical orbit converges to an attracting cycle of $f$.

Section~\ref{sec-onto hyperbolic} is devoted to establish the following theorem.

\begin{introtheorem}
\label{introthm-onto hyperbolic}
Let $\fnot$ be an internally angled post-critically finite polynomial
 such that $\renornot \ne \emptyset$ and
$\chi_\lanot : \renornot \rightarrow \cC(\redschemanot)$ be the associated straightening map.
Denote by $\hyp \cC(\redschemanot)$ the set of hyperbolic maps contained in 
$\cC(\redschemanot)$.
Then  $\chi_\lanot (\cR(\lanot)) \supset \hyp \cC(\redschemanot) $,
 the inverse image $\chi^{-1}_\lanot (\hyp\cC(\redschemanot))$ is a
 complex submanifold of $\poly(d)$ of dimension $d'-1$ and \changetwo
$$\chi_\lanot : \chi^{-1}_\lanot (\hyp\cC(\redschemanot)) \rightarrow \hyp\cC(\redschemanot)$$
is biholomorphic, where $d'$ is the total degree of $\redschemanot$.
\end{introtheorem}

An equivalent condition for $\renornot \ne \emptyset$ is given in
Proposition~\ref{prop-nonempty-renor} (see also
Theorem~\ref{introthm-compactness}). 

\subsection{Compactness}
In Section~\ref{sec-compact} we characterize straightening maps with compact domains.
Compactness properties of the domain of a straightening map are useful to further study its image.

Recall that we say that a  post-critically finite polynomial 
 $f \in \cC(d)$ is {\it primitive} if $f$ has at least one periodic critical point and, for every pair of distinct bounded Fatou components $V_1, V_2$
of $\fnot$, we have that $\partial V_1 \cap \partial V_2 = \emptyset$.

\begin{introtheorem}
\label{introthm-compactness}
  Let $\fnot$ be a hyperbolic post-critically finite polynomial. Then the following statements are equivalent:
  \begin{enumerate}
  \item $\fnot$ is primitive.
  \item $\conecnot = \renornot \neq \emptyset$.
\item $\renornot$ is compact and non-empty.
  \end{enumerate}
\end{introtheorem}

It is worth to mention that for non-hyperbolic post-critically finite polynomials $f_0$ the situation is delicate. In fact, $\cC(f_0)$ might not be compact even if the Fatou components involved in the renormalization do not have common boundary points with other Fatou components (see Example~\ref{exam:last}). 

\subsection{Cubic reduced mapping schemata}
There are exactly four types of  reduced mapping schemata which arise as the reduced mapping schema of a hyperbolic cubic polynomial:  

$\bullet$ $T_{\rm adj}= (\{ 0 \}, \sigma, \delta)$, where $\sigma(0)=0$ and $\delta(0) = 3$, is
called a reduced mapping schema of {\it adjacent type}. For example, if $\fnot$ is a polynomial of the form $z^3+c$ for which
the critical point $z=0$ is periodic, then the reduced mapping schema $\redschemanot$ is  $T_{\rm adj}$. 
Observe that $\poly(T_{\rm adj}) = \poly (3)$ and $\cC(T_{\rm adj}) = \cC(3)$.

$\bullet$ $T_{\rm bit} = (\{v_0,v_1 \}, \sigma, \delta)$, where
$\sigma(v_{j}) = v_{1-j}$ for $j=0,1$, and $\delta$ is constant ($=2$),
is called a reduced mapping schema of {\it bitransitive type}. If
$\fnot$ is a cubic polynomial with distinct critical points $v_0, v_1$
lying in the same periodic orbit, then $\redschemanot = T_{\rm bit}$. Note
that $\poly(T_{\rm bit})$ can be identified with  the family of biquadratic polynomials  $\poly(2\times
2)=\{ (z^2+c_0)^2 + c_1;\ (c_0, c_1) \in \C^2\} \subset \poly(4)$.  Thus
the connectedness locus of $\poly(T_{\rm bit})$ is identified with 
$\{ (c_0, c_1) \in \C^2;~ (z^2+c_0)^2 + c_1 \in \cC(4) \}$.

$\bullet$  $T_{\rm cap} = (\{v_0,v_1 \}, \sigma, \delta)$, where
$\sigma(v_{j}) = v_{0}$ for $j=0,1$, and $\delta$ is constant ($=2$),
is called a reduced mapping schema of {\it capture type}. If
$\fnot$ is a cubic polynomial with distinct critical points $v_0, v_1$
such that $v_0$ is periodic and $v_1$ is not periodic but eventually lands in the orbit of $v_0$, then $\redschemanot = T_{\rm cap}$.
In this case, $\poly(T)$ is naturally identified with
$\C^2$. In fact, given a map $f : \{ v_0, v_1 \} \times \C \to \{ v_0, v_1 \} \times \C$ in $\poly (T)$ there exists $(c_0, c_1)$ 
such that $f(v_j,z) = (v_0, z^2 + c_j)$.
Note that  $K(f)$ is fiberwise connected if and only if $c_0 \in \cM$ and 
$c_1 \in K(z^2+ c_0)$. Therefore, 
$$\cC(T_{\rm cap}) = \cMK = \{ (c , z ) \in \C^2;~ c  \in \cM, z \in  K(z^2+ c) \}.$$ 

$\bullet$  $T_{\rm dis} = (\{v_0,v_1 \}, \sigma, \delta)$ where
$\sigma(v_{j}) = v_{j}$ for $j=0,1$, and $\delta$ is constant ($=2$),
is called a reduced mapping schema of {\it disjoint type}.
If $\fnot$ is a cubic polynomial with distinct critical points $v_0, v_1$ such that
both critical points are periodic but belong to different orbits, then $\redschemanot = T_{\rm dis}$.
It follows that $\poly(T_{\rm dis})$ is identified with $\poly(2) \times \poly(2)$ and 
$\cC(T_{\rm dis}) = \cM \times \cM$.

\[
 \begin{matrix}
  \entrymodifiers={++[o][F-]}
  \xymatrix @-1pc {
  0 \ar@(rd,ru)[]_3
  } &
  \entrymodifiers={++[o][F-]}
  \xymatrix @-1pc {
  0 \ar@/^/[r]^2 & 1 \ar@/^/[l]^2
  } \\
  T_{\rm adj} \mbox{: adjacent type} &
  T_{\rm bit} \mbox{: bitransitive type} \\ \\
  \entrymodifiers={++[o][F-]}
  \xymatrix @-1pc {
  0 \ar@(ld,lu)[]^2 & 1 \ar[l]^2
  } &
  \entrymodifiers={++[o][F-]}
  \xymatrix @-1pc {
  0 \ar@(ld,lu)[]^2 & 1 \ar@(rd,ru)[]_2
  } \\
  T_{\rm cap} \mbox{: capture type} &
  T_{\rm dis} \mbox{: disjoint type} \\
 \end{matrix}
 \]

\subsection{Capture and Disjoint type straightenings}
\label{subsec-cubic}

In Section~\ref{sec-surjectivity}, 
using the compactness given by Theorem~\ref{introthm-compactness} we
are able to extend Theorem~\ref{introthm-onto hyperbolic} and describe
the image of two classes of straightening maps centered at a primitive
hyperbolic post-critically finite polynomial $\fnot$.  The first class
is the generalization of the notion of {\it disjoint type}; described
above in the context of cubic polynomials.  Namely, the reduced
mapping schema $T$ consists of $d-1$ elements which are pointwise
fixed by the schema map.  Equivalently, $\fnot$ has exactly $d-1$
superattracting periodic orbits.  Note that the universal polynomial
model space for $T$, consists of ordered $(d-1)$-tuples of monic
centered quadratic polynomials which act on $d-1$ copies of $\C$.
Thus, the corresponding connectedness locus is $\cM^{d-1}$.

\begin{introtheorem}
 \label{introthm-disjoint-surjectivity}
 Let $\fnot \in \cC(d)$ be a polynomial with exactly $d-1$ superattracting periodic orbits which is primitive. 
Then,  given any internal angle system, the associated
straightening map 
$$
\begin{array}{rccl}
  \chi : & \cR(\lambda_{\fnot}) & \to & \cM^{d-1} \\
         &   f           & \mapsto & (\chi_1(f), \dots, \chi_{d-1}(f)),
\end{array}
$$
is a homeomorphism.
Moreover, for all $1 \le i_1 < \cdots < i_k \le d-1$ and all  $c_1, \dots, c_k \in \cM$, the set
$$\{ f \in \cR(\lanot);~ \chi_{i_j}(f) = c_j \mbox{ for all } j=1, \dots, k \}$$
is contained in a codimension $k$ complex submanifold $\cS$ of $\poly(d)$.
Furthermore, if $k=d-2$, then 
$\chi_j : \cS \cap \cR(\lanot) \to \cM$ extends to a quasiconformal map in a
neighborhood of $\cS \cap \cR(\lanot)$ where $j$ is defined by
$\{i_1,\dots,i_{d-2},j\} = \{1,\dots,d-1\}$. \changetwo
\end{introtheorem}

The second class of straightening maps for which we describe the image arise from hyperbolic cubic polynomials of capture type. 
However,   before giving a precise statement of our result in this case, we need to agree on the notion
of a quasiconformal map defined on a (possibly) singular analytic space.   
\begin{defn}
   Consider a one dimensional  complex analytic space $\cS$ with singular subset $S$
and an open subset $U$ of $\C$. 
We say that a homeomorphism $\psi: \cS \rightarrow U$ is {\it $K$-quasiconformal} if 
$\psi: \cS \setminus S \to \psi(\cS \setminus S)$ is $K$-quasiconformal.
\end{defn}
Note that, in the previous definition,  $\cS \setminus S$ is a one dimensional complex manifold where the standard definition of
quasiconformal maps applies (e.g., see~\cite{Ahlfors})

\begin{introtheorem}
\label{introthm-capture-surjectivity}
Let $\fnot \in \cC(3)$ be an internally angled primitive
hyperbolic post-critically finite  polynomial with rational lamination $\lanot$ and reduced mapping schema of capture type.
 Then $\cR(\lanot)$ is connected, and the associated straightening map
$$\begin{array}{rccl}
  \chi : & \cR( \lanot ) & \to & \cMK \\
         &   f           & \mapsto & (\chi_1(f), \chi_{2}(f)),
\end{array}
$$
    is a  bijection.
Moreover, given $c \in \cM$, the set
$$ \{ f \in \cR(\lanot);~ \chi_1 (f) =c \}$$ 
is contained in a one dimensional complex analytic space
$\cS_c$ which is locally irreducible. Furthermore, 
$$\chi_2 : \cR(\lanot) \cap \cS_c \to K(z^2+c)$$ 
extends locally in $\cS_c$ to a quasiconformal map.
\end{introtheorem}



\section{Polynomial-like maps}
\label{sec-poly-like}

The aim of this section is to prove Theorem~\ref{introstraightening} stated in Section~\ref{subsec-introstraightening}.
That is, to generalize Douady and Hubbard straightening
Theorem (see~\cite{Douady-Hubbard-etude}) to the context of polynomial-like maps over reduced mapping schemata (see Sections~\ref{subsec-introschemata}, \ref{subsec-introuniversal} and~\ref{subsec-intropolylike} for the relevant definitions). 

\subsection{Proof of the straightening theorem}
Let $g: U' \to U$ be a polynomial-like map over a reduced mapping schema $T$. 
Following the ideas of Douady and Hubbard, we ``paste'' $g$ with the dynamics
of the polynomial map over $T$ given by
$f_0(v,z)=(\sigma(v),z^{\delta(v)})$. \change

Let $\bar{\Delta}(r)=\{|z|\le r\}$ denote  the closed disk of radius $r > 0$.
Choose $R>1$ and consider $V=|T| \times ( \C \setminus \bar{\Delta}(R))$.

 Restricting $g$ to a smaller domain, if necessary, we may assume
 the domain $U'$ and the range $U$ have smooth   boundaries.
 Let $h$ be a map between $A=\left( |T| \times \C \right) \setminus U' $ and $ V$ which is (fiberwise) conformal in $A \setminus \overline{U}$ such that  $h(v,z)=(v,h_v(z))$ for some quasiconformal map $h_v$, for all $v \in |T|$ and,
 $h(g(v,z))=f_0(h(v,z))$ for all  $(v,z) \in \partial U'$. \change
 Then consider 
 the map $\tilde{f}:|T|\times \C \to |T|\times \C$ defined by
 \[
  \tilde{f}(v,z) =
   \begin{cases}
    g(v,z)& \mbox{if }(v,z) \in U', \\
    h^{-1} \circ f_0 \circ h (v,z)& \mbox{if }(v,z) \in A.
   \end{cases}
 \]
Applying 
the standard pull-back argument, we obtain 
a $\tilde{f}$-invariant Beltrami differential 
$\mu_0$ on $|T|\times \C$.  
More precisely, let 
\[
 \mu_0 =
  \begin{cases}
   \sigma_0 & \mbox{on } A \setminus \overline{U} \\
   (\tilde{f}^n)^* \sigma_0 & \mbox{on } U \setminus K(g), \\
   \sigma_0& \mbox{on }K(g),
  \end{cases}
\]
where $\sigma_0$ denotes the standard complex structure (i.e. the
Beltrami differential which vanishes everywhere) . \changetwo
Since the pullback of a Beltrami differential under a holomorphic map 
does not change the essential supremum of it, the essential supremum of
$\mu_0$ is strictly less than $1$.  
From the measurable Riemann mapping theorem (see~\cite{Ahlfors}),
we obtain a fiberwise quasiconformal map
$\psi: |T|\times \hat{\C} \to |T|\times \hat{\C}$  which  conjugates
 $\tilde{f}$ with  a polynomial map
 $f$ over $|T|$. More precisely, $\psi (v,z) = (v, \psi_v(z))$ where $\psi_v$ is quasiconformal
and $\psi_v^* \sigma_0 = \mu_0$, so $f= \psi \circ \tilde{f} \circ \psi^{-1}$ is fiberwise a polynomial.
Moreover, after post-composition of $\psi$ by a fiberwise affine transformation, if necessary, we may assume that
$f$ is fiberwise monic and centered. Thus, $f \in \poly(T)$ and $f$ is hybrid conjugate to $g$.

\medskip 
 The uniqueness part of Theorem~\ref{introstraightening} 
follows from the proposition below. See Sections~\ref{subsec-intromarkings} and~\ref{subsec-introstandard} for the notion of external marking.

\begin{prop}
 If two polynomial maps $f_0, f_1 \in \cC(T)$ over a reduced mapping
 schema $T$ are hybrid equivalent, then they are affinely conjugate.

 Moreover, if $f_0, f_1 \in \cC(T)$ and a  hybrid conjugacy respects the standard external
 markings, then $f_0=f_1$.\change
\end{prop}

Before proving the proposition let us give a short discussion about external markings and
affine conjugacies. Below, we show that, modulo affine conjugacy, there is only
one external marking.

Given a polynomial map $f \in \cC(T)$
consider an external marking $\Gamma$ for (a polynomial-like
 restriction of) $f$.
 Since every access contains a unique external ray, by
 Lindel\"of's theorem, we can choose a collection of external rays 
 $(R_f(v,\theta_v))_{v \in |T|}$ as a representative of $\Gamma$.
It follows that, 
 \begin{equation}
  \label{eqn-ext-angle-system}
  \delta(v)\theta_v=\theta_{\sigma(v)}.
 \end{equation}
 Observe that there are only finitely many collections of angles
 $(\theta_v)$ satisfying \eqref{eqn-ext-angle-system}.

Now let
 \[
  A(v,z)=(v,e^{2\pi i \theta_v}z),
 \]
 and $\hat{f}=A^{-1} \circ f \circ A$. 
 It is easy to check that $\hat{f}\in \cC(T)$ and
 $A(R_{\hat{f}}(v,\theta))=R_f(v,\theta+\theta_v)$.
 In particular, $A$ maps the standard external marking for $\hat{f}$ onto 
 the  external marking $(R_f(v,\theta_v))$.

 On the other hand, it is easy to check that if a collection
 $(\theta_v)_{v \in |T|}$ of angles satisfies
 \eqref{eqn-ext-angle-system},
 then $(R_f(v,\theta_v))$ defines an external marking $\Gamma$ for $f$. 

\medskip
Now we prove the proposition and, therefore, we finish the proof of Theorem~\ref{introstraightening}.

\begin{proof}
As in the statement of the proposition we consider two polynomial maps over $T$ 
$f_0,f_1 \in \cC(T)$ 
with hybrid conjugate polynomial-like restrictions
$f_j : U'_j \to U_j$ such that $K(f_j) \subset U'_j$ where $j=0,1$.
Let us denote the hybrid conjugacy $\psi: U_0 \to U_1$.
Replacing $f_1$ by an affinely conjugate polynomial map
we may assume that $\psi$ respects the
 standard external markings.
 
 Define the bijection $\Phi: |T| \times \C \to |T| \times \C$ by
 \[
  \Phi(v,z) = 
  \begin{cases}
   \psi(v,z) & \mbox{if }(v,z) \in K(f_0), \\
   \varphi_{f_1}^{-1} \circ \varphi_{f_0}(v,z) & \mbox{otherwise.}
  \end{cases}
 \]
 Recall that $\varphi_{f_i}: (|T| \times \C) \setminus K(f_i) 
 \to |T| \times \{|z|>1\}$ is the corresponding B\"ottcher map which is biholomorphic and   
defined in the complement of $K(f_i)$ since $f_i \in
 \cC(T)$ for $i=0,1$.
 Clearly, $\Phi \circ f_0 = f_1 \circ \Phi$.
 
As in the proof of \cite[Proposition~6]{Douady-Hubbard-poly-like},
since $\psi$ respects the standard external markings, it follows that
$\Phi$ is biholomorphic, hence affine. 
The affine map $\Phi$ conjugates two monic centered polynomial
maps over $T$
 and it is tangent to the identity at infinity in every fiber.
Therefore, $\Phi$ must be the  identity.
\end{proof}


\section{Laminations}
\subsection{Laminations}
\label{subsec-laminations}
The topological dynamics of locally connected polynomial Julia sets is naturally described by the quotient of the circle under an  
equivalence relation which is invariant under multiplication by $d$ (e.g. see~\cite{Guckenheimer}).
Thurston, in the early 1980's, initiated the study of the basic properties of the equivalence relations which arise from polynomial 
dynamics introducing invariant laminations \cite{Thurston}. Since then they have proven to be
an useful object to encode the landing pattern of external rays. 
In this section, we briefly review the basic notions and results about
laminations. 
For our purposes we will need to consider
the well known $d$-invariant laminations as well
as forward invariant laminations with finite support.
The latter type of laminations are used to define
``combinatorial Yoccoz puzzles'' (compare with~\cite{Inou-lim}).

\medskip
For an integer $d>0$, we let
$$\begin{array}{rccl}
m_d: & \R/\Z & \to & \R/\Z\\
     &\theta & \mapsto & d\theta.
\end{array}$$

Two subsets $A, B \subset \R/\Z$ are said to be {\itshape unlinked}
 if $B$ is contained in a component of $\R/\Z\setminus A$
 (it is equivalent to $A$ being contained in a component of $\R/\Z
 \setminus B$). 
Note that $A,B$ are unlinked if and only if the Euclidean (or hyperbolic)
convex hulls in $\overline{\Delta} = \{ |z| \le 1 \}$ of $\exp(2 \pi i A)$ and $\exp(2 \pi i B)$ are disjoint.

\begin{defn}[Laminations without dynamics]
 Let $E \subset \R/\Z$.
 An equivalence relation $\lambda$ on $E$ is called a {\itshape
 lamination on $E$} if the following three conditions hold.
 \begin{enumerate}
  \item $\lambda$ is closed in $E \times E$.
  \item Every equivalence class is finite.
  \item Equivalence classes are pairwise unlinked.
 \end{enumerate}

The  {\itshape support $\supp(\lambda)$} of $\lambda$ is the union
 of all non-trivial $\lambda$-equivalence classes.

A lamination on $\Q/\Z$ (resp.\ $\R/\Z$) is called a {\itshape rational
 } (resp. {\itshape real}) {\itshape lamination}.
\end{defn}

We will be mainly concerned with laminations that are compatible with the
action of $m_d$.

\begin{defn}[Pushforward and $d$-invariance for laminations]
 Let $d >1$ be an integer and $E' \subset \R/\Z$.
 Let $\lambda$ and  $\lambda'$ be laminations on $E=m_d^{-1}(E')$ and $E'$,
respectively. 

We say that
 $(m_d)_*\lambda=\lambda'$ if for any $\lambda$-equivalence class
 $A$,
 \begin{enumerate}
  \item $m_d(A)$ is a $\lambda'$-equivalence class.
  \item $m_d :A \to m_d(A)$ is {\itshape consecutive preserving}, i.e.,
	for any component $(\theta,\theta')$ of $\R/\Z \setminus A$,
	$(d\theta, d\theta')$ is a component of $\R/\Z \setminus m_d(A)$.
 \end{enumerate}
A lamination $\lambda$ on a $m_d$-completely invariant set $E
 \subset \R/\Z$ 
is called {\itshape $d$-invariant} or simply {\itshape invariant}
if $(m_d)_*\lambda=\lambda$.
\end{defn}

\subsection{Invariant rational laminations: unlinked classes, fibers and critical elements}

Invariant rational laminations will be particularly useful in this
paper.  Recall that $\cC(d)$ denotes the connectedness locus in the space
of monic centered degree $d \ge 2$ polynomials $\poly(d)$.  Given a monic
centered polynomial $f \in \cC(d)$, the {\itshape rational lamination
  $\lambda_f$ of $f$} is the equivalence relation on $\QS$ which
identifies two rational arguments if and only if the corresponding
rays land at a common point.  According to~\cite{Kiwi}, an equivalence
relation $\lambda$ on $\Q/\Z$ is the rational lamination of some $f
\in \cC (d)$ if and only if $\lambda$ is a $d$-invariant
rational lamination.

Rational laminations allow us to establish when distinct polynomials
share certain dynamical features.
Therefore, it is natural to define some subsets of dynamical and parameter space
in terms of rational laminations. 
Given a $d$-invariant rational lamination $\lambda$, let
$$\cC(\lambda) = \{ f \in \cC(d); \lambda_f \supset \lambda \}.$$  
Note that this set is always non-empty~\cite{Kiwi}.

The sets $\cC(\lambda)$ have already deserved a lot of attention.
When $\lambda$ is the rational lamination of a polynomial with all
cycles repelling (equivalently, $\lambda$ is maximal with respect to
the partial order on $d$-invariant rational laminations),
$\cC(\lambda)$ is a ``combinatorial class''. In degree $2$, the MLC conjecture asserts 
that combinatorial classes are singletons. In degree $3$,
there are some combinatorial classes which are non-trivial
continua~\cite{Henriksen}. However, it is natural to conjecture that
in any degree, each combinatorial classes is contained in an analytic
subspace of $\poly(d)$ of codimension (at least) $1$ (see~\cite{Inou-intertwine}).

This paper is mainly concerned with $\cC(\lambda)$ 
for invariant rational laminations $\lambda$ which are not maximal.
In degree $2$, following Douady and Hubbard, these sets coincide with homeomorphic copies of the Mandelbrot set  $\cM$  
contained in itself. 
In order to compare $\cC(\lambda)$ with its model 
we employ a renormalization procedure.
A first step towards renormalization is to cut dynamical space into sectors and the filled Julia set into fibers.

\begin{defn}[Admissible]
  Consider a polynomial $f \in \cC(d)$. We say that a lamination $\lambda$ on $E \subset \QS$ is
{\itshape admissible for $f$} if $\lambda \subset \lambda_f$, namely, 
 for every pair of distinct arguments $\theta$ and $\theta'$ which are
$\lambda$-related, we have that they are also $\lambda_f$-related.
\end{defn}

Note that $\cC(\lambda)$ exactly consists of all
the polynomials for which the invariant rational lamination $\lambda$ is admissible.
Admissible laminations will allow us to cut the dynamical space into sectors and 
the filled Julia set into fibers (compare with~\cite{Schleicher}).

\begin{defn}[Sectors and fibers]
Let $f \in \cC(d)$ and $\lambda$ be a lamination on $E \subset \QS$ that is
{ admissible for $f$}.
  Consider a set $L \subset \S$ which is unlinked with every non-trivial $\lambda$-class. If
$\theta$ and $\theta'$ are distinct 
$\lambda$-equivalent arguments, we denote by
$$\sector(\theta,\thetap;L)$$
the connected component of $\C \setminus \overline{(R_f (\theta) \cup R_f (\thetap))}$ which contains
the external ray $R_f (t)$ for all  $t \in L$.
 (When $\theta=\theta'$, define $\sector(\theta,\theta;L) = \C \setminus
 \overline{R_f(\theta)}$.)
Moreover, we define the {\itshape $\lambda$-fiber} of $L$ by:
$$K_f (L) = K(f) \cap \bigcap_{\theta \sim_{\lambda}\thetap, \theta \neq \thetap} \overline{\sector(\theta,\thetap;L)}.$$
\end{defn}

Our main interest now is when $\lambda$ is an invariant rational lamination and the sets $L$ above are, in a certain sense, 
as large as possible. That is, when $L$ are ``unlinked classes''.

\begin{defn}[Unlinked relation and classes]  
\label{d unlinked class}
 Let $\lambda$ be a \change
rational lamination.
 We say $\theta, \theta'\in \R/\Z \setminus \QS$ are {\itshape
 $\lambda$-unlinked} if
 $\theta=\theta'$ or for any $\lambda$-equivalence class $A$, 
 $\theta$ and $\theta'$ lie in the same component of $\R/\Z \setminus A$.

 The $\lambda$-unlinked relation is in fact an equivalence
 relation.
 Equivalence classes are called {\itshape $\lambda$-unlinked classes}.
\end{defn}

Observe that, by definition,  unlinked classes are contained in $\S
\setminus \QS$.
Infinite $\lambda$-unlinked classes are closely related to the concept
of gaps introduced by Thurston.

In order to visualize the unlinked classes of an invariant rational lamination
it is convenient to consider the real extension of a rational lamination.

\begin{defn}[Real extension]
\label{defn-real extension}
  Given a lamination $\lambda$ on $E \subset \S$
 the {\itshape real extension $\hat{\lambda}$
of $\lambda$} is the smallest equivalence relation that contains the closure
$\overline{\lambda}$ of $\lambda$ in $\S \times \S$.
\end{defn}

According to~\cite{Kiwi}, the real extension of an invariant rational lamination is an invariant real lamination.

For example, consider a polynomial $f$ with locally connected Julia
set and such that the following technical condition is satisfied. No
critical point with infinite forward orbit lies in the boundary of a
bounded Fatou component (e.g., $f$ is post-critically finite).  Then,
the real extension $\widehat{\lambda_f}$ of the rational lamination $\lambda_f$ encodes the 
landing pattern of {\itshape all} the external rays of $f$. That is
$\widehat{\lambda_f}$ identifies two arguments $\theta$ and $\thetap$
if and only if the external rays $R_f(\theta)$ and $R_f (\thetap)$
land at a common point. In this case, $L$ is an infinite
$\lambda$-unlinked class if and only if there exists a bounded Fatou
component $V$ such that $L$ is the set formed by the arguments of
all the irrational rays landing at $\partial V$.  Finite unlinked classes are
irrational equivalence classes of $\widehat{\lambda_f}$ and correspond
to external rays landing at a common point, with infinite forward orbit.

\smallskip
The main properties of $\lambda$-unlinked classes and their fibers are summarized in the proposition below.
Before we state the proposition let us agree that if $X, Y \subset \C$ and $f: X \rightarrow Y$ is a surjective
map defined and holomorphic on a neighborhood of $X$, then we say that the {\itshape  degree of $f: X \rightarrow Y$} is 
$\delta \geq 1$ if every point in $Y$ has $\delta$ preimages in $X$ for
such a holomorphic extension of $f$, counting multiplicities.

\begin{prop}
 \label{prop-unlinked class}
Consider a $d$-invariant rational lamination $\lambda$ with real extension $\hat{\lambda}$ and a polynomial $f \in \cC(\lambda)$.
Let $L$ be a $\lambda$-unlinked class. Then: \change
\begin{enumerate}
\item \label{item-prop-unlinked class-invariance}
      $m_d(L)$ is a $\lambda$-unlinked class.
      In particular, if $\theta$ and $\theta'$ are $\lambda$-unlinked, then
      $d\theta$ and $d\theta'$ are also $\lambda$-unlinked.  
\item \label{item-prop-unlinked class-K}
      $f(K_f (L)) = K_f (m_d (L))$.
\item \label{item-prop-unlinked class-finite}
      If $L$ is finite, then 
  \begin{enumerate}
  \item \label{item-prop-unlinked class-wandering}
	$L$ is a $\hat{\lambda}$-class. In particular, 
	$L$ is wandering ($m_d^n(L) \ne m_d^m(L)$ for any $n \ne m$).
  \item $m_d: L \rightarrow m_d(L)$ is $\delta=\delta(L)$-to-one for some $\delta \geq 1$.
  \item \label{item-prop-unlinked class-finite-deg-f}
	The degree of $f: K_f(L) \rightarrow K_f (m_d (L))$ is well defined and exactly $\delta(L)$.
  \end{enumerate}
\item 
If $L$ is infinite, then
 \begin{enumerate}
  \item \label{item-prop-unlinked class-periodic}
      $L$ is eventually periodic (i.e., there exists $n > 0$ and $\ell
      \ge 0$ such that
$m_d^{n+\ell}(L)=m_d^{\ell}(L)$).
  \item \label{item-prop-unlinked class-complement}
	 $\overline{L}/\hat{\lambda}$ is homeomorphic to $\R/\Z$. Moreover,
if $(\theta, \thetap)$ is a connected component of $\S \setminus \overline{L}$,
then $\theta$ and $\thetap$ are $\lambda$-equivalent and belong to   $\Q/\Z$. Furthermore, either
$(m_d(\theta), m_d (\thetap))$ is a connected component of $\R/\Z
	\setminus m_d (L)$ or $ m_d(\theta)= m_d (\thetap)$.
  \item \label{item-prop-unlinked class-degree}
	$m_d:\overline{L} \to \overline{m_d(L)}$ induces a map
	$\tilde{m}_d: \overline{L}/\hat{\lambda} \to \overline{m_d
	(L)}/\hat{\lambda}$ which is an orientation preserving
	$\delta$-fold covering map for some $\delta=\delta(L) \ge 1$.
	
	In particular, we can choose a collection   of homeomorphisms
	$\alpha_L:\overline{L}/\hat{\lambda} \to \R/\Z$, one for each infinite
	$\lambda$-unlinked class $L$, such that
	$\alpha_{m_d(L)} \circ \tilde{m}_d \circ \alpha_L^{-1}$ is
	equal to the $\delta$-fold covering map
	$m_\delta:\R/\Z \to \R/\Z$ for any $L$.
  \item \label{item-prop-unlinked class-critical return}
	If $m_d^n(L)=L$, then $\delta^n(L)>1$, where
	$\delta^n(L)=\prod_{k=0}^{n-1}\delta(m_d^k(L))$.
  \item \label{item-prop-unlinked class-deg-f}
	Assume further that for any $\lambda_f$-class $A$ 
	and any $\lambda$-unlinked class $L'$ such that $A \cap
	\overline{L} \ne \emptyset$ and $m_d(L')=m_d(L)$
	we have $A \cap \overline{L'} = \emptyset$.  
	Then the degree of $f: K_f(L) \rightarrow K_f (m_d (L))$ is
	exactly $\delta(L)$. 
 \end{enumerate}
\end{enumerate}
\end{prop}

Most of the above proposition is already proved in
 \cite[Section~4]{Kiwi} and \cite[Section~3]{Inou-lim}.
The reader may find the rest of the proof at the end of this section.

\begin{figure}[tb]
 \begin{center}
  \fbox{\includegraphics[width=6cm]{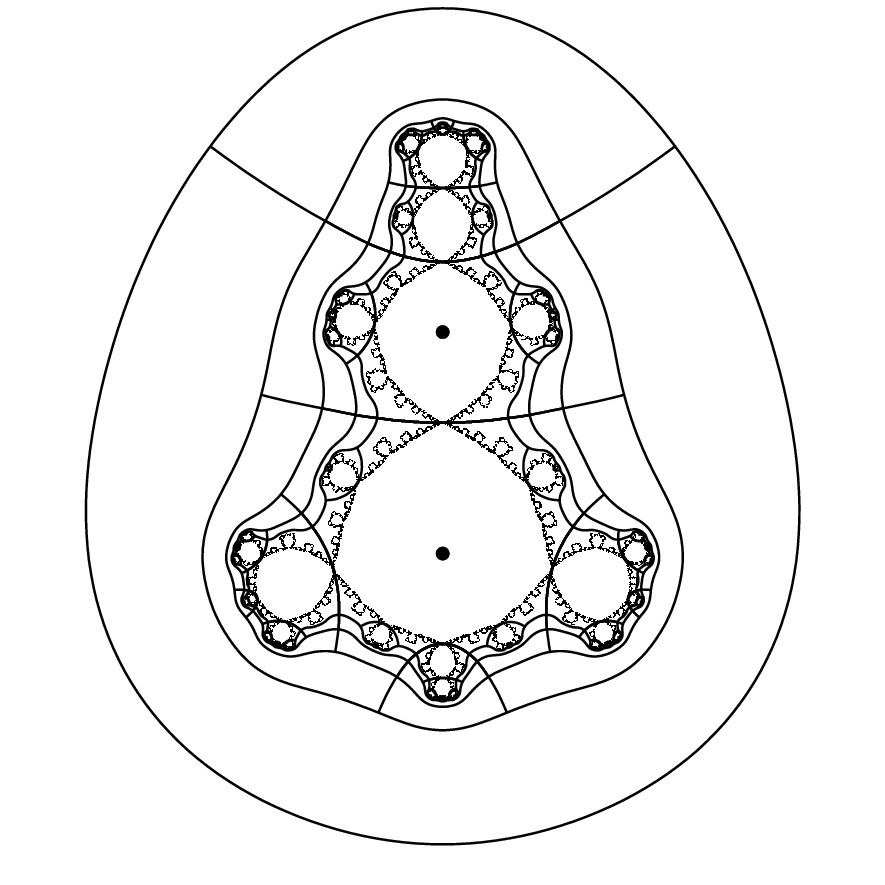}}
  \fbox{\includegraphics[width=6cm]{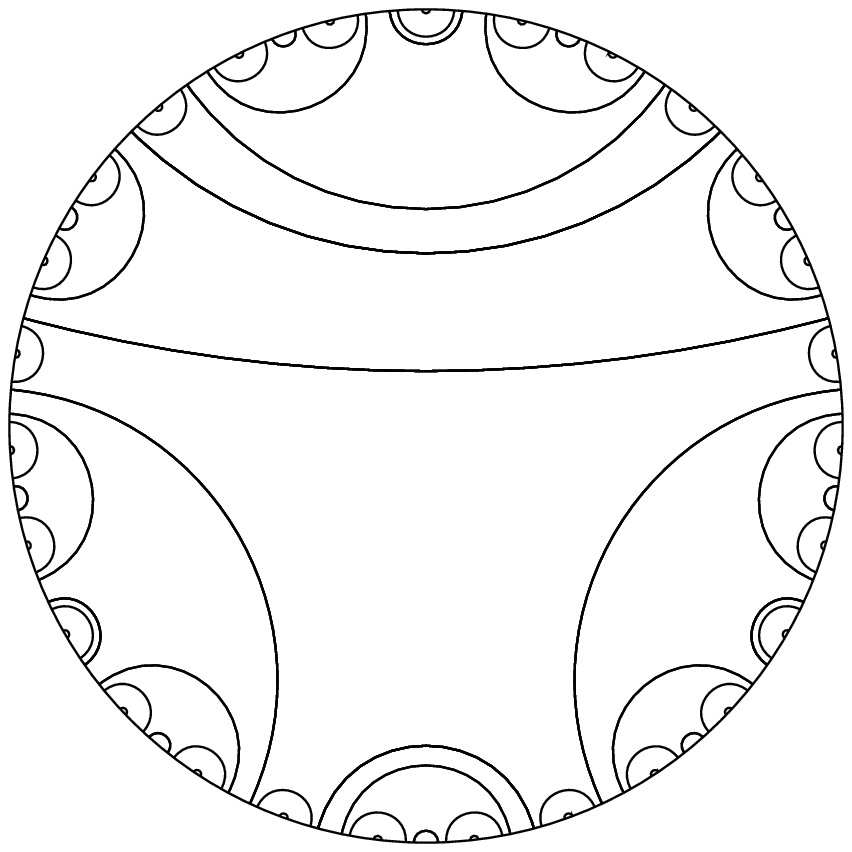}}
  
  \fbox{\includegraphics[width=6cm]{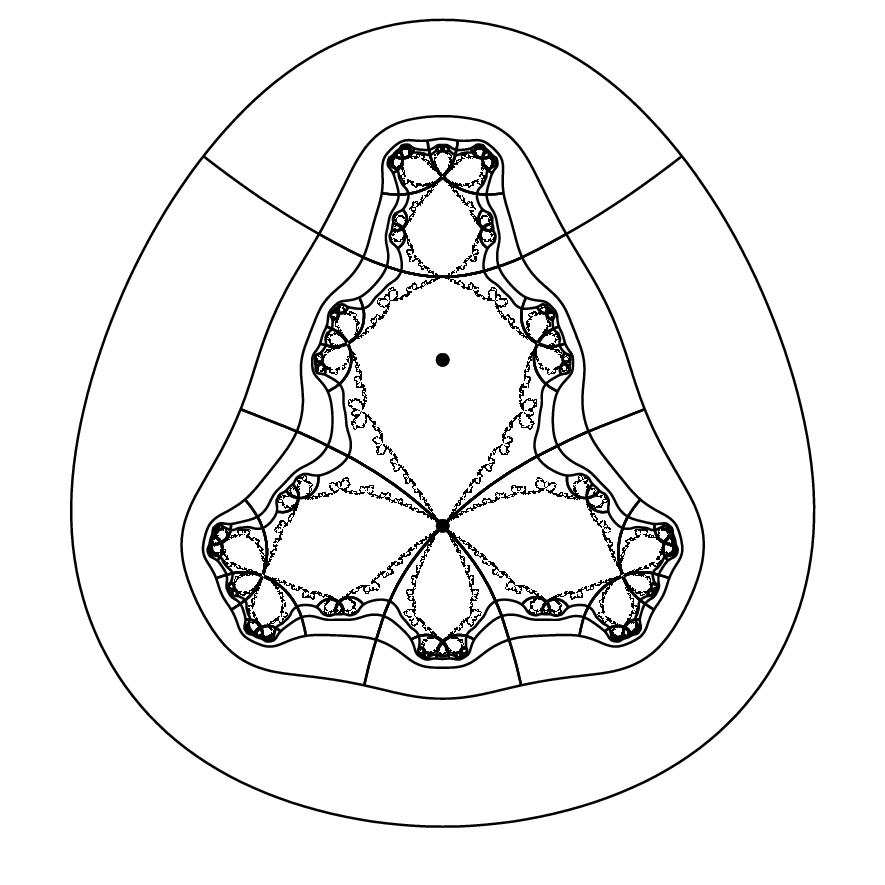}}
  \fbox{\includegraphics[width=6cm]{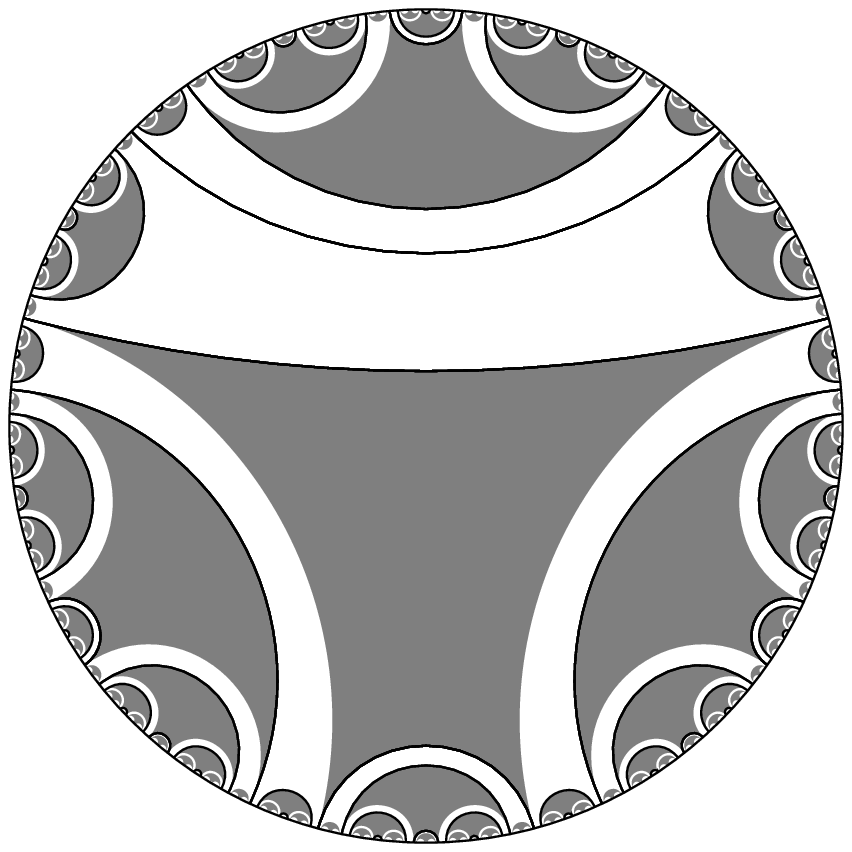}}
  \caption{Yoccoz puzzles and the rational laminations for
   monic centered polynomials affinely conjugate to
  $f_a(z)=az^3-(a+1)z^2+1$ for $a=-1/4$ and 
  $a=\frac{11-3\sqrt{17}}{4}=-0.3423...$
  Dots indicate critical points.}
  \label{fig-capture}
 \end{center}
\end{figure}

\begin{figure}[tb]
 \begin{center}
  \fbox{\includegraphics[width=6cm]{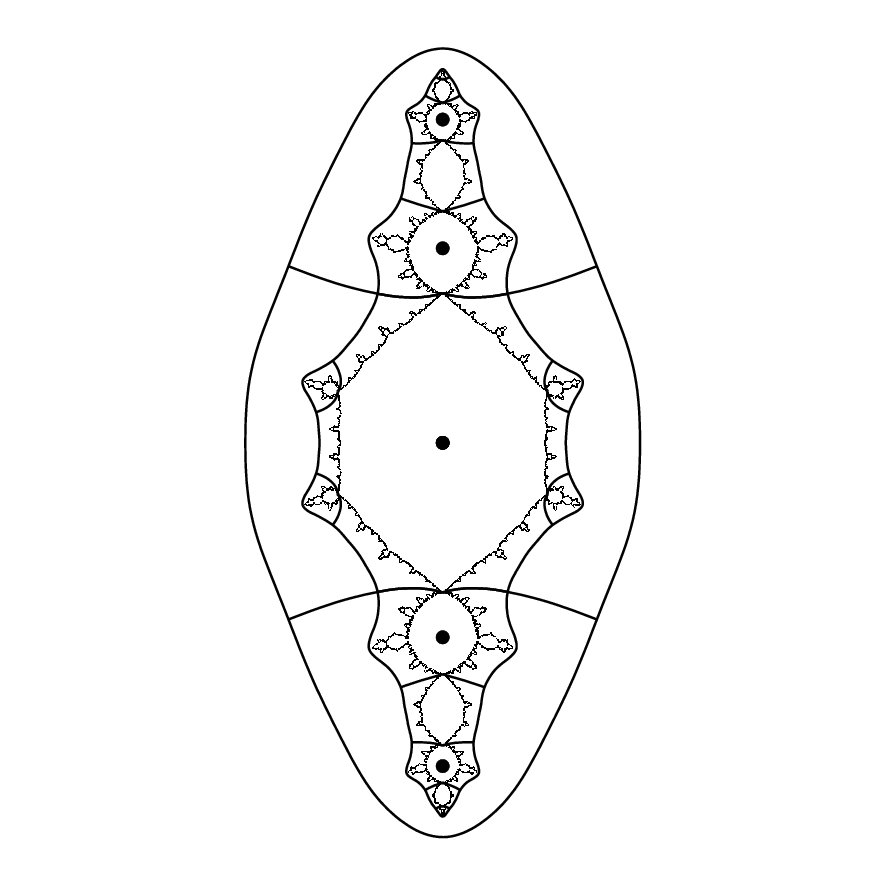}}
  \fbox{\includegraphics[width=6cm]{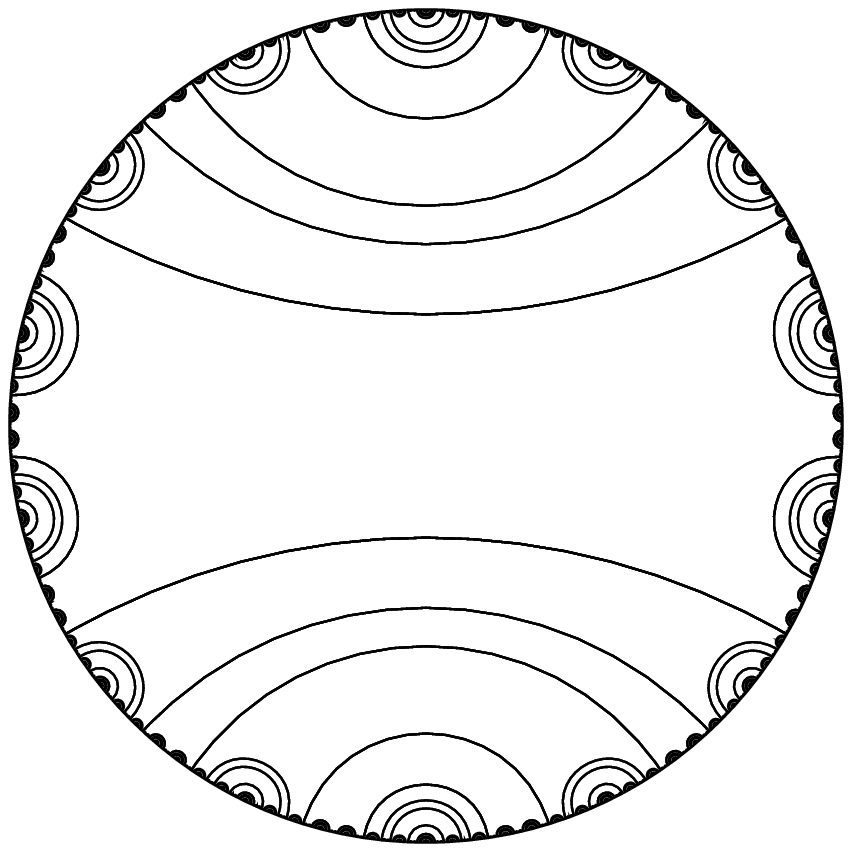}}
  
  \fbox{\includegraphics[width=6cm]{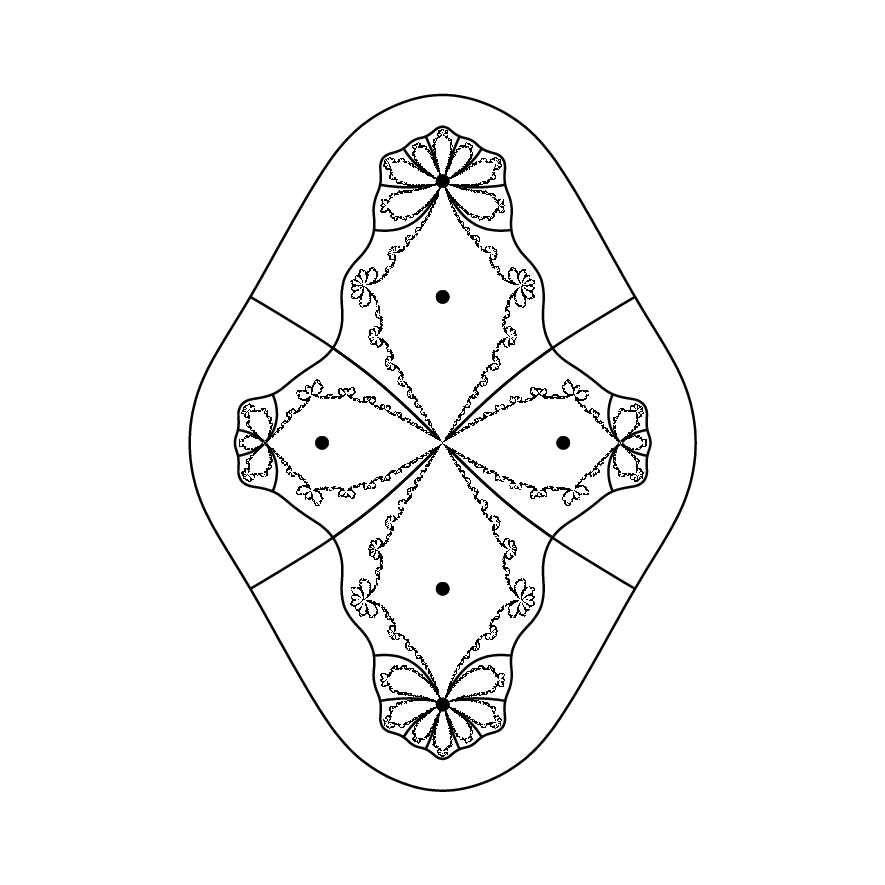}}
  \fbox{\includegraphics[width=6cm]{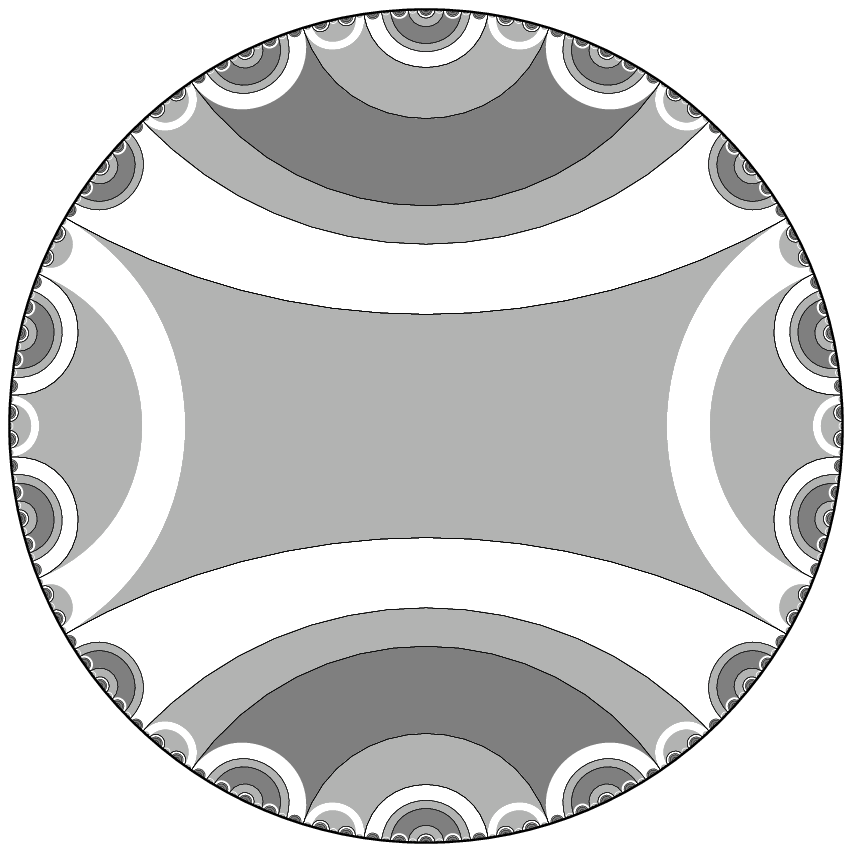}}
  \caption{Yoccoz puzzles and the rational laminations for
  monic centered polynomials affinely conjugate to
  $g_0(z) = 2.959...z^3+2.418...z^5+0.4590...z^7$
  and $g_1(z) = -1.258...z + 0.06411...z^3 +
  0.4047... z^5+0.08182...z^7$. Dots indicate critical points.}
  \label{fig-deg7}
 \end{center}
\end{figure}

The following examples show that the further assumptions in property \ref{item-prop-unlinked
class-deg-f} are necessary.
\begin{exam}
 \label{exam-capture}
 1. Consider the family $P_a(z)=az^3-(a+1)z^2+1$ of cubic polynomials
 having a period two superattractive cycle $\{0,1\}$.
 Let $f_0=P_{a_0}$ and $f_1=P_{a_1}$ where 
 $a_0=-1/4$ and $a_1=\frac{11-3\sqrt{17}}{4}$  (see Figure~\ref{fig-capture}).
 The critical point $c_a=\frac{2(a+1)}{3a}$ is mapped to $0$ for $f_0$
 and for $f_1$, it is mapped to the fixed point which is the
 intersection of the closures of the immediate basins of $0$ and $1$.

 Although $\lanot=\lambda_{f_0}$ is admissible for $f_1$,
 the polynomial $f_1$ is not $\lanot$-renormalizable.
 In fact, let $L$ (resp.\ $L'$) be the $\lanot$-unlinked class whose
 $\lanot$-fiber for $f_0$ contains $0$ (resp.\ $c_{a_0}$).
 Then, we have $K_{f_1}(L) \cap K_{f_1}(L')=\{c_{a_1}\}$.
 Therefore, $K_{f_1}(L)$ contains both of the critical points
and any proper extension of $f_1:K_{f_1}(L) \to K_{f_1}(L')$ to a disk
would have degree $3$. 

 The assumptions of \ref{item-prop-unlinked
 class-deg-f} do not hold. In fact, 
 let $L''$ be the  $\lambda$-unlinked class distinct from $L'$ such that
 $m_d(L)=m_d(L'')$. Note that $\overline{L''} \cap \overline{L'}\ne \emptyset$.
 Also, observe that $A_0=\overline{L} \cap \overline{L'}$ is a $\lanot$-class.
 The $\lambda_{f_1}$-class  
 \[
  A=A_0 \cup (\overline{L''} \cap \overline{L'}),
 \]
formed by the the angles landing at $c_{a_1}$, 
intersects both $\overline{L}$ and $\overline{L''}$.

 \medskip
 2. 
 Let us now illustrate a more complicated situation in which the assumptions of
\ref{item-prop-unlinked class-deg-f} fail to hold.
 Namely, the assumption fails for a $\lanot$-class $A$ intersecting 
the closure of a $\lanot$-unlinked
 class $L$.
However, $L$ is the unique critical  $\lanot$-unlinked class whose closure intersects $A$.
 
 Let $g_0(z) = 2.959...z^3+2.418...z^5+0.4590...z^7$
 and $g_1(z) = -1.258...z + 0.06411...z^3 +
 0.4047... z^5+0.08182...z^7$ (see Figure~\ref{fig-deg7}).
 They are odd polynomials of degree 7.
 The critical points of $g_0$ are
 $0$ (double critical point), $\pm c_1= \pm i$ and $\pm c_2=\pm
 1.662...i$, which satisfy $g_0(0)=0$, $g_0(\pm c_1)=g_0(\mp c_2)=\mp c_1$.
 In particular, $g_0$ is post-critically finite and hyperbolic.
 For $g_1$, the critical points are
 $\pm c_1' = \pm i$, $\pm c_2'= \pm 1.793...i$ and $\pm c_3 = \pm
 0.8265...$ and we have $g_1(\pm c_1)= \mp c_1$, $g_1(\pm c_2)=0$ and
 $g_2(\pm c_3)=\mp c_3$.
 Furthermore, we have $\lambda_{g_1} \supset \lambda_{g_0}$.
 Let $L$ be the $\lambda_{g_0}$-unlinked classes
 such that $c_1 \in K_{g_0}(L)$.
 Then $K_{g_1}(L)$ contains both $c_1'$ and $c_2'$,
 so similar to the previous example,
 this implies the conclusion of \ref{item-prop-unlinked class-deg-f}
 does not hold. 
 
 Observe that there are two ``bifurcations'' from $g_0$ to $g_1$: 
 The first one is that period two points collapse at the origin and 
 new landing relations indicated by the light gray
 regions in the picture of the rational lamination for $g_1$ appear.
 The second is that the critical points $\pm c_2'$ hit the origin,
 the corresponding landing relations are indicated by the dark gray
 regions.
\end{exam}

Below we combinatorially and abstractly describe the ``location'' and ``orbit'' of the critical
points.

  \begin{defn}[Critical elements and orbits]
Let $\lambda$ be a $d$-invariant rational lamination and  
$\hat{\lambda}$ its real extension (Definition~\ref{defn-real extension}).  Given a
$\hat{\lambda}$-class  $A$, denote by $\delta(A)$ the degree of $m_d : A \to
m_d (A)$. If $A$ is a $\hat{\lambda}$-class with $\delta(A) > 1$, then we
say that $A$ is a  {\itshape Julia critical element} of $\lambda$.

If $L$ is an infinite unlinked $\lambda$-class, we say that $\delta(L)$ is the {\itshape degree of $L$} where $\delta(L)$ is
the number obtained in Proposition~\ref{prop-unlinked class}~\ref{item-prop-unlinked class-degree}. 
If $\delta(L) > 1$, then we say that $L$ is a {\itshape Fatou critical element} of $\lambda$.

We denote by
$\crit(\lambda)$ the collection formed by all 
the critical elements for $\lambda$.
The {\itshape post-critical set} and {\itshape critical orbits of $\lambda$} are:
 $$PC(\lambda) = \{m_d^n(C);\ C \in \crit(\lambda), n > 0\},$$ 
 $$CO(\lambda)=\crit(\lambda) \cup
 PC(\lambda).$$
  \end{defn}

It follows that 
$$d= 1+ \sum_{C \in \crit(\lambda)} (\delta(C) -1).$$
(e.g. see~\cite[Lemma~4.10]{Kiwi}).

\begin{defn}[Hyperbolic and post-critically finite laminations]
An invariant rational lamination $\lambda$ is {\itshape hyperbolic} if
it has no Julia critical element. 
An invariant rational lamination $\lambda$ is {\itshape post-critically
  finite} if every
Julia critical element is contained in $\QS$. 
\end{defn}

It follows that a hyperbolic invariant rational lamination is the rational lamination of post-critically finite hyperbolic polynomial.
Also, a post-critically
finite invariant rational lamination is the rational lamination of
a post-critically finite polynomial (see~\cite[Theorem 5.17]{Kiwi}).

\subsection{Renormalizations}
\label{subsec-renormalizations}

 \begin{defn}[Mapping schema associated with $\lanot$]
   Consider an invariant rational lamination $\lanot$. 
\begin{itemize}
 \item 
Let $|\redschema(\lanot)|$ be the collection formed by 
the Fatou critical elements $v$ of $\lanot$.

 \item 
 Let  $\schemamapnot: |\redschema(\lanot)| \to |\redschema(\lanot)|$
 be the map such that $\schemamapnot(v) = v'$ if
 $m_d^{\ell_v} (v)=v'$, for some $\ell_v \ge 0$, and $m_d^k(v) \notin  |\redschema(\lanot)|$, for all $0 < k < \ell_{v}$. 
 We say that  $\ell_{v}$ is the {\itshape return time} of $v$. 
 (Observe that the existence of such $v'$ is guaranteed by
 Proposition~\ref{prop-unlinked class} \ref{item-prop-unlinked
 class-periodic} and \ref{item-prop-unlinked class-critical return}.)

 \item 
 Finally, let
 $\schemadegnot: |\redschema(\lanot)| \to \N$ be the corresponding degree map. That is, $\schemadegnot (v) = \delta (v)$.
\end{itemize}

We call
$$\redschema (\lanot) = (|\redschema(\lanot)|, \schemamapnot, \schemadegnot).$$
the {\itshape reduced mapping schema of $\lanot$}.
 \end{defn}

\begin{defn}[Renormalizable and renormalization]
 \label{defn-renormalization}
We say $f\in \cC(\lanot)$ is {\itshape $\lanot$-renormalizable}
 if, for every $v \in |\redschema(\lanot)|$ 
 there exist topological disks $U'_v$ and $U_v$ such that
 $$g=(f^{\ell_v}:U'_v \to U_{\schemamapnot(v)})_{v\in |\redschema(\lanot)| }$$
is a polynomial-like map
 over $\redschema (\lanot)$ with fiberwise connected filled Julia set
 $$K(g) = \bigcup_{v\in |\redschema(\lanot)| } \{v\} \times K_f(v).$$
 We call $g$ a {\itshape $\lanot$-renormalization} of $f$.
 We denote by $\cR(\lanot)$ the subset of $\cC(\lanot)$ formed by 
all the $\lanot$-renormalizable polynomials.
\end{defn}

The above definition is equivalent, although a priori stronger, than  the one given in the Outline
(Section~\ref{sec-outline}):
\begin{prop}
\label{small-are-small-p}
 Let $f \in \cC(\lanot)$.
 If there exist topological disks $U'_v$ and $U_v$ for each $v \in
 |\redschema(\lanot)|$ such that
 $K_f(v) \subset U'_v$ and
 $g=(f^{\ell_v}:U'_v \to U_{\schemamapnot(v)})_{v \in
 |\redschema(\lanot)|}$ is a polynomial-like map over
 $\redschema(\lanot)$,
 then $g$ is a $\lanot$-renormalization of $f$.
\end{prop}

\begin{proof}
 It suffices to show $K_f(v)=K(g,v)$.
 By construction, we have $K_f(v) \subset K(g,v)$.
 Since $v$ is an infinite $\lanot$-unlinked class, $K_f(v)$ is a
 full continuum. For $v$ periodic of period $k$ under $\sigma_\lanot$,
let $v_0=v, v_1 =\sigma_\lanot(v_0), \dots, v_{k-1}=\sigma^{k-1}_\lanot(v_0)$
and let $\ell_0 = \ell_{v_0}, \dots, \ell_{k-1} = \ell_{v_{k-1}}$. 
Consider
 $U''_v = (f|_{U_{v_0}'}^{-\ell_0} \circ \cdots \circ
 f|_{U'_{v_{k-1}}}^{-\ell_{k-1}})(U_v)$.  
 Then, for $n=\ell_0 + \cdots + \ell_{k-1}$,
we have that $f^n:U''_v \to U_v$ is a polynomial-like map and its filled Julia
 set is $K(g,v)$, by definition.
 It follows that $K(g,v)$ is the smallest full continuum in $U''_v$ completely
 invariant by $f^n:U''_v \to U_v$, since the analogue statement holds for the polynomial map which is the straightening of $f^n$. \change
 Therefore, $K_f(v) \supset K(g,v)$.
 For preperiodic $v \in |\redschema(\lanot)|$, take a
 backward image of $K(g,v')$ where $v'=\sigma^k(v)$ is periodic.
\end{proof}

In order to extract polynomial-like maps with an external marking it is 
convenient to ``internally mark'' the invariant rational lamination $\lanot$.

\begin{defn}[Internally angled]
 An {\itshape internally angled invariant rational lamination} $\lanot$
 is a pair $(\lanot, (\alpha_v:\overline{v} \to \S)_{v \in
 |\redschema (\lanot)|})$ such that
 \begin{enumerate}
  \item $\alpha_v$ induces a homeomorphism from $\overline{v}/\lanot$
	to $\S$.
  \item $\alpha_{\schemamapnot (v)}(d\theta) = \schemadegnot (v)\alpha_v(\theta)$ for all $\theta
	\in \overline{v}$.
 \end{enumerate}
 We call $(\alpha_v)_{v \in |\redschema (\lanot)|}$ an {\itshape internal
 angle system}.
\end{defn}

The existence of an internal angle system is guaranteed by
Proposition~\ref{prop-unlinked class} \ref{item-prop-unlinked class-degree}.
An internally angled invariant rational lamination determines an
external marking of the $\lanot$-renormalization of every $f \in
\cR(\lanot)$ as follows.

\begin{defn}[Induced external marking]
\label{induced-marking-d}
  Let $(\lanot, (\alpha_v:\overline{v} \to \S)_{v \in |\redschema (\lanot)|})$ be an internally angled invariant rational lamination.
  For each $v \in |\redschema (\lanot)|$ choose an argument $\theta_v$ in
$\alpha_v^{-1}(0)$.
  Given a $\lanot$-renormalization $g$ of a polynomial $f \in \cR(\lanot)$,
let $\Gamma_v$ be the access with representative given by the connected component of 
$\overline{R_f (\theta_v) \cap U^\prime_v}$ that   intersects $K_f(v)$. We say that $\Gamma = (\Gamma_v)_{v \in |\redschema (\lanot)|}$ is the {\itshape external marking of $g$ determined by  
the internal angle system  $(\alpha_v:\overline{v} \to \S)_{v \in |\redschema (\lanot)|}$}.
\end{defn}

\begin{rem}
  \label{induced-marking-r}
  The external marking of $g$ determined by an
internal angle system  is independent of the choices involved
in the definition above. In fact, for another choice $\theta_v'$,
$\{\theta_v,\theta_v'\}$ and $v$ are unlinked.
Hence we may assume that the sector $\C \setminus \overline{\sector(\theta_v,\theta_v';v)}$ is
disjoint from $K_f(v)$.
Therefore, $\overline{R_f(\theta_v)}$ and $\overline{R_f(\theta_v')}$ 
define the same access because $\C \setminus \overline{\sector(\theta_v,\theta_v';v)}$ is
simply connected and its boundary is $\overline{R_f(\theta_v)} \cup
\overline{R_f(\theta_v')}$.
Moreover, the external marking is invariant, since $d^{\ell_v} \theta_v$ is an argument $\theta_{\sigma_{\lambda_0} (v)}$
of a ray landing at the same point as the
rays with arguments in $d^{\ell_v} \alpha_v^{-1}(0) = \alpha_{\sigma_{\lambda_0} (v)}^{-1}(0)$, and the access to $K_f({\sigma_{\lambda_0} (v)})$ \change
determined by this ray   only depends on $\alpha_{\sigma_{\lambda_0} (v)}^{-1}(0)$.
 \end{rem}

\begin{lem}
  Let $\lanot$ be an invariant rational lamination.
  Let $g$ be a $\lanot$-renormalization of $f \in \cR(\lanot)$ and $\Gamma$ an external marking of $g$.
Then there exists an internal angle system $\alpha$ such that the external marking determined by $\alpha$
is exactly $\Gamma$. 
\end{lem}

\begin{proof}
  Let $\beta = ( \beta_v )$ be an internal angle system for $\lanot$.
  Consider the straightening of $g$ that maps the external marking induced by $\beta$ onto the standard marking of $P=\chi_\lanot(g)$.
  It follows that $\Gamma$ is mapped onto an external marking of $P$ determined by a collection of external rays $( R_P(\theta_v,v) )$.
  Now the internal angle system  $\alpha = ( \alpha_v = \beta_v + \theta_v )$ is such that the induced external marking $\Gamma_\alpha$ of $g$
maps onto the external marking of $P$ determined by $( R_P(\theta_v,v) )$. Therefore, $\Gamma = \Gamma_\alpha$.
\end{proof}

\subsection{Proof of Proposition~\ref{prop-unlinked class}}
Most of this proposition is already proved in
\cite[Section~4]{Kiwi} and \cite[Section~3]{Inou-lim}.
More precisely, in \cite{Kiwi}: for \ref{item-prop-unlinked class-invariance} see~Lemma~4.6, for \ref{item-prop-unlinked class-finite} see Lemma~4.8  and Proposition 4.3 (proofs of R3 and R4), for \ref{item-prop-unlinked class-complement} and \ref{item-prop-unlinked class-degree} see Lemma~4.22 and Corollary~4.23 and for \ref{item-prop-unlinked class-critical return} see Corollary~4.18.  
In \cite{Inou-lim}: for \ref{item-prop-unlinked class-complement} and \ref{item-prop-unlinked class-degree} see Proposition~3.6.  
It only remains to establish properties~\ref{item-prop-unlinked class-K},
 \ref{item-prop-unlinked class-finite-deg-f} and \ref{item-prop-unlinked
 class-deg-f}.

It is convenient to consider a (countable) set $M$ contained in $\lanot \subset
\QS \times \QS$ formed by the ``image'' of non-trivial relations. More precisely, 
$$M = \{ (d\theta, d\thetap) \in \lanot; \theta \ne \thetap,
(\theta, \thetap) \in \lanot \}.$$
Observe that $M \supset \lambda_0 \setminus \{$diagonal$\}$.

Given a $\lanot$-unlinked class $L$,
 the invariance of $\lanot$ implies that:
 \begin{align*}
  K_f(m_d(L)) &= K(f) \cap \bigcap_{(\theta, \thetap) \in M}
  \overline{\sector(\theta,\thetap; m_d(L))}, \\
  K_f(L) &= K(f) \cap \bigcap_{(\theta, \thetap) \in M} 
   \overline{S'({(\theta, \thetap)})},
 \end{align*}
 where
 \[
  S'{(\theta, \thetap)} = 
 \bigcap_{(\theta_1,\theta_2) \in \lanot,~
  d\theta_1=\theta,~ d\theta_2=\thetap} 
  \sector(\theta_1,\theta_2; L)
 \]
 or equivalently, $S'(\theta,\thetap)$ is the union of all the
 components of $f^{-1}(\sector(\theta, \thetap; m_d(L)))$ each of which
 contains $R_f(\theta'')$ for some $\theta'' \in L$.  

Now \ref{item-prop-unlinked class-K} follows, after countable intersection 
of the equation,
 \[
  f\left(\overline{S'{(\theta, \thetap)}} \right)
 = \overline{\sector(\theta,\thetap; m_d(L)),}
 \]
which holds for all $(\theta, \thetap) \in M$.

To prove \ref{item-prop-unlinked class-finite-deg-f}, we assume that
$L$ is finite and recall that $L \subset (\R \setminus \Q)/\Z$ by definition.
 Hence $K_f(L) \subset S'(\theta,\thetap)$
 for all $(\theta, \thetap) \in M$. 
Since, 
 $f:S'(\theta,\thetap) \to \sector(\theta, \thetap; m_d(L))$ is a proper map
for all $(\theta, \thetap) \in M$, property~\ref{item-prop-unlinked class-finite-deg-f} follows.

 Now assume $L$ is infinite and satisfies the assumption of
 \ref{item-prop-unlinked class-deg-f}. 
 Let $M_L = 
\{(\theta, \thetap) \in M;\ \theta, \thetap \in \overline{m_d(L)}\}$
 and for each $(\theta, \thetap) \in M_L$, let
 \[
  M'_L (\theta, \thetap) = \{(\theta_1, \theta_2) \in M;\ 
d\theta_1=\theta,\ d\theta_2=\thetap,\ \theta_j \in
 \overline{L},\ (j=1,2)\}.
 \]
 Then we have
\begin{align*}
  K_f(m_d(L)) &= K(f) \cap \bigcap_{(\theta, \thetap) \in M_L}
  \overline{\sector(\theta,\thetap; m_d(L))}, \\
  K_f(L) &= K(f) \cap \bigcap_{(\theta, \thetap) \in M_L} 
   \overline{S''(\theta, \thetap)},
\end{align*}
 where
 \[
  S''(\theta, \thetap) = 
 \bigcap_{(\theta_1, \theta_2) \in M'_L(\theta,\thetap)} 
  \sector(\theta_1,\theta_2; L).
 \]
 The assumption guarantees that, for all $(\theta, \thetap) \in M_L$, 
 each component $K$ of $f^{-1}(K_f(m_d(L)))$ is either
 contained in or disjoint from $\overline{S''(\theta, \thetap)}$.
Hence, for all  $(\theta, \thetap) \in M_L$, $$f:f^{-1}(K_f(m_d(L))) \cap \overline{S''(\theta, \thetap)} \to K_f(m_d(L))$$
 has a proper holomorphic extension of degree at least $\delta(L)$.
 After taking an intersection to exclude all other critical elements,
 property \ref{item-prop-unlinked class-deg-f} follows.
 \qed


\section{Yoccoz puzzles}
\label{subsec-puzzles}

First, we introduce combinatorial Yoccoz puzzles in terms of
laminations, then we define Yoccoz puzzles corresponding to a combinatorial
Yoccoz puzzle. See~\cite{Inou-lim} for more details.

Throughout this section we fix an integer $d \ge 2$.

\subsection{Combinatorial Yoccoz Puzzle}

\begin{defn}
 A {\itshape combinatorial Yoccoz puzzle $\Lambda=(\Lambda_k)_{k \ge 0}$} is a sequence of rational laminations with finite support such that:
 \begin{enumerate}
  \item If $\theta \in \supp(\Lambda_0)$, then $\theta$ is periodic under $m_d$.
  \item If $A$ is a $\Lambda_0$-class, then
	$m_d(A)$ is a $\Lambda_0$-class and $m_d:A \to m_d(A)$ is
	consecutive preserving.
  \item $(m_d)_*\Lambda_{k+1} = \Lambda_k$ for all $k\ge 0$.
  \item If $A$ is a non-trivial $\Lambda_{k+1}$-class, then $m_d(A)$ is also
	a non-trivial $\Lambda_k$-class.
 \end{enumerate}
We say that $\supp(\Lambda)=\bigcup_k \supp(\Lambda_k)$ is the {\itshape support of $\Lambda$}.

A $\Lambda_k$-unlinked class (see~Definition~\ref{d unlinked class}) is a {\itshape
 combinatorial (puzzle) piece of depth $k$}.
The collection of all depth $k$ combinatorial puzzle pieces is denoted by  $\cL_k(\Lambda)$.
 
 We say that a combinatorial Yoccoz puzzle $\Lambda$ is {\itshape admissible for a lamination $\lambda$}
 if $\Lambda_k \subset \lambda$ for all $k \ge 0$.
 We say $\Lambda$ is {\itshape admissible for a polynomial $f \in
 \cC (d)$}
 if $\Lambda$ is admissible for $\lambda_f$.
\end{defn}

We state, without proof, the following rather straightforward properties of combinatorial puzzles.

\begin{prop}
 \label{prop-comb-puzzles}
 Let $\Lambda = (\Lambda_k)_k$ be a combinatorial Yoccoz puzzle.
 Then 
 \begin{enumerate}
  \item $\{\supp(\Lambda_k)\}_k$ forms an increasing sequence in
	$\Q/\Z$ and $\supp(\Lambda)$ is completely invariant under $m_d$.
  \item For all $k \ge 0$, if  $L$ is a combinatorial piece of depth $k+1$, then
	there exists a unique  combinatorial piece $L'$ of depth $k$ such that $L
	\subset L'$.
  \item 
    \label{prop-comb-puzzles-degree}
    For all $k \ge 0$, if $L'$ is a  combinatorial piece of  depth $k+1$, then
	$m_d(L')$ is a combinatorial piece of depth $k$ and $m_d:L' \to m_d(L')$ is $\delta$-to-one for some $\delta \ge 1$.
  \item Let $L$ be a combinatorial piece of depth $k$, for some $k \ge 0$. If $\theta \in \partial L$, then
	for all $k'>k$, there exists a unique  combinatorial piece $L'$ of depth $k'$ such
	that $L' \subset L$ and $\theta \in \partial L'$.
  \item \label{prop-comb-puzzles-total-degree}
        Let $L_1, \dots, L_n$ be a complete list (without repetitions) of the combinatorial puzzle pieces of depth $k$. 
Let $C_1, \dots, C_m$ be a complete list (without repetitions) of the $\Lambda_k$-critical classes. Then,
$$\sum_{i=1}^n (\delta(L_i)-1) + \sum_{j=1}^m (\delta(C_j)-1) = d-1,$$
where, for $A=L_i$ or $C_j$, the number $\delta (A)$ denotes the degree of $m_d : A \rightarrow m_d (A)$. 
 \end{enumerate}
\end{prop}

For a combinatorial puzzle piece $L$, we denote by $\delta(L)$ the degree of $m_d : L \rightarrow m_d(L)$.

\subsection{Yoccoz Puzzle}

Let $f$ be a monic centered polynomial with connected Julia set and associated B\"ottcher map $\varphi: \C \setminus K(f) \to \CDC$ where
$\Delta = \{ |z| < 1 \}$.
For $r > 0$, the Jordan curve
$$E_f (r) = \varphi^{-1}(\{ z \in \CDC ; |z| = \exp (r) \})$$
is called the {\itshape equipotential of level $r$}.
The open topological disk bounded by $E_f (r)$ will be denoted by $D_f (r)$.

Fix $r>0$.
Let $\Lambda=(\Lambda_k)$ be an admissible combinatorial
Yoccoz puzzle for $f \in \cC (d)$.
We define the {\itshape depth $k$ Yoccoz puzzle for $(f,\Lambda)$} as follows.
Given a combinatorial puzzle piece $L$ of depth $k$, let
\[
  P_k(f,L) = \overline{D_f(r/d^k)} \cap  \bigcap_{\theta \sim_{\Lambda_k}\thetap, \theta \neq \thetap} \overline{\sector(\theta,\thetap;L)}
\]
be the {\itshape puzzle piece of depth $k$ for $(f,\Lambda)$ associated with $L$}.
We omit $k$ and/or $f$ when no confusion arises.
Recall that we denote the collection of combinatorial puzzle pieces of depth $k \ge 0$ 
by $\cL_k(\Lambda)$.
Similarly we let
 $$\cP_k(f,\Lambda)= \{P(L);\ L \in \cL_k(\Lambda)\}$$
be the collection of all puzzle pieces of depth $k$. That is, the {\itshape Yoccoz puzzle of depth $k$ for $(f,\Lambda)$}.

\begin{rem}
 The interior of a puzzle piece is not connected in general.
 A ``degenerate'' puzzle appears
 when two distinct $\Lambda_k$-classes $A, A'$ are contained in the same
 $\lambda$-class.
 Therefore, our definition of  puzzle piece 
does not coincide with the usual one (i.e., the closure of a
 bounded component of the complement of a suitable graph constructed 
with equipotentials and external rays).
\end{rem}

 \begin{figure}[bt]
  \begin{center}
   \fbox{\includegraphics[width=6cm]{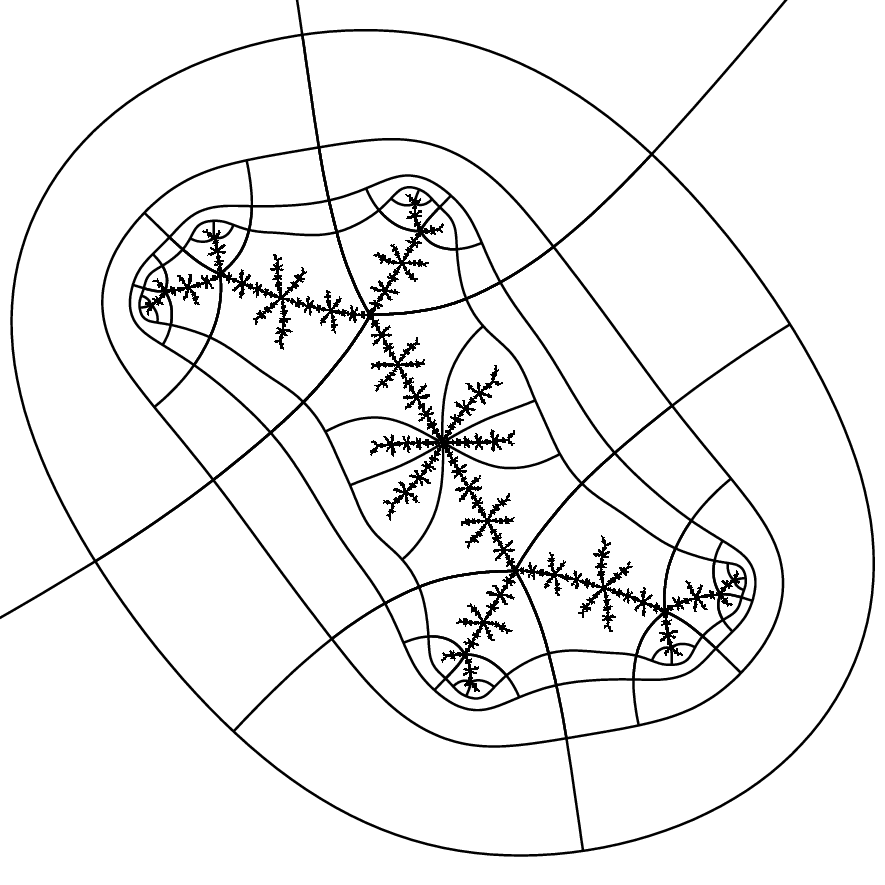}}
   \fbox{\includegraphics[width=6cm]{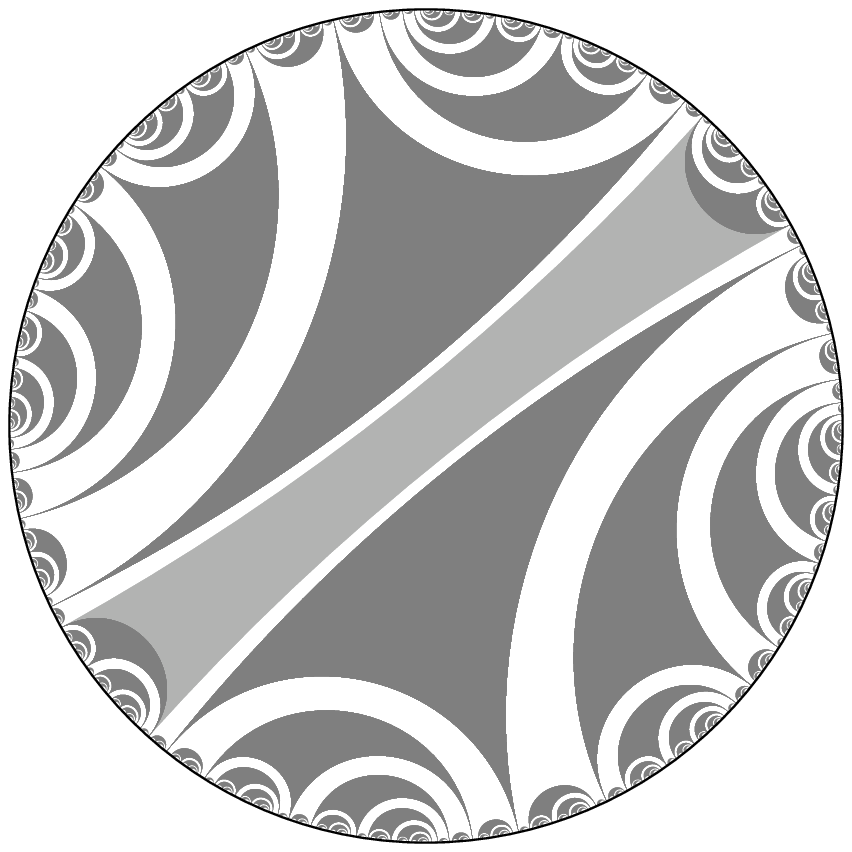}}
   \caption{The Julia set and a ``degenerate'' puzzle for $z^2+c$ with
   $c=-0.1010...+0.9562...\imunit$.}
   \label{fig-degpuzzle}
  \end{center}
 \end{figure}  

 \begin{figure}[tb]
  \begin{center}
   \fbox{\includegraphics[width=6cm]{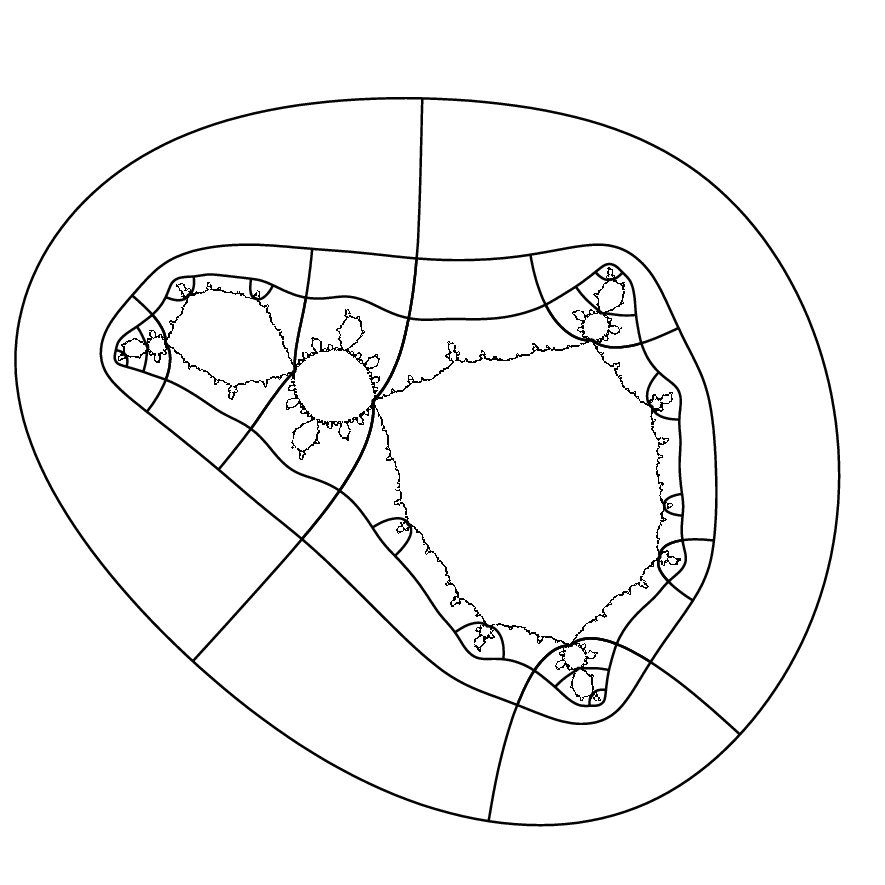}}
   \fbox{\includegraphics[width=6cm]{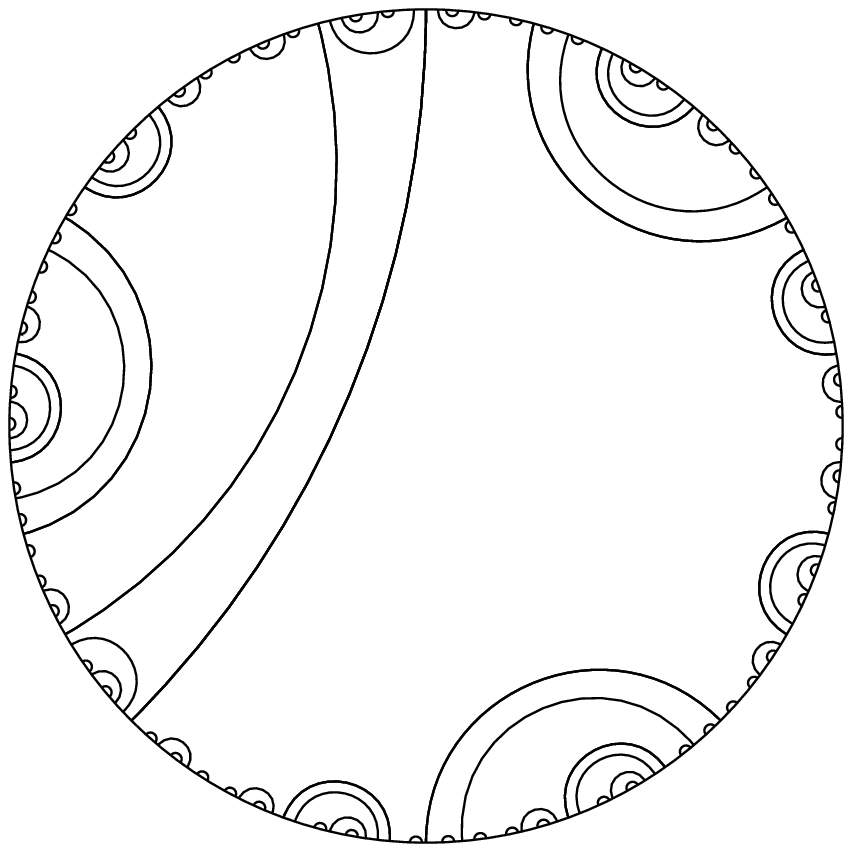}}
   
   \fbox{\includegraphics[width=6cm]{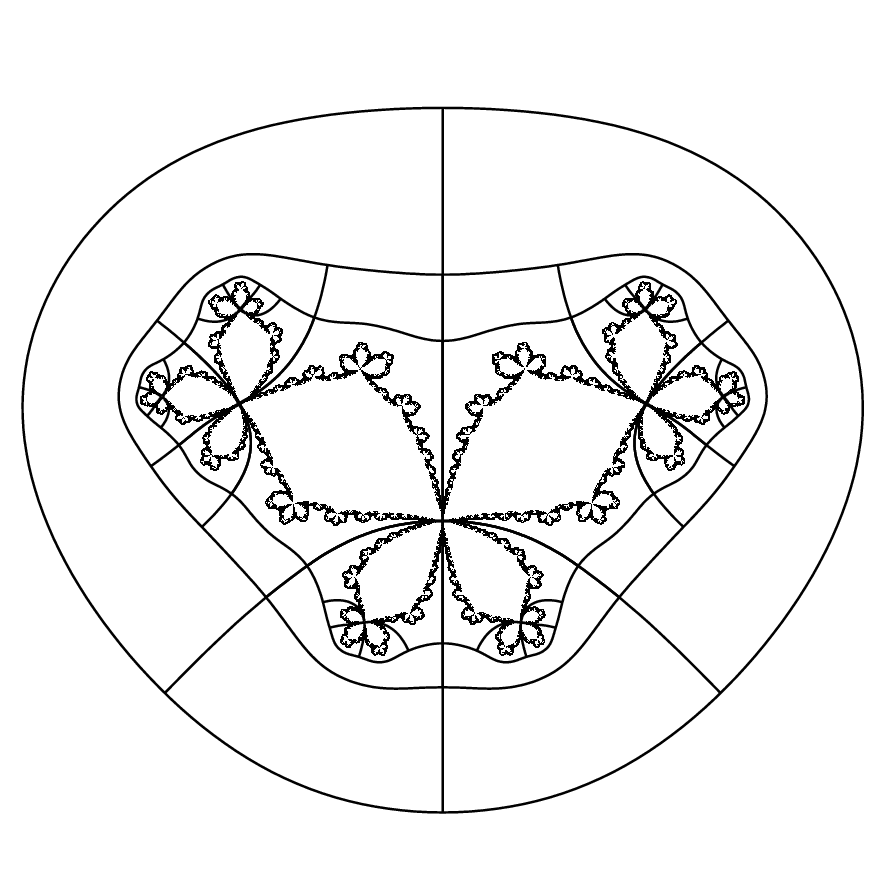}}
   \fbox{\includegraphics[width=6cm]{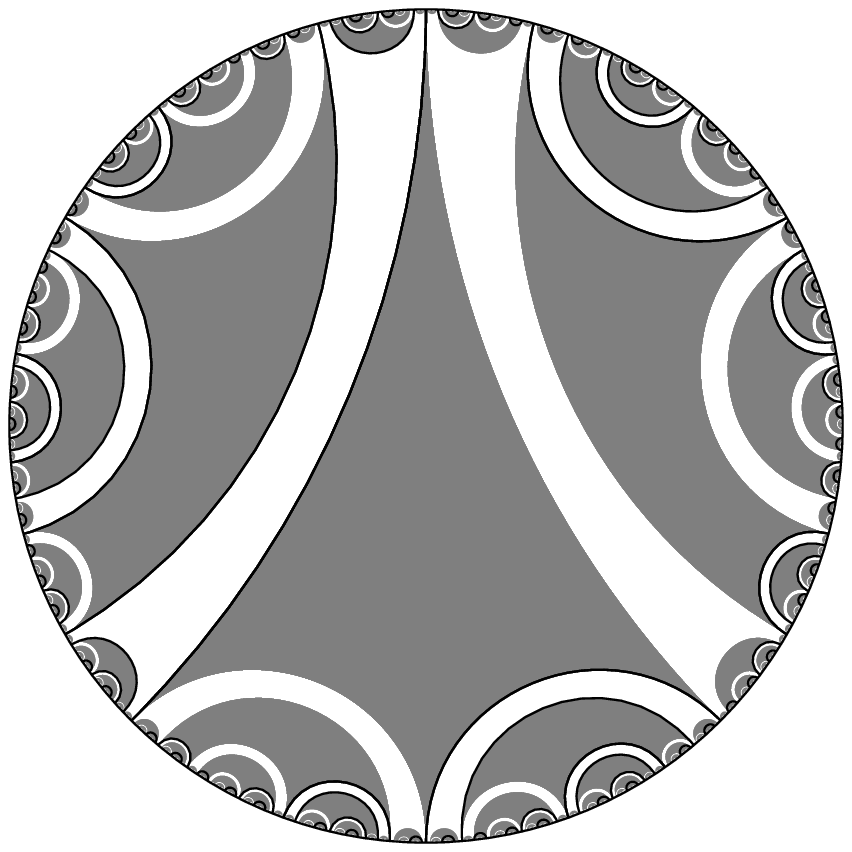}}
   \caption{Yoccoz puzzles and the rational laminations for
   $f_0(z)=z^3+az^2$ ($a=1.502...-0.7790...\imunit$) and
   $f_1(z)=z^3-3/4z-\sqrt{7}\imunit/4$.}
   \label{fig-cubic2+2}
  \end{center}
 \end{figure}

\begin{exam}
 1. 
 Let $c=-0.1010...+0.9562...\imunit$ and $f(z)=z^2+c$,
 so that $\alpha=f^4(c)$ is the alpha fixed point (i.e., the fixed point
 which is not a landing point of $R_f(0)$).
 Figure~\ref{fig-degpuzzle} shows its Julia set with its Yoccoz puzzles
 and (a part of) its lamination.
 The rational lamination $\lambda_f$ of $f$ properly contains 
the rational lamination of 
 Douady's rabbit (say $\lanot$).
 In fact, $f^4(c)$ has the landing angles 
 of the alpha fixed point, i.e., $1/7,2/7,4/7$. Thus the critical point $0$ has
 six landing angles, consisting of two $\lanot$-equivalence classes.
 
 The light gray region in the right picture indicates the landing
 relation for the critical point $0$,
 which is not contained in the lamination of Douady's rabbit.
 The Yoccoz puzzles for $f$, determined by the 
 lamination of Douady's rabbit (namely, let $E_k$ be the set of
 landing angles for $f^{-k}(\alpha)$ and $\Lambda_k=\lambda_0|_{E_k}$),
 is such that the interior of the puzzle piece of depth 4 containing $\{0,\alpha,-\alpha\}$ has two \change
 components.

 \medskip 
 2.
 Let $f_0(z)=z^3+az^2$ where $a=1.502...-0.7790...\imunit$
 and let $f_1(z)=z^3-3/4z-\sqrt{7}\imunit/4$ 
 (Figure~\ref{fig-cubic2+2}).
 There are two superattractive cycles of periods one and two for
 $f_0$, and $f_1$ has two superattractive cycles both of period two.
 Observe that $\lambda_{f_0} \subset \lambda_{f_1}$. Every non-trivial 
class of $\lambda_{f_0}$ eventually maps onto the class formed 
by the angles $1/4, 5/8$.
 As above, we consider the  Yoccoz puzzles $(\Lambda_k)$ determined
by $\lambda_{f_0}$. The puzzle piece for $f_1$ of any depth
 corresponding to the  fixed superattractive basin for $f_0$ has
 disconnected interior (and the number of components increases as
 depth increases).
\end{exam}

\begin{prop}
 \label{prop-puzzles}
 Let $\Lambda = (\Lambda_k)$ be an admissible combinatorial Yoccoz
 puzzle for $f \in \cC (d)$.
For all $k \ge 0$,  let 
 \[
  \Gamma_k = E_f(r/d^k) \cup \bigcup_{\theta \in \supp(\Lambda_k)}
 \overline{R_f(\theta)}.
 \]
 Then the following statements hold:
 \begin{enumerate}
  \item If $L$ is a combinatorial puzzle piece, then $P(L)$ is compact, connected and full. Moreover,
the interior of $P(L)$ is the union of all the bounded components $W$ of 
	$\C \setminus \Gamma_k$
	for which there exists $\theta \in L$ such that  $R_f(\theta)$ has non-empty intersection with $W$.
  \item $\cP_k(f,\Lambda)$ is a partition of
	$\overline{D_f(r/d^k)}$, i.e., depth $k$ puzzle pieces have mutually disjoint interior and 
	\[
	 \overline{D_f(r/d^k)} 
	= \bigcup_{P \in \cP_k(f,\Lambda)} P.
	\]
  \item If combinatorial puzzle pieces $L$ and $L'$ satisfy $L \subset L'$, then
	$P(L) \subset P(L')$.
  \item If $L$ is a combinatorial puzzle piece of depth $k \ge 1$, then $f(P(L)) = P(m_d(L))$ and $f:\Int P(L) \to \Int P(m_d(L))$
	is a proper map of degree $\delta(L)$.
 \end{enumerate}
\end{prop}
For the proof, see \cite[Proposition~4.1]{Inou-lim}.
This proposition implies that our puzzles have similar properties as
the usual Yoccoz puzzles (defined by the closures of bounded components
of $\C \setminus \Gamma_k$) and moreover, their properties can be
described only in terms of $\Lambda$.

\subsection{From combinatorial puzzles to laminations}

Combinatorial Yoccoz puzzles are particularly suited to 
study the smallest invariant rational lamination for which
they are admissible.

\begin{defn}
  Let $\Lambda = (\Lambda_k)$ be a combinatorial Yoccoz puzzle.
  The {\itshape rational lamination $\lambda(\Lambda)$ generated by $\Lambda$} 
is the smallest equivalence relation in $\QS$ containing the closure of
$\bigcup_k \Lambda_k$ in $\QS \times \QS$.
\end{defn}

The previous definition is justified with the following result.

\begin{prop}
 \label{prop-comb-puzzle-to-lamination}
 For a combinatorial Yoccoz puzzle $\Lambda=(\Lambda_k)_{k \ge 0}$,
 the rational lamination $\lambda(\Lambda)$ generated by $\Lambda$ is 
an invariant rational lamination.
\end{prop}

\begin{proof}
 First, it is clear that $\lambda (\Lambda)$-classes are pairwise unlinked.
 Let $\lambda = \bigcup_k \Lambda_k$.
 Applying the argument for the claim in the proof of
 \cite[Lemma~4.7]{Kiwi} to our case, we have that the following assertion holds:

  If $\theta,\theta' \in \Q/\Z$ are $\overline{\lambda}$-related but not
  $\lambda$-equivalent, and if $\theta$ is periodic under  $m_d$,  
  then $\theta'$ is periodic of the same period.

 If  $A$ is a $\lambda (\Lambda)$-class, then  there exists an iterate $ \ell \geq 0$ such that $m_d^\ell(A)$ contains a periodic angle under iterations of  $m_d$.
 Therefore, $m_d^\ell(A)$ (hence $A$) is a finite set since, by the assertion above,  $m_d^\ell(A)$ consists of periodic angles of
 the same period.

 Now let us prove that $\lambda (\Lambda)$ is closed.
 Assume that, for all $n \ge 0$,  $\theta_n,\theta_n'$ are $\lambda (\Lambda)$-equivalent arguments such that
 $\theta_n \to \theta$ and $\theta_n' \to \theta'$, as $n \to \infty$.
 Then, since $\lambda (\Lambda)$-classes are pairwise unlinked,
 we may assume that these sequences are monotone. That is,  $\theta_n \searrow \theta$ and $\theta_n' \nearrow
 \theta'$. Moreover, we may also assume that,  for all $n$, the  $\lambda (\Lambda)$-class containing $\theta_n$ is contained
 in $[\theta_n,\theta_{n-1}) \cup (\theta_{n-1}',\theta_n']$.
It follows  that there exists some $\tilde{\theta}_n \in
 [\theta_n,\theta_{n-1})$ and $\tilde{\theta}_n' \in
 (\theta_{n-1},\theta_n]$ such that $\tilde{\theta}_n,
 \tilde{\theta}_n'$ are $\lambda$-equivalent.
 Hence, $\theta, \theta'$ are $\overline{\lambda}$-related and, therefore, $\lambda (\Lambda)$-equivalent.

 It remains to prove that $(m_d)_*\lambda (\Lambda)=\lambda (\Lambda)$.
The image $m_d(A)$ of a $\lambda (\Lambda)$-class $A$ is clearly contained
 in a $\lambda (\Lambda)$-class $B$.
Consider two elements $\theta_1, \theta_2$ of $A$ such that $(\theta_1,\theta_2)$ is a connected component of
$\R/\Z \setminus A$.  
We suppose that $\theta_0' \in B \cap (d\theta_1,d\theta_2)$
and proceed to obtain a contradiction in order to prove that $B \cap (d\theta_1,d\theta_2) = \emptyset$.
 We may assume $\theta_0'$ and $d\theta_1$ are
 $\overline{\lambda}$-related (because either $d\theta_1
 \sim_{\overline{\lambda}} \min[(d\theta_1,d\theta_2)\cap B]$
 or $d\theta_2 \sim_{\overline{\lambda}} \max[(d\theta_1,d\theta_2) \cap
 B]$ holds). 
 Then there exists $\theta_{0,n}' \nearrow \theta_0'$  and $\theta_{1,n}'
 \searrow d\theta_1$ such that $\theta_{0,n}'$ and $\theta_{1,n}'$ are 
 $\lambda$-equivalent.
 Let $\theta_{1,n} \searrow \theta_1$ be such that $d\theta_{1,n} =
 \theta_{1,n}'$.
 Taking preimages we obtain $\theta_{0,n}$ such that $d\theta_{0,n}
 = \theta_{0,n}'$ and $\theta_{0,n}, \theta_{1,n}$ are
 $\lambda$-equivalent.
By the pairwise unlinked property of $\lambda (\Lambda)$-classes, we have
$\theta_{0,n}, \theta_{1,n} \in (\theta_1,\theta_2)$.
 Passing to a subsequence, we may assume $\theta_{0,n} \to \theta_0$
 and $d\theta_0 = \theta_0'$.
 Since  $\theta_0 \in [\theta_1,\theta_2]$ and $\theta_0
 \sim_{\overline{\lambda}} \theta_1$, either $\theta_0 = \theta_1$
 or $\theta_0 = \theta_2$, which is a contradiction.
 Hence $m_d:A \to m_d(A)$ is
 consecutive preserving and $B=m_d(A)$ is a $\lambda (\Lambda)$-class 
(i.e. $(m_d)_* \lambda (\Lambda) = \lambda (\Lambda)$).
\end{proof}

\begin{defn}
\label{def-generator}
 Let $\Lambda=(\Lambda_k)$ be an admissible combinatorial Yoccoz
 puzzle for a $d$-invariant rational lamination $\lambda$.  We say
 $\Lambda$ is a {\itshape generator of $\lambda$} if there exists $k \ge
 0$ such that, for any combinatorial piece $L$ of depth $k$ which contains
 a $\lambda$-unlinked class $v$ in the forward orbit of a Fatou critical element of $\lambda$, we have that  $\delta(L) = \delta(v)$.
  We call such $k$ a {\itshape separation depth}
 for $\lambda$.
\end{defn}

The next proposition shows that a generator actually ``generates'' a
substantial part of the corresponding rational lamination.

\begin{prop}
 \label{prop-generator}
 Let $\Lambda=(\Lambda_k)_{k \ge 0}$ be a generator of an invariant rational lamination $\lambda$.
If $L$ is a Fatou critical element of $\lambda$, then $L$ is a Fatou critical element of $\lambda (\Lambda)$.
\end{prop}

\begin{proof}
 Let $\lambda'={\lambda}(\Lambda)$ and observe that $\lambda' \subset \lambda$.
Therefore,  given any infinite $\lambda$-unlinked class $L$, there exists a unique $\lambda'$-unlinked class $L'$
such that $L' \supset L$. 

Let $L$ be a Fatou critical element of $\lambda$ and $L'$ the
$\lambda'$-unlinked class containing $L$. 
We claim that $L' = L$. For this purpose, let $k$ be a separation depth
 for $\Lambda$. 
Denote by $L_k$ the combinatorial puzzle piece of depth $k$ containing $L'$.
Thus, $\delta(L_k) = \delta(L') = \delta(L)$.
Now assume that $L$ is periodic under $m_d$, say of period
 $p$. Repeating this argument, along the $m_d$-orbit of $L \subset L'$, it follows that $m_d^p:L' \rightarrow L'$ and 
$m_d^p:L \rightarrow L $ have the same degree. Moreover, under both of these maps, the grand orbits of a given $\theta \in L$  coincide.
By Proposition~\ref{prop-unlinked class} \ref{item-prop-unlinked
 class-degree},
 such a grand orbit is dense in $L$ and in $L'$. Therefore, we conclude that $L=L'$. 
In the case that $L$ is preperiodic, then there exists $\ell \geq 1$ such that $m_d^\ell (L) = m_d^\ell (L')$ is periodic 
and the degree of $m_d^\ell:L \rightarrow m_d^\ell (L)$ also coincides with the degree of $m_d^\ell:L' \rightarrow m_d^\ell (L')$.
Hence, $L=L'$.
\end{proof}

\begin{lem}
 \label{lem-strict-inclusion}
 Consider rational laminations $\lambda \subset \lambda'$.
 If $\lambda \ne \lambda'$ ($\lambda'$ is strictly stronger than
 $\lambda$),
 then
 there exist an infinite $\lambda$-unlinked class $L$ which is not a
 $\lambda'$-unlinked class.

 In particular, if every infinite $\lambda$-unlinked class is a
 $\lambda'$-unlinked class, then $\lambda=\lambda'$.
\end{lem}

\begin{proof}
Assume that $\lambda \ne \lambda'$.
Let $A$ and $B$ be distinct $\lambda$-classes which are contained in the same
$\lambda'$-class. Denote by $I$ (resp. $J$) the connected component
of $\R/\Z \setminus A$ (resp. $\R/\Z \setminus B$) containing $B$ (resp. $A$).
Since every $\lambda'$-class is the union of finitely many $\lambda$-classes, we
may assume that $I \cap J$ is free of elements $\lambda'$-equivalent to elements of
$A$ (and to elements of $B$). 
Observe that $I \cap J$ has two connected component $V,W$ and any pair of arguments, one from
each component, is linked with $\partial (I \cap J)$.
 
Now we identify $\R/\Z$ with the unit circle $\partial \Delta$ in
 the complex plane by the map $e(\theta)=e^{2\pi i \theta}$. \changetwo
Denote by $C$ ($\subset \Delta$) the interior of the convex hull in $\overline{\Delta}$
 of $\partial (I\cap J)$ ($\subset A \cup B$) with respect to the hyperbolic metric.
It follows that $C$ is pairwise disjoint with the convex hull of every
$\hat\lambda$-class where $\hat\lambda$ is the real extension of $\lambda$. 
Denote by $\cL$ the union of the convex hull in $\overline{\Delta}$
of all $\hat\lambda$-classes. Hence,
$C$ is contained in a connected component
$U$ of $\Delta \setminus \cL$. By~\cite[Lemmas~4.17~and~4.22]{Kiwi}, $\partial U \cap \partial \Delta$ is a Cantor set
and the endpoints of each connected component of $\R/\Z \setminus \partial U$ are
$\lambda$-equivalent.
Thus $(\partial U \cap \partial \Delta) \setminus \Q/\Z$ is contained in an infinite
$\hat\lambda$-unlinked class which in turn is contained in an infinite
$\lambda$-unlinked class $L$.
Moreover, the extreme points
of $V$ (resp. $W$) are not isolated points
of  $ \partial U \cap \overline{V}$ (resp. $ \partial U \cap \overline{W}$).
Thus, there exist a pair of elements of $L$, one from $V$ and another one from $W$,
 which is linked with $\partial (I \cap J) \subset A \cup B$. 
Hence $L$ is not a $\lambda'$-unlinked class.
\end{proof}

\begin{cor}
 \label{cor-generator}
 Let $\Lambda=(\Lambda_k)_{k \ge 0}$ be a generator of a hyperbolic invariant
 rational lamination $\lambda$. Then $\lambda(\Lambda)=\lambda$.
\end{cor}

\begin{proof}
 It suffices to show that every infinite $\lambda$-unlinked class is a
 $\lambda(\Lambda)$-unlinked class.
 Let $L$ be an infinite  $\lambda$-unlinked class and $L'$ be the
 $\lambda(\Lambda)$-unlinked class containing $L$.
 Since $L'$ is eventually periodic and any periodic orbit of unlinked
 classes contains a critical element, 
 $m_d^k(L')$ is critical for some $k\ge 0$.
 Let $n$ be the smallest such $k$.
 Then $L'$ contains a critical element for $\lambda$.
 Since $\lambda$ is hyperbolic, such a critical element must be a 
 $\lambda$-unlinked class. Hence it follows by the above proposition
 that $m_d^n(L')$ is a critical $\lambda$-unlinked class.
 In particular, $m_d^n(L')=m_d^n(L)$.
 Since $m_d^n:L' \to m_d^n(L')$ is one-to-one,
 we have $L=L'$.
\end{proof}

\begin{lem}
\label{lem-exists generator}
If $\lambda$ is an invariant rational lamination, then there exists a generator of $\lambda$.
\end{lem}

\begin{proof}
Let $v$ be a $\lambda$-unlinked class in the forward orbit of a Fatou critical element of $\lambda$.
By definition, for every critical element $C \neq v$ of $\lambda$ there exists
a $\lambda$-class $A_C$ such that the connected component of $\S \setminus A_C$ that contains $v$ is disjoint from $C$.
Hence, there exists a collection $A_1, \dots, A_n$ of $\lambda$-classes such that the following holds.
If $v$ is a $\lambda$-unlinked class in the forward orbit of a Fatou critical element of $\lambda$ and $C \neq v$  is a critical element of
$\lambda$, then the connected component of $ \S \setminus A_j$ that contains $v$ is disjoint from $C$, for some $j=1, \dots, n$.
Let $E_0$ be the periodic arguments in the $m_d$-forward orbit of $A_1 \cup \cdots \cup A_n$.
For $k \ge 0 $, let $E_k = m_d^{-k}(E_0)$ and $\Lambda_k$ be the restriction of $\lambda$ to $E_k \times E_k$. 

By construction, we may consider a sufficiently large $k$ such that  $A_1,\dots,A_n$ are
$\Lambda_k$-classes.
Given  $v$ in the forward orbit of a Fatou critical element, let  $L \in \cL_k(\Lambda)$ be the combinatorial puzzle piece containing
$v$. If  $L$ contains a critical element $C$, then $v$ must be critical and $C=v$.
Thus,
$$\delta(v) -1 = \sum_{C \subset L} (\delta(C)-1)$$
where the sum is taken over all critical elements $C$ contained in $L$ and
is equal to $0$ if there is no such critical element in $L$.
Moreover, from the proof of~\cite[Lemma~3.9]{Inou-lim}, we have that
\begin{equation}
 \label{eqn-counting-degree}
 \sum_{C \subset L} (\delta(C)-1) = \delta(L)-1.
\end{equation}
It follows that $\delta(v) = \delta(L)$ and that  the puzzle $\Lambda = (\Lambda_k)_{k \ge 0}$
is a generator for $\lambda$.
\end{proof}


\section{Combinatorial surjectivity of straightening}
\label{sec-tuning}
In this section we study the action of 
straightening on rational laminations and establish that, in a certain sense (see~Theorem~\ref{thm-combinatorial onto}), straightening maps
are (almost) combinatorially onto.

\begin{defn}
 A {\itshape rational lamination over a mapping schema
 $T=(|T|,\sigma,\delta)$} is a collection of rational laminations
 $\lambda=(\lambda_v)_{v \in |T|}$ with 
 $\delta(v)_*\lambda_v=\lambda_{\sigma(v)}$.

 For $f \in \cC(T)$, the {\itshape rational lamination}
  of $f$ is the collection $\lambda_f=(\lambda_{f,v})_{v \in |T|}$ of equivalence relations in 
$\QS$ such that:
 $\theta$ and $\theta'$ are $\lambda_{f,v}$-equivalent if
 $R_f(v,\theta)$ and $R_f(v,\theta')$ land at the same point.
\end{defn}

We can similarly define unlinked classes and  the reduced mapping schema of a rational lamination over $T$.
Hence we also define critical elements, post-critically finiteness and hyperbolicity of a rational lamination over $T$.
Note that the total number of critical elements of $\lambda_v$, counting multiplicities, is $\delta(v)-1$.

\begin{thm}
\label{thm-combinatorial onto}
Let $\lanot$ be an invariant rational lamination such that $\cR(\lanot) \neq \emptyset$.
Consider an internal angle system of $\lanot$ and let $\chi$ be the
 associated straightening map. \changetwo

There exists an algebraic set $Y=Y(\lanot) \subset
	\poly(\redschema(\lanot))$ of pure codimension one such that:
        \begin{itemize}
        \item 	if $P \in \cC(\redschema(\lanot))$, 
	then there exists $f \in \cR(\lanot)$ such that $\lambda_P=\lambda_{\chi(f)}$ 
	or $P$ is a non-hyperbolic map contained in $Y$.
        \item 	if $P \in \cC(\redschema(\lanot))$ is post-critically finite,
	then there exists $f \in \cR(\lanot)$ such that $P=\chi(f)$ 
	or $P$ is a non-hyperbolic map contained in $Y$.
        \end{itemize}
\end{thm}

Each irreducible component of the algebraic set $Y$ 
 is contained in one defined by the presence of parabolic periodic
 points of a given period or a preperiodic critical point of given
 preperiod and eventual period. \changetwo


 The proof of the theorem is at the end of this section. In~\ref{combinatorial straightening}
we discuss the effect of straightening on rational laminations (combinatorial straightening).
In~\ref{combinatorial tuning}, we introduce an inverse operation (combinatorial tuning).
In~\ref{subsec-combinatorial inverse} we show that combinatorial tuning and straightening are the inverse of each other.
In~\ref{subsec-renormalizable set} we characterize laminations with non-empty renormalizable set. Finally, in~\ref{subsec-proof combinatorial onto}
we prove the theorem above.


\subsection{Combinatorial straightening}
\label{combinatorial straightening}

\begin{defn}
 Let $\lanot$ and $\lambda$ be invariant rational laminations.
 Assume $\lanot$ has non-empty reduced schema $T(\lanot)$ 
and $\lanot \subset \lambda$.
 Let $\alpha = (\alpha_v: \overline{v} \to \S)$ be an internal angle
 system for $\lanot$.

 For each $v \in |\redschemanot|$, let $\lambda_v'$ be the equivalence
 relation that identifies $\theta$ and $\thetap$ if and only if there
 exist $t \in \alpha^{-1}_v (\theta)$ and $s \in \alpha^{-1}_v
 (\thetap)$ such that $t$ and $s$ are $\lambda$-equivalent.  Then
 $\lambda'=(\lambda_v')_{v \in |\redschemanot|}$ is a rational
 lamination over $\redschemanot$.  We say $\lambda'$ is the {\itshape
   combinatorial straightening of $\lambda$ with respect to
 $(\lanot,\alpha)$.}
\end{defn}

We let the reader check  that $\lambda'$ is in fact a rational lamination over  
$\redschemanot$.

%
%

\begin{lem}
  \label{lem-renormalized lamination}
Assume that $\lanot$ is an invariant rational lamination with non-empty reduced schema,
fix an internal angle system 
and denote by $\chi : \cR(\lanot) \to
 \cC(\redschema(\lanot))$ the corresponding straightening map. 

If  $f \in \cR(\lanot)$, then the rational lamination of $\chi(f)$ is the combinatorial
straightening of $\lambda_f$, with respect to $\lanot$ with the same internal angle system.
\end{lem}

\begin{proof}
Denote by $\alpha = (\alpha_v: \overline{v} \to \S)$ the internal angle system,
consider $f \in \cR(\lanot)$ and let $(\lambda_v)$ be  the combinatorial straightening  of
$\lambda_f$.

Let $g=(f^{\ell_v}:U'_v \to U_{\schemamap(v)})$ be a $\lanot$-renormalization.
Let $P = \chi(g): |\redschema(\lanot)| \times \C  \rightarrow |\redschema(\lanot)| \times \C$
be the straightening of $g$ via the fiberwise quasiconformal conjugacy
$\psi_v : U_v \to \psi_v(U_v) \subset \{v \} \times \C$. 
Let $\varphi_v : \C \setminus K(P,v) \to \C \setminus \overline{\Delta}$ 
be the B\"ottcher map.
Given $v \in |\redschema(\lanot)|$ and $\theta \in \overline{v} \cap \QS$,
denote by $\beta_v (\theta) \in \S \cong\partial \Delta$ the 
``landing point'' of the arc $\varphi_v \circ \psi_v (R_f (\theta) \cap U_v)$.
That is, the closure of this arc intersects $\partial \Delta$ at the point with
argument $\beta_v (\theta)$.

We claim that $\alpha_v$  and  $\beta_v$  coincide in $\overline{v} \cap \QS$.
In fact, these maps coincide at least in one point, since $\psi_v$ respects 
external markings. That is, if $\alpha_v (\theta_0) =0$, then
$\beta_v (\theta_0) =0$. Now we observe that the set of arguments where 
$\alpha_0$ and $\beta_0$ coincide is backward
invariant. More precisely, let $A \subset \overline{v}$ be the largest set such that 
$m^{\ell_v}_d (A) = \{ \theta \}$ where $\alpha_{\schemamap(v)} (\theta) =
\beta_{\schemamap(v)} (\theta) =0$.
Notice that
$\beta_{\schemamap(v)} \circ m^{\ell_v}_d = m_{\schemadeg (v)} \circ \beta_v$
and $\alpha_{\schemamap(v)} \circ m^{\ell_v}_d = m_{\schemadeg (v)} \circ \alpha_v$.
Therefore, $\alpha_v (A) = \beta_v (A)$. 
Also, 
the cyclic order of the arguments in $A$ is preserved by each of the
maps $\alpha_v$ and $\beta_v$. 
Therefore, $\alpha_v (\theta') = \beta_v (\theta')$, for all $\theta' \in A$. 
We may recursively repeat this argument in order to conclude that
 $\alpha_v$ and $\beta_v$ coincide at a set $C$ such that
$\alpha_v (C)$ is dense. It follows that $\alpha_v$ and $\beta_v$ agree on 
$\overline{v} \cap \QS$.

By Lindel\"of Theorem,  $\lambda_{P,v} \supset \lambda_v$, for all $v \in |\redschema(\lanot)|$.

Now consider two distinct $\lambda_v$-classes $A_1$ and $A_2$ and denote by $z_1$ and $z_2$ the landing points, in $K_f (v)$, of the rays with arguments 
in $\alpha^{-1}_v (A_1)$ and $\alpha^{-1}_v (A_2)$, respectively. Since $\psi_v$ is injective, the landing points  $\psi_v(z_1)$ and $\psi_v (z_2)$ of the external rays \change 
of $P$ in the $v$-fiber, with arguments in $A_1$ and $A_2$ are distinct. Therefore, $\lambda_{P,v} = \lambda_v$.
\end{proof}

\subsection{Combinatorial Tuning}
\label{combinatorial tuning}
In the classical context of quadratic polynomials, ``tuning'' is the inverse of straightening. 
Below we will introduce  the inverse of combinatorial
straightening.

\medskip
Let us first define the {\it pull-back of a relation}.
For a map $\varphi:A \to B$ and a relation $\lambda$ on a set $B$,
we say that  $a_1 \in A$ is $\varphi^*\lambda$-related to $ a_2 \in A$ 
if and only if
$\varphi(a_1)$ is $\lambda$-related to  $\varphi(a_2)$.

\begin{defn}
 \label{defn-comb tuning}
Consider a $d$-invariant rational lamination  $\lanot$  with an internal angle system 
$\alpha = (\alpha_v: \overline{v} \rightarrow \S)$ and an invariant rational lamination $(\lambda_v)_{v \in |\redschema (\lanot)|}$ 
over $\redschema (\lanot)$.

Given an infinite $\lanot$-unlinked class $v$, denote by $n_v \geq 0$ the smallest 
integer such that $v'= m_d^{n_v} (v) \in |\redschema (\lanot)|$. 
Let $\pi_{v} = \alpha_{v'} \circ m_d^{n_v}:  \overline{v} \rightarrow  \S$ and 
$\lambda^*_v = \pi_v^* (\lambda_{v'})$.

We define the {\itshape combinatorial tuning} $\cT_\alpha(\lanot,(\lambda_v))$ 
as the smallest equivalence relation in $\QS$ containing:
\begin{enumerate}
\item $\lanot$ and,
\item $\lambda_v^*$, for all infinite $\lanot$-unlinked classes $v$.
\end{enumerate}
\end{defn}

\begin{prop}
  Let $\lanot$ be a $d$-invariant rational lamination, $\alpha$ be an internal angle system of $\lanot$ and
 $(\lambda_v)_{v \in |\redschema (\lanot)|}$ be an invariant rational lamination
 over $T(\lanot)$. Then their combinatorial tuning
 $\lambda'= \cT_\alpha (\lanot,(\lambda_v))$ is a $d$-invariant
 rational lamination containing $\lanot$. Moreover, the following statements hold:

 \begin{enumerate}
 \item If $\lanot$ and $(\lambda_v)$ are hyperbolic, then $\lambda'$ is
 hyperbolic.
\item If $\lanot$ and $(\lambda_v)$ are post-critically finite, then $\lambda'$ is
 post-critically finite.
 \end{enumerate}
\end{prop}

\begin{proof}
   Let $\cA$ be the collection formed by the non-trivial $\lanot$-equivalence classes and
 let $\cB$ be the collection formed by all the 
 non-trivial $\lambda_v^*$-equivalence classes of all infinite
 $\lanot$-unlinked classes $v$.
 A non-trivial $\lambda'$-equivalence class $E$ is a maximal set of
 the form 
 \begin{equation}
  \label{eqn-lambda' class}
   E = \bigcup_{n=1}^N A_n
 \end{equation}
 for $A_n \in \cA \cup \cB$
 such that for any $n,m$, there exists a sequence 
 $n_0=n, n_1,\dots, n_k=m$ such that $A_{n_j}$ and $A_{n_{j+1}}$
 intersect for all $j=1,\dots,k-1$.
  Finiteness of $E$ follows from the fact that there exist  $\ell$ and $p$ such that
$m_d^\ell (\theta) = m_d^{\ell+p}(\theta)$ for all $\theta \in E$, since the same holds
for the elements of $\cA \cup \cB$.
 Thus,  every $\lambda'$-equivalence class is
 finite.

 Observe that if $A,B \in \cA \cup \cB$ are disjoint,
 then $A$ and $B$ are unlinked.
 In fact, consider $A_1, A_2, B_1, B_2 \in \cA \cup \cB$ 
 such that $A_1$ intersects $A_2$ and $B_1$ intersects $B_2$, and
 $A_1 \cup A_2$ and $B_1 \cup B_2$ are disjoint.
 Since $A_j$ is unlinked with both $B_1$ and $B_2$,
 we have that $B_1$ and $B_2$ are contained in the same component of $\R/\Z \setminus
 A_j$. Therefore, $A_1 \cup A_2$ and $B_1 \cup B_2$ are unlinked.
 We repeat this argument and conclude that $\lambda'$-equivalence
 classes are pairwise unlinked.

 To see that $\lambda'$ is closed, consider a sequence of pairs
 $(\theta_n,\theta_n')$ such that $\theta_n \sim_{\lambda'}
 \theta_n'$ and $\theta_n \to \theta$, $\theta_n' \to \theta'$ as $n \to
 \infty$. 
 If $\theta_n \sim_\lanot \theta_n'$ for infinitely many $n$,
 then $\theta \sim_\lanot \theta'$ because $\lanot$ is closed, 
 so it follows that $\theta \sim_{\lambda'} \theta'$.
 If there exists an infinite $\lanot$-unlinked class $v$ such that $\theta_n,
 \theta_n' \in v$ for infinitely many $n$, then similarly we have
 $\theta \sim_{\lambda'} \theta'$ because $\lambda_v^*$ is closed.
 Otherwise, by passing to a subsequence,
 we may assume that the convergences $\theta_n
 \nearrow \theta$ and $\theta_n' \searrow \theta'$ are monotone,
 that $\theta_n$ and $\theta_n'$ lie in an infinite $\lanot$-unlinked
 class $v_n$, and that $v_n$ ($n=1,2,\dots$) are pairwise disjoint.
 Then for each $n$, there exist $\iota_n \in (\theta_n,\theta_{n+1}),\
 \iota_n' \in (\theta_n',\theta_{n+1}')$ such that $\iota_n
 \sim_\lanot \iota_n'$ and that $v_n$ and $v_{n+1}$ lie in different
 components of $\R/\Z \setminus \{\iota_n,\iota_n'\}$.
 Therefore $\iota_n \nearrow \theta$ and $\iota_n' \searrow
 \theta'$, thus $\theta \sim_{\lanot} \theta'$.

 To prove the $d$-invariance of $\lambda'$,
 we start with the following observation:
 Let $v$ be an infinite $\lanot$-unlinked class.
 Then $\pi_v:\overline{v} \to \R/\Z$ in Definition~\ref{defn-comb tuning} 
induces a homeomorphism $\overline{v}/\lanot \to \R/\Z$.
 Furthermore, $\pi_{m_d(v)} \circ m_d \circ \pi_v^{-1}:\R/\Z \to \R/\Z$
 is well-defined equal to $m_{\delta(v)}$,
 and $\delta(v)_*\lambda_{v'}=\lambda_{w'}$,
 where $v'=m_d^{n_v}(v)$, $w=m_d(v)$ and $w'=m_d^{n_w}(w)$.
 Therefore, a $\lambda_v^*$-class (resp.\ $\lambda_v^*$-unlinked
 class) $A \subset \overline{v}$ is mapped by $m_d$ to a
 $\lambda_{m_d(v)}^*$-class 
 (resp.\ $\lambda_v^*$-unlinked class) and it is one-to-one
 if $v$ is not critical, and $\delta(\alpha_v(A))$-to-one
 if $v$ is critical.
 
 Let $E$ be a $\lambda'$-class of the form \eqref{eqn-lambda' class}.
 Then its image
 \[
 m_d(E) = \bigcup_{n=1}^N m_d(A_n)
 \]
 is contained in some $\lambda'$-class $E'$.
 If $E' \setminus m_d(E)$ is nonempty, then there exist some $A' \in
 \cA \cup \cB$ intersecting both $E' \setminus m_d(E)$ and
 $m_d(E)$.
 If $A' \in \cA$ then there exists a $\lanot$-class $A$ such that
 $m_d(A)=A'$ and $A$ intersects $E$. Moreover, $A$ is not contained in
 $E$ because $m_d(A)=A' \not \subset m_d(E)$.
 This contradicts that $E$ is maximal.
 Similarly, if $A'$ is a $\lambda_{v'}^*$-class for an infinite
 $\lanot$-unlinked class $v'$, then there exist an infinite
 $\lanot$-unlinked class $v$ and a $\lambda_v^*$-class $A$ such that
 $m_d(A)=A'$ and $A$ intersects $E$ but not contained in $E$,
 that is a contradiction.
 Therefore, $m_d(E)=E'$ is a $\lambda'$-equivalence class.
 
 Before proving that  $m_d: E \to m_d(E)$ is consecutive preserving
we establish that $\lambda'$-unlinked classes are invariant.
 So, consider a $\lambda'$-unlinked class $L$.
 Then either $L$ is a finite $\lanot$-unlinked class,
 or $L$ is contained in an infinite $\lanot$-unlinked class.
 If $L$ is a finite $\lanot$-unlinked class, then 
 $m_d(L)$ is also a finite $\lanot$-unlinked class, which is again
 a $\lambda'$-unlinked class. 
 On the contrary, if $L$ is contained in an infinite $\lanot$-unlinked
 class $v$, then $L$ is a $\lambda_v^*$-unlinked class.
 By construction, $m_d(L)$ is a $\lambda_{m_d(v)}^*$-unlinked class,
 which is also a $\lambda'$-unlinked class.
 
Now we prove that $m_d: E \to m_d(E)$ is consecutive preserving.
 Let $(\theta_1,\theta_2)$ be a component of $\R/\Z \setminus E$.
 Then there exist $\theta_i' \in (\theta_1,\theta_2)$ arbitrarily close
 to $\theta_i$ such that $\theta_1'$ and $\theta_2'$ are
 $\lambda'$-unlinked.
 If otherwise, there exist $\varepsilon>0$ such that $\theta_1'$ and $\theta_2'$ are  $\lambda'$-linked,
 for all  $\theta_1' \in (\theta_1,\theta_1+\varepsilon)$ and for all $\theta_2' \in
 (\theta_2-\varepsilon, \theta_2)$.
 Hence there exist  $t_n \in (\theta_1,\theta_1+1/n)$
 and $s_n \in (\theta_1+\varepsilon,\theta_2-\varepsilon)$ such that
 $t_n$ and $s_n$ are $\lambda'$-equivalent.
 By passing to a subsequence, we may assume $s_n$ converges to some $s
 \in [\theta_1+\varepsilon,\theta_2-\varepsilon] \subset
 (\theta_1,\theta_2)$.
 Since we already proved that $\lambda'$ is closed, it follows that
 $\theta_1$ and $s$ are $\lambda'$-equivalent and $s \in E$,
 which is a contradiction.
 Thus for any $\theta \in (d\theta_1,d\theta_2)$, there exists some 
 $\theta_1', \theta_2' \in (\theta_1,\theta_2)$ such that $\theta_1'$
 and $\theta_2'$ are $\lambda'$-unlinked and $\{\theta_1',\theta_2'\}$
 and $\{\theta',\theta_1\}$ are linked (not unlinked) for any $\theta'
 \in m_d^{-1}(\theta) \cap (\theta_1,\theta_2)$.
 Since $d\theta_1'$ and $d\theta_2'$ are also $\lambda'$-unlinked,
 $\theta$ is not $\lambda'$-equivalent to $d\theta_1$.
 Therefore $(d\theta_1,d\theta_2)$ is a component of $\R/\Z \setminus
 m_d(E)$ and we have proved $m_d:E \to m_d(E)$ is consecutive preserving.
 
It follows that, $\lambda'$ is a $d$-invariant rational lamination.

For (i), if $\lanot$ and $(\lambda_v)$ are hyperbolic,
 then for each Fatou critical element $\tilde{v}'$ of $\lambda_v$,
 $v'=\alpha_v^{-1}(\tilde{v}')$ is a critical element of $\lambda'$ 
 with $\delta(v')=\delta(\tilde{v}')$. Moreover, the total number of critical elements of $\lambda_v$, counting multiplicities, is $\delta(v)-1$.
 Therefore,
 \[
  \sum_{v \in \in |T(\lanot)|} \sum_{\tilde{v}' \in \crit(\lambda_v)} (\delta(\tilde{v}')-1) = \sum_{v \in |T(\lanot)|} (\delta(v)- 1) = d-1.
 \]
It follows that  all the critical elements of $\lambda'$ are infinite $\lambda'$-unlinked classes.
That is, $\lambda'$ is also hyperbolic.
 
For (ii), let us further assume that $\lanot$ and $(\lambda_v)$ are
 post-critically finite.
 Let $v'$ be a finite $\lambda'$-unlinked class.
 Then $v'$ is either a finite $\lanot$-unlinked class,
 or a finite $\lambda_v^*$-unlinked class for some infinite
 $\lanot$-unlinked class $v$. Hence $v'$ is not critical by assumption.
 Therefore, $\lambda'$ is post-critically finite.
\end{proof}

\subsection{Combinatorial straightening and combinatorial tuning}
\label{subsec-combinatorial inverse}
Here we prove that combinatorial tuning is the inverse of combinatorial
straightening:
\begin{thm}
 \label{thm-comb-inverse}
 Let $\lanot$ be a $d$-invariant rational lamination and 
 $\alpha=(\alpha_v)$ be an internal angle system for $\lanot$.
 Then:
 \begin{enumerate}
  \item \label{item-comb-str-tun=id}
	For any invariant rational lamination $(\lambda_v)_{v \in
	|T(\lanot)|}$ over $T(\lanot)$, the combinatorial straightening of
	$\cT_\alpha (\lanot,(\lambda_v))$ with respect to $(\lanot,\alpha)$
        is $(\lambda_v)_{v \in
	|T(\lanot)|}$.
  \item \label{item-comb-tun-str=id}
	If $\lambda \supset
	\lanot$ is a  $d$-invariant rational lamination, then
$\lambda= \cT_\alpha(\lanot,(\lambda_v))$ 
where  $(\lambda_v)$ is  the combinatorial
straightening of $\lambda$ (with respect to $(\lanot,\alpha)$).
 \end{enumerate}
\end{thm}

\begin{cor}
 \label{cor-inj-comb-str}
 Let $\lanot$ be a $d$-invariant rational lamination.
 For $f, g \in \cR(\lanot)$, if the rational laminations of 
 the straightenings $\chi_{\lanot}(f)$ and $\chi_{\lanot}(g)$ are equal,
 then $\lambda_f=\lambda_g$.
\end{cor}

\begin{proof}
 This follows from the fact $\lambda_f$ and $\lambda_g$ are the
 combinatorial tunings of $\lanot$ with the rational laminations of
 $\chi_{\lanot}(f)$ and of $\chi_{\lanot}(g)$, respectively, by
 Lemma~\ref{lem-renormalized lamination} and the above theorem. 
\end{proof}

For the proof of 
Theorem~\ref{thm-comb-inverse} it is convenient to state (without proof) the following direct consequences of the definitions of tuning and straightening.

\begin{lem}
 \label{lem-str/tun-unlinked-class}
 Let $\lanot$ be a $d$-invariant rational lamination and
 $\alpha=(\alpha_v)$ be an internal angle system for $\lanot$.
 Then:
 \begin{itemize}
  \item For any invariant rational lamination $(\lambda_v)_{v \in
	|T(\lanot)|}$ over $T(\lanot)$, let
	$\lambda'=\cT_\alpha(\lanot,(\lambda_v))$ be the
	combinatorial tuning of $\lanot$.
	Then any $\lambda'$-unlinked class contained in $v \in
	T(\lanot)$ has the form $\alpha_v^{-1}(L)$ for some
	$\lambda_v$-unlinked class $L$.
  \item For any $d$-invariant rational lamination $\lambda \supset
	\lanot$, let $(\lambda_v)_{v \in |T(\lanot)|}$ be its
	combinatorial straightening with respect to $\lanot$.
	Then any $\lambda_v$-unlinked class has the form $\alpha_v(L)$
	for some $\lambda$-unlinked class $L \subset v$.
 \end{itemize}
\end{lem}

\begin{proof}[Proof of Theorem~\ref{thm-comb-inverse}]
 \ref{item-comb-str-tun=id} Denote by $(\lambda_v')$ the combinatorial straightening of
$\cT_\alpha (\lanot,(\lambda_v))$. From the definitions it follows that $\lambda_v' \supset \lambda_v$. Now consider a $\lambda_v$-unlinked class $L$, then $\alpha_v^{-1}(L)$ is 
a $\cT_\alpha (\lanot,(\lambda_v))$-unlinked class, by the previous lemma. Again from the lemma above, $L=\alpha_v(\alpha_v^{-1}(L))$
is a $\lambda_v'$-unlinked class. Now by Lemma~\ref{lem-strict-inclusion}, $\lambda_v=\lambda_v'$ for all 
$v \in |T(\lanot)|$.

 \ref{item-comb-tun-str=id} Let
 $\lambda'=\cT_\alpha(\lanot,(\lambda_v))$. \changetwo
 By definition each  $\lambda'$-class is a 
(finite) union of $\lanot$-classes and $\lambda^*_v$-classes where $v$ is
a $\lanot$-unlinked class and $\lambda^*_v$ is as in Definition~\ref{defn-comb
 tuning}.  
 Since  the equivalence classes of these two relations
 are contained in $\lambda$-classes, we have $\lambda' \subset \lambda$.
 
 On the other hand, for a $\lambda$-unlinked class $L$ contained in an
 infinite $\lanot$-unlinked class $v$, 
 we have $L=(m_d^{n_v}|_v)^{-1} \circ \alpha_{v'}^{-1} \circ \alpha_{v'}
 \circ m_d^{n_v}(L)$, where notations are as in Definition~\ref{defn-comb
 tuning} (observe that $m_d^{n_v}|_v:v \to v'$ is a one-to-one
 order-preserving bijection). 
 Therefore, again by Lemma~\ref{lem-str/tun-unlinked-class}, $L$ is a
 $\lambda'$-unlinked class.
 Hence it follows that $\lambda'=\lambda$ by
 Lemma~\ref{lem-strict-inclusion}. 
\end{proof}

\subsection{Renormalizable set}
\label{subsec-renormalizable set}
Our aim now is to 
characterize laminations with non-empty renormalizable set.

\begin{prop}
\label{prop-nonempty-renor}
  Consider a $d$-invariant rational lamination such that $\redschema (\lanot) \neq \emptyset$.
Then the following statements are equivalent:
\begin{enumerate}
\item \label{item-nonempty-cR}
      $\cR(\lanot) \neq \emptyset$.
\item \label{item-no crit class on bdry}
If $A$ is a critical $\lanot$-class and $L$ is an infinite unlinked class in the forward orbit of a  critical element, then
$A \cap \overline{L} = \emptyset$.
\end{enumerate}
\end{prop}

Before proving the above result let us introduce an important class of rational laminations.

\begin{defn}[Primitive invariant laminations]
 We say that a $d$-invariant rational lamination is {\itshape primitive}
 if there do not exist infinite $\lambda$-unlinked classes $L$ and $L'$
 and a $\lambda$-class $A$ such that $L \ne L'$ and both $\overline{L}
 \cap A$ and $\overline{L'} \cap A$ are nonempty.
\end{defn}

As a direct consequence of the proposition we obtain the following result.

\begin{cor}
  If $\lanot$ is a primitive rational lamination such that $T(\lanot)\neq \emptyset$, 
then $\cR (\lanot) \neq \emptyset$.
\end{cor}

To prove the proposition, we start with the lemma below, which gives  sufficient conditions for a polynomial $f \in \cC(\lanot)$ to be in $\cR(\lanot)$.
The reader may find the proof of the proposition after the proof of the lemma. 
In the lemma, we simultaneously obtain an analytic family of polynomial like maps over
$\redschemanot$ parametrized by a neighborhood of $f$ in $\poly(d)$ (see Section~\ref{subsec-intropolylike} for the definition of analytic families).  We will use this analytic family to prove Theorem~\ref{introthm-onto hyperbolic} in Section~\ref{sec-onto hyperbolic}.

\begin{lem}
\label{lem-sufficient renor}
  Consider an invariant rational lamination $\lanot$.
Assume that (ii) of 
 Proposition~\ref{prop-nonempty-renor} holds.
Given a combinatorial Yoccoz puzzle $\Lambda = (\Lambda_k)$
  which is a generator of $\lanot$, there exists $K >0$ such that $f \in \cR(\lanot)$ if  $f \in \cC(\lanot)$ and $f$ satisfies the following two conditions: 

\begin{itemize}
 \item For all $v \in CO(\lanot)$, if   $\theta \in \supp \Lambda_0 \cap \overline{v}$, then
the external ray $R_f(\theta)$ lands at a repelling periodic point  and,
 \item For all $v \in CO(\lanot)$, if  $\theta \in \supp \Lambda_K \cap \overline{v}$, \change
then the landing point
       of the external ray $R_f(\theta)$ is not a critical point,
\end{itemize}
 Moreover, if $f \in \cC(\lanot)$ satisfies the above conditions, then  there exist a neighborhood $\cV$ of $f$, in $\poly(d)$,   that parametrizes 
  an analytic family of polynomial-like maps over $\redschemanot$
 \[
 (h_g = (g^{\ell_v}:U_{g,v}' \to U_{g,\sigma(v)})_{v \in
 |\redschemanot|})_{g \in \cV}
 \]
 such that the boundaries of the Jordan domains $U_{g,v}$ and
 $U_{g,v}'$ are smooth,
 the fiberwise connectedness locus contains $\cC(\lanot) \cap \cV$, and
 $h_g$ is a $\lanot$-renormalization for $g \in \cC(\lanot) \cap \cV$.
\end{lem}

\begin{proof}
  We claim that under the hypothesis of the lemma $K_f (v) = \bigcap P_{k} (f, v)$. 
Directly from the definitions it follows that $K_f (v) \subset \bigcap P_{k} (f, v)$.
Now if $z \notin K_f (v)$, then there exists a sector $S$ bounded by $\lanot$-equivalent rays with arguments in $\overline{v}$
 that is disjoint from $K_f (v)$ and such that $z \in S$. Since $\Lambda$ is a generator, there exists $k$ and a sector $S'$ bounded
by $\Lambda_{k}$-equivalent arguments such that $z \in S' \subset S$, by Proposition~\ref{prop-generator}.
Thus, $z \notin P_{k} (f, v)$.
 
 Now in order to show that $f \in \cR(\lanot)$ 
 we construct polynomial-like mapping over $T(\lanot)$ 
 by modifying puzzle pieces
 (cf. \cite[Section 3]{Haissinsky-Mandelbrot}, \cite[Lemma
 12.5]{Kiwi-thesis} and \cite{Milnor2}). 
 More precisely, let $k_0$ be a separation depth. 
 By Proposition~\ref{prop-generator} and the assumption,
 we can take $k_0$ sufficiently
 large so that for all $k>k_0$ and for all $v$ in the forward element of
 a Fatou critical element of $\lanot$, there exists no critical
 $\Lambda_k$-class $A$ which intersects $\overline{L_k(v)}$,
 where $L_k(v) \in \cL_k(\Lambda)$ is the combinatorial puzzle
 piece of depth $k$ containing $v$.

 Let $K=k_0+1$
 and let $f \in \cC(\lanot)$ satisfy the assumption,
 i.e., the landing point of $R_f(\theta)$ for every
 $\theta \in \partial L_K(v) \cap \overline{v}$ is neither parabolic nor
 critical for any $v \in CO(\lanot)$.
 By the choice of $k_0$, all the critical elements in
 $\overline{L_K(v)}$ are contained in $\overline{v}$, so it follows that
 the same property holds for any $\theta \in \partial L_K(v)$.
That is,  the landing point of $R_f(\theta)$ for every
 $\theta \in \partial L_K(v) \cap \overline{v}$ is neither parabolic nor
 critical for any $v \in CO(\lanot)$.

 We prove that there exists a desired analytic family of
 polynomial-like maps over $\redschemanot$ on a neighborhood of $f$.
 For $v \in CO(\lanot)$, let $Z(v)$ be the set of all points which are
 the landing point of $R_f(\theta)$ for some $\theta \in \partial
 L_K(v)$ and let $Z = \bigcup_{v \in CO(\lanot)} Z(v)$. 
 If $z \in Z$ is periodic, then it is repelling, by assumption, hence it has a holomorphic
 continuation in a neighborhood of $f$. More precisely, there
 exists a repelling periodic point $z_g$ for $g \in \poly(d)$
 sufficiently close to $f$ such that $z_g$ moves continuously on $g$ and
 $z_f=z$.
 Moreover, the set of landing angles of $z_g$ for $g$ contains that of $z$
 for $f$, by the continuity of the landing angle for repelling periodic
 points (note that this statement includes that even if $g \not \in
 \cC(d)$, such external rays do not bifurcate).
 
 Also, if $z \in Z$ is eventually periodic, we have a similar
 continuation; let $f^n(z) \in Z$ be a repelling periodic point.
 Since the local degree of $f^n$ at $z$ is one, there exists a
 unique preperiodic point $z_g$, for any $g$ close to $f$, of the same
 period and preperiod as $z$.
 Since landing angles of $f^n(z)$ are landing angles of $f^n(z_g)$
 (for $f$ and $g$ respectively), so do $z$ and $z_g$.

 It follows that there exists a neighborhood $\cV$ of $f$ such that
 the graph $\Gamma_g(v) = \overline{\bigcup_{\theta \in \partial
 L_K(v)} R_g(\theta)}$ depends continuously on $g \in \cV$.
 (More precisely, since the B\"ottcher coordinates and eventually periodic
 points depend holomorphically on $g \in \cV$, it admits a holomorphic
 motion $\cV_f \times \Gamma_f(v) \to \Gamma_g(v)$.)
 

 Fix small $\delta>0$. 
 For each $\theta$ let 
 \[
  N_f(\theta) = \varphi_f^{-1}(\{e^{r+2\pi i \eta};\
 |\eta-\theta| < \delta r^2\}).
 \]
 We may assume $N_f(\theta) \cap D_f(r_0)$ ($\theta \in
 \supp(\Lambda_K)$) are mutually disjoint 
 where $r_0$ is the level of the equipotential used for level $0$ puzzle
 pieces.
 For each $v \in CO(\lanot)$, let
 \[
  \tilde{P}_k(f,v)= P_k(f,v) \setminus \bigcup_{\theta \in L_k(v)}
 N_f(\theta).
 \]
 By the continuity of $\Gamma_g(v)$ and $\varphi_g$, we may assume
 $P_k(g,v)$ and $\tilde{P}_k(f,v)$ are defined similarly for all $g \in
 \cV_f$ and $k \le K$.
 
 For each $v \in CO(\lanot)$ and each $z \in Z(v)$, we may choose a small
 simply connected open neighborhood 
 $D_z$ and $D_z'$ around $z$ satisfying the following:
 \begin{itemize}
  \item $\partial D_z$ is smooth,
  \item $D_z \cap D_w  = \emptyset$ if $z \ne w$,
  \item $D_z' \Subset D_z$,
  \item $f(D_z') = D_{f(z)}$, and
  \item $f|_{D_z}$ is injective.
 \end{itemize}
 In fact, for each periodic orbit $\cO \subset Z$,
 let $p= \# \cO$ be the period of $\cO$ and $k$ be the largest integer
 such that there exists $w \in Z$ satisfying
 $f^{k-1}(w) \not \in \cO$ and $f^k(w) \in \cO$.
 Take $z \in \cO$ and choose a sequence of $p$ nested topological disks centered at $z$
 \[
  D_z^{(k)} \Subset \cdots \Subset D_z^{(k+p-1)}
 \]
 such that $f^p(D_z^{(k)}) \Supset D_z^{(k+p-1)}$.
 For each $w \in Z \cap f^{-k}(\cO)$, let $D_w$ be the component of
 $f^{-n}(D_z^{(n)})$ containing $w$, where $n \ge 0$ is the smallest
 integer such that $f^n(w)=z$.
 If $D_z^{(k+p-1)}$ is sufficiently small, then $f^n:D_w \to D_z^{(n)}$
 is univalent. Let $D_w' = (f^n|_{D_w})^{-1}D_z^{(n-1)}$.
 Then it is easy to check that $D_z$ and $D_z'$ satisfies the desired
 condition if $D_z^{(k+p-1)}$ is sufficiently small.
 
 We may assume $D_z$ also moves continuously on $g \in \cV_f$.
 In fact, for $g$ sufficiently close to $f$, $z_g$ is still repelling.
 Hence for such $g$, let $D_z(g)=D_z$ and $D_z'(g) =
 g^{-1}(D_{f(z)}) \cap D_z(g)$.

 Let us denote $D_z^* = D_z \setminus \{z\}$.
 \begin{claim*}
  For each $v \in CO(\lanot)$ and $z \in Z(v)$,
  there exists a finite union of disjoint $C^1$ arcs $\gamma(v,z) \subset
  D_z^*$ such that:
  \begin{enumerate}
   \item \label{thickening-claim:endpoints}
	 The endpoints of each component of $\gamma(v,z)$ are in
	 $\Int \tilde{P}_{K-1}(f,v)$.
   \item \label{thickening-claim:components}
	 For each component $B$ of $D_z \setminus \tilde{P}_{K-1}(f,v)$
	 with $z \in \partial B$,
	 exactly one component of $\gamma(v,z)$ intersects $B$ and $B
	 \setminus \gamma(v,z)$ consists of two components.
	 In particular, let $S_f(v,z)$ be the union of components of
	 $D_z \setminus (\tilde{P}_{K-1}(f,v) \cap \gamma(v,z))$ which
	 contain $z$ in the boundary. Then $S_f(v,z) \cup
	 \tilde{P}_{K-1}(f,v)$ is a neighborhood of $z$.
   \item \label{thickening-claim:transversality}
	 $\gamma(v,z)$ intersects an arc $\partial N_f(\theta)
	 \cap P_{K-1}(f,v)$ exactly once and transversally.
   \item \label{thickening-claim:expanding}
	 Let $S_f'(v,z) = (f|_{D_z})^{-1}(S_f(m_d(v),f(z)))$.
	 Then $S_f(v,z) \supset \overline{S_f'(v,z)} \setminus
	 \tilde{P}_{K}(f,v)$.
  \end{enumerate}
 \end{claim*}
 Note that $\partial N_f(\theta)$ consists of two arcs
 $\varphi_f^{-1}(\{e^{r+2\pi i \eta};\ \pm \delta(\eta-\theta) = r^2\})$,
 only one of which intersects $P_{K-1}(f,v)$,
 which we denote by $\tilde{R}_f(\theta,v)$.

 \begin{proof}[Proof of the claim]
  First consider the case when $v \in CO(\lanot)$ and $z \in Z(v)$ are
  periodic. 
  Let $n>0$ be the period of $v$
  (it is the period of $z$ multiplied by the denominator of the
  combinatorial rotation number of $z$).
  Let $T(v,z)=D_z^*/(w \sim f^n(w))$.
  It is isomorphic to a torus $\C^*/(w \sim (f^n)'(z)w))$ by
  the linearizing coordinate at $z$.
  There is a natural isomorphism $\iota_{z,v}:T(v,z) \to T(m_d(v),f(z))$
  induced by $f:f^{-1}(D_{f(z)}^*) \cap D_z^* \to D_{f(z)}^*$.
  
  Let $Q(v,z) = (P_{K-1}(v,z) \cap D_z^*)/{\sim} \subset T(v,z)$.
  The boundary $\partial Q(v,z)$ consists of simple closed curves
  of the form $(R_f(\theta) \cap D_z^*)/{\sim}$ in the same homotopy class.
  In particular, both $Q(v,z)$ and $T(v,z) \setminus Q(v,z)$ are unions 
  of disjoint annuli, whose core curves are also mutually homotopic.

  Let $A$ be a component of $T(v,z) \setminus Q(v,z)$.
  Observe that the Green potential $g_f=\log|\varphi_f|$ naturally
  induces a cyclic order on each component of $\partial A$.
  Take $n$ disjoint $C^1$ arcs
  $\hat{\gamma}_{A,0},\dots,\hat{\gamma}_{A,n-1} \subset  \overline{A}$
  such that 
  \begin{itemize}
   \item $\Int \hat{\gamma}_{A,k} \subset A$ and $\hat{\gamma}_{A,k}$
	 connects the two boundary components of $A$;
   \item at the endpoints, $\hat{\gamma}_{A,k}$ has a smooth
	 extension which intersects $\partial A$ transversally;
   \item The endpoints of $\hat{\gamma}_{A,k}$ are ordered with respect
	 to the cyclic orders on components of $\partial A$, where index
	 $k \in \{0,\dots,n-1\}$ is understood modulo $n$.
  \end{itemize}
  Let $\theta_1$, $\theta_2$ satisfy 
  \begin{equation}
   \label{eqn:theta_i}
   \partial A = ( (R_f(\theta_1) \cup R_f(\theta_2)) \cap D_z^* )/{\sim}.
  \end{equation}
  Then $(\tilde{R}_f(\theta_i, v) \cap D_z^*)/{\sim}$ accumulates 
  $(R_f(\theta_i) \cap D_z^*)/{\sim}$ in the $C^1$ topology.
  Therefore, we fix the extension above so that it intersects 
  $(\tilde{R}_f(\theta_i, v) \cap D_z^*)/{\sim}$ transversally
  if $\delta>0$ is sufficiently small.
  
  Now we take representatives $\gamma_{A,0},\dots,\gamma_{A,n-1} \subset
  D_z^*$ of the extension of
  $\hat{\gamma}_{A,0},\dots,\hat{\gamma}_{A,n-1}$ as follows:
  Take a point $w \in R_f(\theta_1) \cap D_z^*$ such that the
  equivalence class $[w]$ is an endpoint of $\hat{\gamma}_{A,n}$
  and $f^n(w) \in D_z^*$.
  Then $\gamma_{A,k}$ is a lift of the extension of $\hat{\gamma}_{A,k}$
  intersecting the subarc $[w,f^n(w))$ of $R_f(\theta_1)$.

  If we take $w$ sufficiently close to $z$ and $\delta>0$ is
  sufficiently small,
  then $f^k(\gamma_{A,n-k})$ intersects $\tilde{R}_f(m_d^k(\theta_i),
  m_d^k(v))$  transversally.
  For $0 \le k < n$, let us define
  \[
   \gamma(m_d^k(v),f^k(z))= \bigcup_{A \in \mathcal{A}(v,z)}
  f^k(\gamma_{A,n-k}), 
  \]
  where $\mathcal{A}(v,z)$ is the set of connected components of $T(v,z)
  \setminus Q(v,z)$.
  By considering all such periodic orbits for $(v,z) \mapsto (m_d(v),f(z))$,
  we can define $\gamma(v,z)$ for all periodic $v \in CO(\lanot)$ and $z
  \in Z(v)$.

  For such $v,z$ and $A \in \mathcal{A}(v,z)$, 
  there exists exactly one component of $\gamma(v,z)$ (say $\gamma_A$)
  such that $\gamma_A/{\sim}$ intersects $A$.
  Take $\theta_1$, $\theta_2$ satisfying \eqref{eqn:theta_i}.
  Then $\tilde{R}_f(\theta_i,v)$ intersects $\gamma_A$ once,
  so it follows that there is a unique bounded complement of
  $\overline{\tilde{R}_f(\theta_1,v)} \cup
  \overline{\tilde{R}_f(\theta_2,v)} \cup \gamma_A$, which we denote by
  $S_A$.
  Let $B$ be the component of $D_z^* \setminus \tilde{P}_{K-1}(f,v)$
  with $B/\sim = A$.
  Since $\partial B \subset \overline{\tilde{R}_f(\theta_1,v)} \cup
  \overline{\tilde{R}_f(\theta_2,v)} \cup \partial D_z $
  and $\gamma_A$ does not intersect $\partial D_z$,
  $B \setminus \gamma_A$ consists of two components,
  one of which is equal to $S_A$.
  Hence we have proved \ref{thickening-claim:endpoints},
  \ref{thickening-claim:components} and
  \ref{thickening-claim:transversality} because $S_f(v,z) = \bigcup_{A
  \in \mathcal{A}(v,z)} S_A$.
  
  By construction, $f(\gamma_A)$ and $\gamma_{\iota_{z,v}(A)}$ are
  disjoint and the Green potential of $f(R_f(\theta_1) \cap \gamma_A)$ 
  is greater than $R_f(d\theta_1) \cap \gamma_{\iota_{z,v}(A)}$.
  Hence it follows that the domain bounded by $\overline{R_f(d\theta_1)}$,
  $\overline{R_f(d\theta_2)}$ and $\gamma_{\iota_{z,v}(A)}$ is contained in
  the image by $f$ of the domain bounded by $\overline{R_f(\theta_1)}$,
  $\overline{R_f(\theta_2)}$ and $\gamma_{\iota_{z,v}(A)}$.
  Therefore, we have $f(S_A) \supset \overline{S_{\iota_{z,v}(A)}} \setminus
  \tilde{P}_{K-1}(f,m_d(v))$.
  By taking the inverse image by $f|_{D_z}$, the property
  \ref{thickening-claim:expanding} follows.

  For $v \in CO(\lanot)$ and preperiodic $z \in Z(v)$,
  there exists some $j>0$ such that $f^j(z)$ is periodic.
  By the above argument, we may assume that $\gamma_{A,k}$ as above are defined
for
  $(m_d^j(v),f^j(z))$.
  Taking the inverse image by $f^j|_{D_z}$ of $f^{nl}(\gamma_{A,k})$
  for some $l$ and $k$, where $n$ is the period of $m_d^j(v)$,
  we can similarly define $\gamma(v,z)$ with the required properties.
  Thus, we have concluded with the proof of the claim.
 \end{proof}
 
 Since $\Gamma_g(v)$ and $D_z$ depends holomorphically on $g \in \cV$,
 we may assume the claim holds not only for $f$, 
 but also for $g \in \cV$ with the same $\gamma(v,z)$,
 by replacing $\cV$ with a smaller neighborhood, if necessary.
 
 For $g \in \cV_f$, let $\hat{P}_{K-1}(g,v)$ and $\hat{P}_K(g,v)$ be the
 interiors of 
 \[
 \begin{matrix}
  \displaystyle \tilde{P}_{K-1}(g,v) \cup 
  \bigcup_{z \in Z(v) \cap \partial P_{K-1}(g,v)} S_g(v,z), \\
  \displaystyle (g^{-1}(\tilde{P}_{K-1}(g,m_d(v))) \cap P_{K}(g,v)) \cup
  \bigcup_{z \in Z(v)} S_g'(v,z).
 \end{matrix}
 \]
 respectively.
 Then $g:\hat{P}'_K(g,v) \to \hat{P}_{K-1}(g,m_d(v))$ is a proper map of
 degree $\delta(v)$. 
 Moreover, since
 \[
  g^{-1}(\tilde{P}_{K-1}(g,m_d(v))) \cap P_K(g,v) \subset \Int
 \tilde{P}_{K-1}(g,v) \cup Z(v)
 \]
 by definition,  it follows that $\hat{P}'_K(g,v)$ is a relatively
 compact subset of $\hat{P}_{K-1}(g,v)$.
 Therefore, we obtained a polynomial-like map
 \[
  h^U=(g:\hat{P}_K'(g,v) \to \hat{P}_{K-1}(f,m_d(v)))
 \]
 over the 
 {\itshape unreduced} mapping schema $T^U(\lanot)=(CO(\lanot),m_d,\delta)$.
 (Namely, it is an {\itshape unreduced} $\lanot$-renormalization of
 $f$. See Section~\ref{sec-injectivity} for more details of unreduced
 mapping schema and renormalization.) 
 
 By construction, every point in $\hat{P}'_K(g,v) \setminus P_K(g,v)$
 is contained in $D_z \setminus P_K(g,v)$ for some $z \in Z(v)$,
 hence it eventually escapes from the domain of $h^U$.
 Since $P_K(g,v) \setminus \hat{P}'_K(g,v) \subset \C \setminus K(g)$,
 it follows that $K_g(v) \subset \hat{P}'_K(g,v)$ for $g \in \cR(\lanot)$.
 Therefore, the filled Julia set $K(h^U,v)$ is equal to $K_g(v)$.
 In particular, it is connected.

 We may now extract a polynomial-like map over the reduced mapping
 schema from the one over the unreduced mapping schema.
 That is, for $v \in |T(\lanot)|$, let $U_{g,v}=\hat{P}_{K-1}(g,v)$ and
 let
 \[
  U_{g,v}' = (g|_{\hat{P}_{K-1}(g,v)})^{-1} \circ
 (g|_{\hat{P}_{K-1}(g,m_d(v))})^{-1} \circ \cdots \circ 
 (g|_{\hat{P}_{K-1}(g,m_d^{\ell_v-1}(v))})^{-1} (U_{g,v}).
 \]
 Therefore, every $g \in \cV \cap \cC(\lanot)$ is $\lanot$-renormalizable.
 
 It remains to show that the obtained family of polynomial-like mappings
 $h_g = (g^{\ell_v}:U'_{g,v} \to U_{g,\sigma(v)})_{v \in |T(\lanot)|}$
 forms an analytic family.
 Let 
 \begin{align*}
  \cU_v &= \bigcup_{g \in \cV} \{g\} \times U_{g,v}, &
  \cU_v' &= \bigcup_{g \in \cV} \{g\} \times U_{g,v}'.
 \end{align*}
 Since $\partial U_{g,v}$ and $\partial U_{g,v}'$ 
 are piecewise $C^1$ curves which move continuously on $g \in \cV$,
 $\cU_v$ and $\cU_v'$ are homeomorphic to
 $\cV \times \Delta$ over $\cV$ for any $v \in |\redschemanot|$
 by \cite[Theorem~2.17]{Riedl}, where $\Delta$ is
 the unit disk.
 Moreover, we have $\overline{\cU_v'}=\bigcup_{g \in \cV}
 \{g\} \times \overline{U}_{g,v}'$ , hence it follows that 
 the natural projection $\cU_v' \to \cV$ and $\cU_v' \ni (g,z) \mapsto 
 (g, g^{\ell_v}(z)) \in \cU_{\sigma(v)}$ are proper.
 Therefore, $h_g$ is an analytic family of polynomial-like mappings.

\end{proof}

As a consequence of the proof we obtain the following result:

\begin{cor}
 \label{cor-renor K}
 Under the assumption and notations of Lemma~\ref{lem-sufficient renor},
 for any $g \in \cV$, the following statements hold:
 \begin{enumerate}
  \item $\Lambda_K \subset \lambda_g$. In particular, the puzzle piece
	$P_k(g,v)$ is well-defined for $k \le K$ and $v \in CO(\lanot)$.
  \item $K(h_g,v) = \{z \in \C;\ g^n(z) \in P_K(g,m_d^n(v)) 
	\mbox{ for all }n \ge 0\}$.
 \end{enumerate}
\end{cor}

For later use, we further need to study the case when $\lanot$ is
hyperbolic.
\begin{lem}
\label{lem-hyperbolic thickening}
 Under the notations in Lemma~\ref{lem-sufficient renor},
 assume that $\lanot$ is hyperbolic.
Given $g \in \cV$,  if $K(h_g)$ is fiberwise connected,
 then $g \in \cC(\lanot)$.
\end{lem}

\begin{proof}
Note that if $K(h_g)$ is fiberwise connected, then $g \in \cC(d)$ since $\lanot$ is hyperbolic all the critical points of $g$ lie in a fiber of the domain of
$h_g$. It follows that all external rays are well defined for $g$. 
 Thus, it is sufficient to show that $\Lambda_n \subset \lambda_g$ for all $n
\ge K$, since the smallest rational lamination that contains the puzzle
is $\lanot$, by Corollary~\ref{cor-generator}.

For all $k \ge K$, let $\Lambda_k^g$ be the restriction of $\lambda_g$
to the $\supp \Lambda_k$ and let $\Lambda^g = (\Lambda^g_k)$.  We
proceed by induction to show that $\Lambda_n \subset \lambda_g$ and that
for all $v \in |T(\lanot)^U|$, 
$$K(h_g, v) \subset \bigcup_{L \in \cL_n(\Lambda^g),\ L \cap v \neq
 \emptyset} P_n(g,L).$$
By Corollary~\ref{cor-renor K}, this statement holds for $n=K$.
Hence we assume that the above statements are true for $n$ and 
proceed to establish the
corresponding statements for $n+1$.

Let $A$ be a $\Lambda_{n+1}$-class which is not in the support of
$\Lambda_n$.  Let $L'_n$ be the $\Lambda_n$-unlinked class containing
$A$.

If $m_d: L'_n \to m_d(L'_n)$ is injective, then from the unlinked and
invariant property of $\lambda_g$ it follows that $A$ is contained in a
$\lambda_g$-class, since $m_d(A)$ is a $\Lambda_n (\subset \lambda_g)$-class
and elements of $A$ cannot be $\lambda_g$-equivalent points in $\partial
L'_n$ (they have different preperiod).

Now assume that $m_d: L'_n \to m_d(L'_n)$ is $\delta$-to-one for some
$\delta > 1$. 
Since $K$ is greater than the separation depth of $\Lambda$,
it follows that there exists a critical Fatou element $v$ of
$\lanot$ such that $v \subset L'_n$ and $\delta(v)=\delta$.
Thus, there exists a $\Lambda_{n+1}$-unlinked class $L'_{n+1}$ 
such that $v \subset L'_{n+1}$.

Let $(\theta_0, \theta_1)$ be a connected component of $\R/\Z \setminus
L'_{n+1}$. We claim that $\theta_0$ and $\theta_1$ are
$\lambda_g$-equivalent.  By contradiction, assume that $\theta_0$ and
$\theta_1$ are not $\lambda_g$-equivalent. Then there exists a
$\Lambda^g_{n+1}$-unlinked class $L$ that separates $\theta_0$ and
$\theta_1$ such that $m_d(L)$ is disjoint from $L_n=m_d(L'_{n+1})$. 
In fact, $m_d(L)$ is contained in a $\Lambda_n$-class $L'$.
Since $m_d^{-1}(L_n) \cap L'_n = L_{n+1}'$ and $L$ is not
contained in $L'_{n+1}$, $L'$ is disjoint from $L_n$.
Thus, by the inductive hypothesis, the interior of $P_n(g,m_d(L))$ is free of points in
$K(h_g ,m_d(v))$, it follows that the interior of $P_{n+1}(g,L)$ is
disjoint from $K(h_g,v)$. 
However, $v$ must intersect the component of $\R/\Z \setminus L$
containing $\theta_0$ and the one containing $\theta_1$.  Thus, the
interior of $P_{n+1}(g,L)$ separates $K((h_g,v)$ which contradicts the
assumption that $K((h_g,v)$ is connected. We conclude that $\theta_0$ and $\theta_1$ are $\lambda_g$-equivalent.

We may assume $A$ is contained in $[\theta_0,\theta_1]$ above.
Since $m_d([\theta_0,\theta_1] \cap L'_n)$ is disjoint from
$L_n$ and $m_d:L_n \to m_d(L_n)$ has degree $\delta(v)$ we have that 
$m_d$ is an order-preserving bijection on
$[\theta_0,\theta_1] \cap L'_n$.
Therefore, as in the case when $L'_n \to m_d(L'_n)$ is injective,
$A$ is contained in a $\lambda_g$-class.

Let $v \in |T(\lanot)^U|$ and $L'_n \in \cL_n(\Lambda)$ and $L'_{n+1} \in
 \cL_{n+1}(\Lambda)$ satisfy $v \subset L'_{n+1} \subset L'_n$.
We must  prove that for $L \in \cL_{n+1}(\Lambda^g)$ disjoint with $v$, 
we have that $K(h_g,v)$ and the interior of $P_{n+1}(g,L)$ are disjoint.
Since we have already proven $\Lambda_{n+1} \subset \Lambda_{n+1}^g$,
$L$ and $L_{n+1}$ are disjoint. 
Hence $L$ is contained in a component $(\theta_0,\theta_1)$ of $\R/\Z
\setminus L'_{n+1}$. Therefore, $\Int P_{n+1}(g,L)$ and $K(h_g,v)$ are
separated by $\overline{R_g(\theta_0) \cup R_g(\theta_1)}$, hence
they are disjoint.
\end{proof}

As another  consequence of Lemma~\ref{lem-sufficient renor} we obtain that $\cR(\lanot)$ is a large subset of 
$\cC(\lanot)$ provided that $\cC(\lanot)$ is not contained in a proper algebraic subset
of $\poly(d)$.

\begin{cor}
\label{cor-existence-X}
 Let $\lanot$ be an invariant rational lamination.
 Then there exists an algebraic set $X \subset \poly(d)$ of pure
 codimension one such that $\cR(\lanot) \supset \cC(\lanot) \setminus X$.
\end{cor}

\begin{proof}
Assume that (ii) of Proposition~\ref{prop-nonempty-renor} holds for $\lanot$.
Let $\Lambda = (\Lambda_k)$
  be a combinatorial Yoccoz puzzle which is a generator of $\lanot$
and  let $K$ be as in Lemma~\ref{lem-sufficient renor}.
Denote by $m$ the minimum common multiple of the periods of the periodic arguments which belong to $\Lambda_K$.
Let $X$ be the codimension one algebraic subset of $\poly(d)$ formed by the union of all polynomials $f$ such that $f^m$ has a multiple fixed point
with all polynomials $f$ which have a critical point $\omega$ such that $f^K(\omega)$ is a fixed point of $f^m$.  
It follows that, if $f \in \cC(\lanot) \setminus X$, then the external rays with arguments in $\Lambda_K$  do not land at critical points 
or at multiple periodic points. From Lemma~\ref{lem-sufficient renor} we conclude that $f \in \cR(\lanot)$.

 Now if (ii) of Proposition~\ref{prop-nonempty-renor} does not hold for $\lanot$, then
 $\cC(\lanot)$ is contained in a proper algebraic subset $X$ of
 $\poly(d)$. In fact, every $f \in \cC(\lanot)$ has a preperiodic critical
 point of preperiod and eventual period depending only on $\lanot$.
 (Indeed we have $\cR(\lanot) = \cC(\lanot) \setminus X = \emptyset$ by
 Proposition~\ref{prop-nonempty-renor}.)
\end{proof}

The last ingredient in the proof of Proposition~\ref{prop-nonempty-renor} 
is the following version of the main result stated in~\cite{Kiwi}.

\begin{thm}
\label{thm-realization}
Let $\lambda$ be a degree $d$ invariant rational lamination. Then there exists $f \in \poly(d)$ without neutral cycles such that $\lambda_f = \lambda$.
Moreover, $f$ can be chosen so that every Fatou critical point is periodic or eventually periodic. \change
\end{thm}

\begin{proof}
  Although the version above is not explicitly stated in~\cite{Kiwi} all the necessary elements are contained there.
  In fact, the proof of Theorem 1.1. (pp. 149--150 of~\cite{Kiwi}) starts establishing that the map $f$ is free of neutral periodic point.
  Also, from the statement of Corollary 7.3 and the definitions in Section 7.1 of~\cite{Kiwi}, it follows that all the Fatou critical points
of $f$ have finite forward orbit. 
\end{proof}

\begin{proof}[Proof of Proposition~\ref{prop-nonempty-renor}]
\ref{item-nonempty-cR} $\implies$ \ref{item-no crit class on bdry}. Assume that $f \in \cR(\lanot)$. Then there exists a 
polynomial-like map $g$ over $\redschema (\lanot)$ such that $K(g,v) = K_f (v)$ 
for all $v \in |\redschema (\lanot)|$. For each $v \in |\redschema (\lanot)|$,  there exist neighborhoods 
$U'_v \Subset U_v$ such that $f^{\ell_v} : U'_v \rightarrow U_{\sigma_\lanot (v)}$ is a proper map of degree $\schemadeg (v)$.
In order to conclude that \ref{item-no crit class on bdry} holds, we proceed by contradiction and suppose that there exists a critical $\lanot$-class $A$
such that $A \cap {m_d^\ell (\overline{v})} \neq \emptyset$ for some $0 \le \ell < \ell_v$.
Note that there are exactly $\schemadegnot (v)$ $\lanot$-classes $B_1, \dots, B_{\schemadegnot (v)}$ 
intersecting $\overline{v}$ that map into the $\lanot$-class $m^{\ell_v - \ell}_d (A)$ under $m^{\ell_v}_d$.
Let $A_1, \dots, A_{\schemadeg(v)}$ be the corresponding $\lambda_f$-classes and consider a reference element $t \in m^{\ell_v - \ell}_d (A)$.
It follows that the set 
$$ \{ t' \in A_1 \cup \cdots \cup A_{\schemadegnot (v)}; m^{\ell_v}_d (t') =t \}$$
has cardinality strictly greater than $\schemadegnot (v)$, since   $m^{\ell}_d (A_j) = A$ for at least one $j$.
Therefore, the neighborhood $U'_v$ of $K_f(v)$ contains the portion inside $D_f(r)$, for $r$ sufficiently small,
of at least $\schemadegnot (v)+1$ rays which map
onto the same ray under $f^{\ell_v}$. Hence $f^{\ell_v} : U'_v \rightarrow U_{\sigma_\lanot (v)}$ is not a degree $\schemadegnot (v)$ map.

\ref{item-no crit class on bdry} $\implies$ \ref{item-nonempty-cR}.
Assume that \ref{item-no crit class on bdry} holds, according to the previous theorem, there
exists a polynomial $f \in \cC(d)$ with rational lamination $\lambda_f
= \lanot$ and no neutral cycles.
Hence, we may apply Lemma~\ref{lem-sufficient renor} to a generator of $\lanot$ given by Lemma~\ref{lem-exists generator}
to conclude that $f \in \cR(\lanot)$.
\end{proof}

\subsection{Combinatorial surjectivity}
\label{subsec-proof combinatorial onto}
We will need the fact that post-critically finite maps over a reduced mapping schema are uniquely determined  by their rational lamination.
The corresponding fact is well known in the context  of  monic and centered polynomials. 
We give a precise statement and sketch a proof using P. Jones~\cite{Jones} result about holomorphic
removability of post-critically finite Julia sets.

\begin{thm}
  \label{rigid-pcf-thm}
  Let $T$ be a reduced mapping schema. Consider $P_1, P_2 \in \cC (T)$ post-critically finite polynomial maps over $T$ with
 rational laminations $\lambda_1= (\lambda_{1,v})_{v \in |T|}$ and  $\lambda_2= (\lambda_{2,v})_{v \in |T|}$, respectively. 

If $\lambda_1 = \lambda_2$, then $P_1 =P_2$.
\end{thm}

\begin{proof}
  For $j=1,2$ and for $v \in |T|$ periodic, the corresponding fiber $J(P_j,v)= \partial K(P_j,v)$ of the Julia set  of $P_j$  coincides with the Julia set of a monic centered post-critically finite polynomial. Thus, $J(P_j,v)$ is locally connected (e.g., see~\cite{Milnor}) and holomorphically removable (see~\cite{Jones}).
Since every fiber $\{v\} \times \C $ is eventually periodic under $P_j$, it follows that $J(P_j,v)$ is locally connected and holomorphically removable for 
all $v \in |T|$. 
  
  Now denote by $\psi_{j,v}: \C \setminus \overline{\Delta} \to \C \setminus K(P_j,v)$ the inverse of the B\"ottcher map.
Since $J(P_j,v)$ is locally connected, we have that  $\psi_{j,v}$ extends continuously to $\partial \Delta \equiv \R / \Z$.
Denote by $\lambda'_{j,v}$ the equivalence relation in $\R / \Z$ which identifies $s,t$ if and only if $\psi_{j,v} (s) = \psi_{j,v}(t)$.

It follows that $\lambda_{j,v} = \lambda'_{j,v} \cap (\Q/\Z \times  \Q/\Z)$. Moreover, since $P_j$ has no critical point with infinite 
forward orbit (in the boundary of a Fatou component), Lemma~4.17 in~\cite{Kiwi} generalizes to establish that 
the real lamination $\lambda'_{j,v}$ is the smallest equivalence relation in $\R/\Z$ which contains the closure of $\lambda_{j,v}$.

Now we assume that $\lambda_1 = \lambda_2$. From the previous paragraph we conclude that  $\lambda'_{1,v} =  \lambda'_{2,v}$. Thus, 
$\psi_{2,v} \circ \psi_{1,v}^{-1} : \C \setminus K(P_1,v) \to \C \setminus K(P_2,v)$ extends continuously to 
a map $h_v: \C \setminus \operatorname{int} K(P_1,v) \to  \C \setminus \operatorname{int} K(P_2,v)$ for all $v \in |T|$.

Note that the maps $h_v$ induce a conjugacy between $P_1$ and $P_2$ in the complement of the bounded Fatou components.
We may extend $h_v$ continuously to each bounded Fatou component $U$ as a conformal map. In fact, observe that $h_v (\partial U)$ must be the
boundary of a bounded Fatou component of $P_2$, say $V$. Consider the unique 
 conformal isomorphism $g_U$ between $U$ and $V$ such that its continuous extension agrees with $h_v$ at one  point in $\partial U$ and such that
$g_U$ maps the unique point in $U$ which lies in the grand orbit of a critical point of $P_1$ to the corresponding point for $P_2$ in $V$.
It is not difficult to check that the continuous extension to $\partial U$ of such a conformal isomorphism $g_U$ agrees with $h_v$ at every point.
It also follows that after  extending $h_v$ to every bounded Fatou component $U \subset K(P_1,v)$ by such a map $g_U$ we
obtain continuous map $h_v : \C \to \C$, which is conformal off the Julia set, for all $v \in |T|$.
Moreover, $(v,z) \mapsto (v,h_v(z))$ is a conjugacy between
$P_1$ and $P_2$. Since the Julia set of $P_1$ is holomorphically removable  we have that  $h_v$ is in fact holomorphic and invertible, 
thus affine, for all $v$.
It follows that $P_1=P_2$ since $h_v$ is tangent to the identity at infinity for all $v$.
\end{proof}

\begin{proof}[Proof of Theorem~\ref{thm-combinatorial onto}]
Let $(\Lambda_k)$ be a combinatorial Yoccoz puzzle which is a
generator for $\lanot$.  Consider $K$ as in Lemma~\ref{lem-sufficient  renor} and let $ m \in \N$ be such that $d^K \theta$ has
period dividing $m$ for all $\theta \in \Lambda_K$.  Let $Y \subset
\poly(\redschemanot)$ be the analytic set formed by all $P
\in \poly(\redschemanot)$ which possess a critical point $c$ such that
$P^K(c)$ is periodic of period dividing $m$.  

Now consider  $P
\in \cC(\redschemanot)$ and, by Theorem~\ref{thm-realization}, there exists a polynomial $f$ 
without neutral cycles such that $\lambda_f = \cT_\alpha (\lanot,
\lambda_P)$. It follows that $f \in \cR(\lanot)$ or $f$ has a Julia
critical point $c'$ such that $f^K(c')$ is periodic of period
dividing $m$ which is the landing point of an argument $\theta$ in
$\supp \Lambda_K \cap \overline{v}$.  
We claim  that $P \in Y$ in the latter case.
In fact, note that the $\lanot$-class of $\theta$
is not critical, for otherwise $\cR(\lanot)=\emptyset$, by
Proposition~\ref{prop-nonempty-renor}.  Hence, the $\lambda_f$-class
of $\theta$ must contain an argument $\theta' \neq \theta$ such that $d \theta = d
\theta'$.  It follows that $\{ \theta, \theta'\}$ is pairwise unlinked with
every non-trivial $\lanot$-class.  Therefore, the $\lanot$-class $A$  of $\theta$ as well as the 
$\lanot$-class $A'$  of $\theta'$ have non-empty intersection with 
$\overline{v'}$ for some critical element $v'$ of
$\lanot$.   We conclude that the elements of the disjoint sets $\alpha_{v'}(A \cap \overline{v'})$ and
$\alpha_{v'}(A' \cap \overline{v'})$ are $\lambda_{P,{v'}}$-equivalent.  Hence, for any $t \in \alpha_{v'}(A \cap \overline{v'})$ the
landing point of $R_P(v',t)$ is a critical point
$c$ of $P$ such that $P^{K} (c)$ is periodic of period
dividing $m$. That is, $P \in Y$.
We conclude that $P$ is a non-hyperbolic element of $Y$ or $f \in \cR(\lanot)$ is such that $\lambda_{\chi(f), v} = \lambda_{P,v}$.

Now if $P \in \cC(\redschema(\lanot))$ is post-critically finite,  then by the previous theorem $P$ is uniquely determined by
its rational lamination. Taking $f$ as in the previous paragraph,  we have that $f \in \cR(\lanot)$  and $\chi(f) = P$ or; $P$ is a non-hyperbolic element of $Y$.
\end{proof}


\section{Injectivity of straightening maps}
\label{sec-injectivity}
This section is devoted to establish that straightening maps are injective.
That is, we prove Theorem~\ref{introthm-injectivity}.
The proof relies on the following main result.

\begin{keylem}
 \label{keylem-zero-meas}
 Let $\lanot$ be a post-critically finite $d$-invariant rational lamination with reduced schema $T(\lanot)$.
If $f \in \cC(\lanot)$, then the set 
 \[
  F = K(f) \setminus \bigcup_{n \ge 0} \bigcup_{v \in |\redschemanot|} f^{-n}( K_f(v))
 \]
 has zero area.
\end{keylem}

Section~\ref{subsec-key} contains the proof of this lemma. 
In Section~\ref{subsec-proof injective} we prove Theorem~\ref{introthm-injectivity} assuming the Key Lemma. 

\begin{rem}
 Epstein and Yampolsky \cite{Epstein-Yampolsky} first proved this key
 lemma for cubic polynomials constructed by intertwining surgery of two
 quadratic polynomials.
 Ha\"{i}ssinsky \cite{Haissinsky-intertwine} also claimed the key lemma for
 a general intertwining surgery, but his proof contains a gap.
 (He claimed that a set containing $F_0$ in the proof below is
 empty. However, $F_0$ is always nonempty in the case of intertwining
 surgery.)
\end{rem}

\begin{rem}
 Theorem~\ref{introthm-injectivity} often implies that
the polynomial obtained through a quasiconformal construction from a collection
of polynomials is uniquely determined by  the original polynomials, up to affine conjugacy.
 For example, consider the case of intertwining surgery by Epstein and
 Yampolsky \cite{Epstein-Yampolsky}.
Starting with two quadratic polynomials with connected Julia sets they construct 
a cubic polynomial (called an intertwining) having two quadratic-like restrictions 
 whose Julia sets intersect at a repelling fixed point. The resulting cubic polynomial
only depends on the choice of initial quadratic polynomials and not on the choices made throughout the surgery.

 The first author \cite{Inou-intertwine} used the intertwining construction to obtain a
 polynomial having a capture type renormalization.
 For such a construction, Theorem~\ref{introthm-injectivity} can be applied 
to show that the resulting polynomial is
 uniquely determined by the initial data,  up to affine conjugacy.
 Ha\"{i}ssinsky \cite{Haissinsky-intertwine} also constructed intertwinings 
 at parabolic fixed points  which produce a non-renormalizable polynomial.
 However, the Key Lemma~\ref{keylem-zero-meas} still holds and, therefore,
the resulting polynomial of the construction is uniquely determined by the corresponding initial data, up to affine conjugacy.
\end{rem}

\subsection{Proof of Theorem~\ref{introthm-injectivity}}
\label{subsec-proof injective}
First, we assume the key lemma and prove Theorem~\ref{introthm-injectivity}.

 Let $\lanot$ be a post-critically finite $d$-invariant rational
 lamination and let $\alpha$ be an internal angle system for $\lanot$. 
 Consider two renormalizable polynomials $f_1, f_2 \in \cR(\lanot)$  such that  
 $\straightening_\lanot(f_1)=\straightening_\lanot(f_2)$, where $\straightening_\lanot$ is the corresponding
straightening map. We must show that $f_1 = f_2$.

For the purpose of this proof it is better to extract ``polynomial-like maps'' $g_1, g_2$  over the {\itshape unreduced 
mapping schema} for $\lanot$. 
More precisely, let $|T^U (\lanot)|$ be the (forward) orbit of the Fatou critical elements of $\lanot$
and, for all $v \in |T^U (\lanot)|$, let $\delta (v)$ be the degree of $m_d : v \to m_d(v)$.
It follows that $T^U (\lanot) = (|T^U (\lanot)|, m_d, \delta)$ is a mapping schema. 

It is easy to construct a polynomial-like map over the unreduced mapping
schema $T^U(\lanot)$ from a $\lanot$-renormalization:
\begin{lem}
 For $f \in \cR(\lanot)$, there exist domains $U'_v \Subset U_v$ ($v \in
 |T^U(\lambda_0)|$) such that
 $g = (g_v =f: U'_v \to U_{m_d(v)})$ is a polynomial-like map over
 $T^U(\lambda_0)$ and   $K(g,v)=K_f(v)$ for all $v \in |T^U(\lanot)|$.
\end{lem}
The proof is left to the reader.
We say that $g=(f:  U'_v \to U_{m_d(v)})_{v \in |\schema(\lanot)|}$ is
an {\itshape unreduced $\lanot$-renormalization} of $f$.

Using a hybrid conjugacy between unreduced renormalizations for $f_1$ and
$f_2$, we construct the following quasiconformal homeomorphism.

\begin{lem}
 There exist a generator $\Lambda=(\Lambda_k)$ of separation depth $k_0$
 for $\lanot$ and a quasiconformal homeomorphism $\Phi_0:\C \to \C$
 such that:
 \begin{enumerate}
  \item The Yoccoz puzzle of depth $k_0$ for $f_1$ is mapped to that for
	$f_2$, i.e., for any puzzle piece $P_{k_0}(f_1,L)$ of depth
	$k_0$ for $f_1$, we have $\Phi_0
	(P_{k_0}(f_1,L))=P_{k_0}(f_2,L)$.
  \item $\Phi_0=\varphi_2^{-1} \circ \varphi_1$ on $\C \setminus
	D_{f_1}(r/d^{k_0})$ where $\varphi_j:\C \setminus K(f_j) \to \C
	\setminus \{|z| \le 1\}$ is the B\"ottcher map for $f_j$.
  \item $\Phi_0 \circ f_1 = f_2 \circ \Phi_0$ on $K_{f_1}(v)$ for $v \in
	|T^U(\lanot)|$ and on $R_{f_1}(\theta)$ for $\theta \in \supp
	\Lambda_{k_0}$.
  \item $\frac{\partial \Phi_0}{\partial \bar{z}} \equiv 0$ a.e.\ on
	$K_{f_1}(v)$ for $v \in |T^U(\lanot)|$.  
 \end{enumerate}
\end{lem}

\begin{proof}
 Let $\Lambda = (\Lambda_k)$ be a generator of $\lanot$ and 
 $k_0 \ge 0$ be a separation depth. 
 We may assume that every Julia critical element of $\lanot$ 
 is a $\Lambda_k$-class for all $k \ge k_0$.  
 As in Section~\ref{subsec-puzzles}, for each $v \in |\schema(\lanot)|$, 
 let $L_k(v)$ be the combinatorial puzzle piece of depth $k \ge 0$ for
 $\Lambda$ containing $v$. Denote by $P_k(f_j,v)$  the corresponding puzzle
 piece for $(f_j,v)$.

 Let $E_v = \overline{v} \cap \partial L_{k_0}(v)$ and let
 $I_{f_j}(v)$ be the set of landing points of $R_{f_j}(\theta)$ for $\theta
 \in E_v$.
 It consists of repelling or parabolic (pre)periodic points in $K_{f_j}(v)$.
 For $x \in I_{f_1}(v)$, let $x_2 = x_2(x)$ be the landing point of
 $R_{f_2}(\theta)$ for any $\theta \in \ang_{f_1}(x) \cap E_v$   where $\ang_{f_1}(x)$ denotes the subset of $\Q/\Z$ formed by the
arguments of the rays of $f_1$ landing at $x$.
 The function $x_2:I_{f_1}(v) \to I_{f_2}(v)$ is well-defined because
 the rational laminations of $f_1$ and $f_2$ are the same by
 Corollary~\ref{cor-inj-comb-str}. 

 \begin{claim*}
  There exist unreduced $\lanot$-renormalizations $g_j=(f_j:U_{j,v}'
  \to U_{j,\sigma(v)})_{v \in |T^U(\lanot)|}$ for $f_j$ ($j=1,2$) and
  a quasiconformal homeomorphism $\psi_v:U_{1,v} \to U_{2,v}$ for $v \in
  |T^U(\lanot)|$ such that
  \begin{enumerate}
   \item $\frac{\partial \psi_v}{\partial \bar{z}} \equiv 0$ a.e.\ on
	 $K_{f_1}(v)$,  
   \item $\psi_{\sigma(v)} \circ f_1 = f_2 \circ \psi_v$ on $U_{1,v}'$,
   \item $\overline{U_{j,v}} \subset D_{f_j}(r/d^{k_0})$ for $j=1,2$.
   \item $\partial U_{1,v}$ and $\partial U_{2,v}$ are $C^1$ Jordan
	 curves and $\psi_v$ extends to a $C^1$-diffeomorphism between
	 them,
   \item for $x \in I_{f_1}(v) \cap K_{f_1}(v)$ and $\theta \in \overline{v}
	 \cap \ang_{f_1}(x)$, we have that
	 $\psi_v = \varphi_2^{-1} \circ \varphi_1$ on a neighborhood of  
	 $R_{f_1}(\theta) \cap U_{1,v}$.
  \end{enumerate}
 \end{claim*}

 In fact, to obtain $g_j$ and $\psi_v$ as above, first take an unreduced $\lanot$-renormalization $g_1=(f_1:U_{1,v}' \to
 U_{1,\sigma(v)})_{v \in |T^U(\lambda_0)|}$ defined in sufficiently small domains $U_{1,v}'$   such that
 $\partial U_{1,v}'$ and $\partial U_{1,v}$ are $C^1$ Jordan curves
 and $\overline{U_{1,v}} \subset D_{f_1}(r/d^{k_0})$.
 
 Continue by considering  an unreduced $\lanot$-renormalization $g_2=(f_2:U_{2,v}' \to
 U_{2,\sigma(v)})_{v \in |T^U(\lambda_0)|}$ such that  
 $\partial U_{2,v}'$ and $\partial U_{2,v}$ are $C^1$ Jordan curves and
 $\partial U_{2,v}$ and $\varphi_2^{-1} \circ \varphi_1(\partial U_{1,v})$
 coincide on a neighborhood of $R_{f_2}(\theta)$ for any $x \in
 I_{f_1}(v)$ and $\theta \in \overline{v} \cap \ang_{f_1}(x)$.
 
 Now let $A_{j,v} = \overline{U_{j,v}} \setminus U_{j,v}'$ be the
 fundamental annulus and we proceed to   define $\psi_v$ on $A_{1,v}$.
 Let $\psi_v = \varphi_2^{-1} \circ \varphi_1$ on a neighborhood of $A_{1,v} \cap R_{f_1}(\theta)$
 for each $v \in |T^U(\lambda_0)|$, $x \in I_{f_1}(v)$ and $\theta \in
 \ang_{f_1}(x)$.
 By construction and since $\varphi_2^{-1} \circ \varphi_1$ is a
 holomorphic conjugacy from $f_1$ to $f_2$ on the basin of infinity, 
 we can extend it to a $C^1$-diffeomorphism between $A_{1,v}$ and
 $A_{2,v}$ such that $\psi_{\sigma(v)} \circ f_1 = f_2 \circ \psi_v$ on
 $\partial U_{1,v}'$ and $\psi_{v}(A_{1,v} \cap P_{k_0}(f_1,v)) =
 A_{2,v} \cap P_{k_0}(f_2,v)$. 
 
 As in the proof of the straightening theorem (in
 Section~\ref{sec-poly-like}),
 we can extend $\psi_v|_{A_{1,v}}$ to a hybrid conjugacy on $U_{1,v}$ to
 $U_{2,v}$ by successive pullbacks.
 By construction, we have $\psi_v = \varphi_2^{-1} \circ \varphi_1$ on
 $U_{1,v} \cap R_{f_1}(\theta)$ for any $x \in I_{f_1}(v)$ and $\theta \in
 \ang_{f_1}(x)$, which implies the combinatorial condition needed to
 obtain a quasiconformal map $\psi_v$ after gluing the definition of
 $\psi_v$ in $U_{1,v} \setminus K(f_1,v)$
 with the original hybrid conjugacy.
 (If $I_{f_1}(v)$ is empty, then $K_{f_1}(v)$ is disjoint from $\partial
 P_{k_0}(v)$. Hence we need only choose $U_{1,v}$ sufficiently small
 such that $\overline{U_{1,v}} \subset P_{k_0}(v)$ and choose a
 diffeomorphism $\phi_v:A_{1,v} \to A_{2,v}$ properly so that the
 combinatorial condition is satisfied).

 Furthermore, since $\psi_v$ is defined by successive
 pullbacks, it coincides with $\varphi_2^{-1} \circ \varphi_1$ on a
 neighborhood of $R_{f_1}(\theta) \cap
 U_{1,v}$ for $\theta \in \overline{v} \cap \ang_{f_1}(\theta)$.
 Thus we have proved the claim.

 \medskip
 
 Let $\cI_{f_j} = \bigcup_{v \in |T^U(\lambda_0)|} I_{f_j}(v)$.
 For each $x \in \cI_{f_1}$, take small neighborhoods $O_x$ and
 $O_{x_2}$ of $x$ and $x_2=x_2(x)$ and a quasiconformal homeomorphism
 $\tau_x:O_x \to O_{x_2}$ such that 
 \begin{itemize}
  \item $(O_x)_{x \in \cI_{f_1}}$
	and $(O_{x_2(x)})_{x \in \cI_{f_1}}$
	are collections of pairwise disjoint open sets.
  \item $\tau_x(x)=x_2$,
  \item $\tau_x = \varphi_2^{-1} \circ \varphi_1$ on $O_x \cap
	f_1^{-1}(O_{f(x)}) \cap R_{f_1}(\theta)$ for any $\theta \in
	\ang_{f_1}(x)$. 
  \item if a puzzle piece $P$ of depth $k_0$ is not of the form
	$P_{k_0}(f_1,v)$ for any $v \in |T^U(\lanot)|$,
	then $\tau_x$ is $C^1$ on $(\Int P) \cap O_x$.
  \item $\tau_x=\psi_v$ on $P_{k_0}(f_1,v) \cap O_x$ if $x \in K_{f_1}(v)$,
 \end{itemize}
 Such homeomorphisms are obtained after lifting a $C^1$-diffeomorphism
 between two quotient annuli.
 When $x$ is periodic, these annuli are obtained from  punctured
 neighborhoods of $x$ and $x_2$ intersected with a puzzle piece $P$ of
 depth $k_0$ after identification of $z$ with $f_1^p(z)$, and $w$ with
 $f_2^p(w)$ respectively, where $p$ is the period of external
 rays landing at $x$.
 (Observe that even if $x$ (or $x_2$) is parabolic, such a puzzle piece $P$
 which is not of the form $P_{k_0}(f_i,v)$ is contained in a repelling
 petal because all immediate parabolic basins are contained in
 small filled Julia sets of $\lambda_0$-renormalizations.)
 For preperiodic $x$, $\tau_x$ is obtained by the pullback of
 $\tau_{f_1^n(x)}$ by the dynamics.

 Now we define $\Phi_0$. First we define it on the following three parts:
 \[
  \Phi_0 =
  \begin{cases}
   \varphi_2^{-1} \circ \varphi_1 & \mbox{on }\C \setminus
   D_{f_1}(r/d^{k_0}), \\
   \psi_v & \mbox{on }P_{k_0}(f_1,v) \cap U_{1,v}, \\
   \tau_x & \mbox{on } O_x \setminus \left(\bigcup_{v \in
   |T^U(\lambda_0)|}P_{k_0}(f_1,v)\right)\ (x \in \cI_{f_1}).
  \end{cases}
 \]
 By construction, $\Phi_0$ is quasiconformal and $C^1$ in a neighborhood 
 of the boundary of its domain of definition.   
 Furthermore, $\Phi_0$ maps $R_{f_1}(\theta)$ into $R_{f_2}(\theta)$ for
 $\theta \in \supp \Lambda_{k_0}$.

 Therefore, we can extend it first to $\Gamma_{k_0}= \bigcup_{\theta \in
 \supp \Lambda_{k_0}} R_{f_1}(\theta)$ and then to the rest of the plane by a
 $C^1$-diffeomorphism so that it  
 preserves the Yoccoz puzzle of depth $k_0$ and conjugates $f_1$ to
 $f_2$ on $\Gamma_{k_0}$. 
 
 On $K_{f_1}(v)$, we have
 \[
  \Phi_0 \circ f_1 = \psi_{\sigma(v)} \circ f_1 = f_2 \circ \psi_v 
  = \Phi_0 \circ f_2,
 \]
 and $\frac{\partial \Phi_0}{\partial \bar{z}} =
 \frac{\partial \psi_v}{\partial \bar{z}} \equiv 0$ a.e.  
\end{proof}

Note that $\Phi_0$ is a conjugacy in an arcwise connected set $S_0$ which includes all the critical orbits of $f_1$.
Also keep in mind that $S_0$ contains $\C \setminus D_{f_1}(r/d^{k_0})$ as well as 
$K_{f_1}(v)$ for all $v \in |\schema(\lanot)|$.
Therefore, we may recursively define $\Phi_n$, for $n \ge 1$, as the unique
(quasiconformal) homeomorphism such that
$\Phi_n \circ f_1 = f_2 \circ \Phi_{n-1}$
which agrees with $\Phi_{n-1}$ in the arcwise connected set  $S_{n-1} = f_1^{-(n-1)} (S_0)$. 
It follows that the quasiconformal dilatation constant of $\Phi_n$ agrees with that of 
$\Phi_0$ and, $\Phi_n$ is conformal in $\C \setminus
D_{f_1}(r/d^{k_0+n})$ as well as in $f^{-n}(K_{f_1}(v))$, for all $v \in
|\schema(\lanot)|$.
Since $\Phi_n$ is eventually independent of $n$ on a dense subset on a plane,
$\Phi_n$ converges to a quasiconformal map $\Phi$
that conjugates $f_1$ and $f_2$ in 
the dense set $\bigcup S_n$, and therefore in $\C$. 
Moreover, $\Phi$ is conformal in 
$$\bigcup_{n \geq 0} f^{-n} (\cK(|\schema(\lanot)|) \cup \left( \C \setminus K(f_1) \right),$$
which is a full measure set in $\C$, in view of  the Key Lemma.
Therefore $\Phi$ is affine. In fact, $\Phi$ is the identity since it is a conjugacy
between monic centered polynomials which is tangent to the identity at infinity.
\qed

\subsection{Proof of the key lemma}
\label{subsec-key}
\newcommand{\CO}{{CO}}

The rest of this section is devoted to prove Key
Lemma~\ref{keylem-zero-meas}.
Let $\Lambda=(\Lambda_k)$ be a generator of $\lanot$ and let $k_0\ge 0$
be a separation depth. Denote by $\schema(\lanot)$ the unreduced mapping schema of $\lanot$ as in the previous proof.

We introduce the following subsets of the Fatou critical orbit elements of $\lanot$ (i.e.,  $|\schema(\lanot)|$):
\begin{align*}
 \CO^0 &= \{v \in |\schema(\lanot)|;\ K_f(v) \cap
 \supp(f,\Lambda_{k_0+1})=\emptyset\},\\ 
 \CO^1 &= |\schema(\lanot)| \setminus \CO^0.
\end{align*}

\begin{figure}
 \begin{center}
  \fbox{\includegraphics[width=10cm]{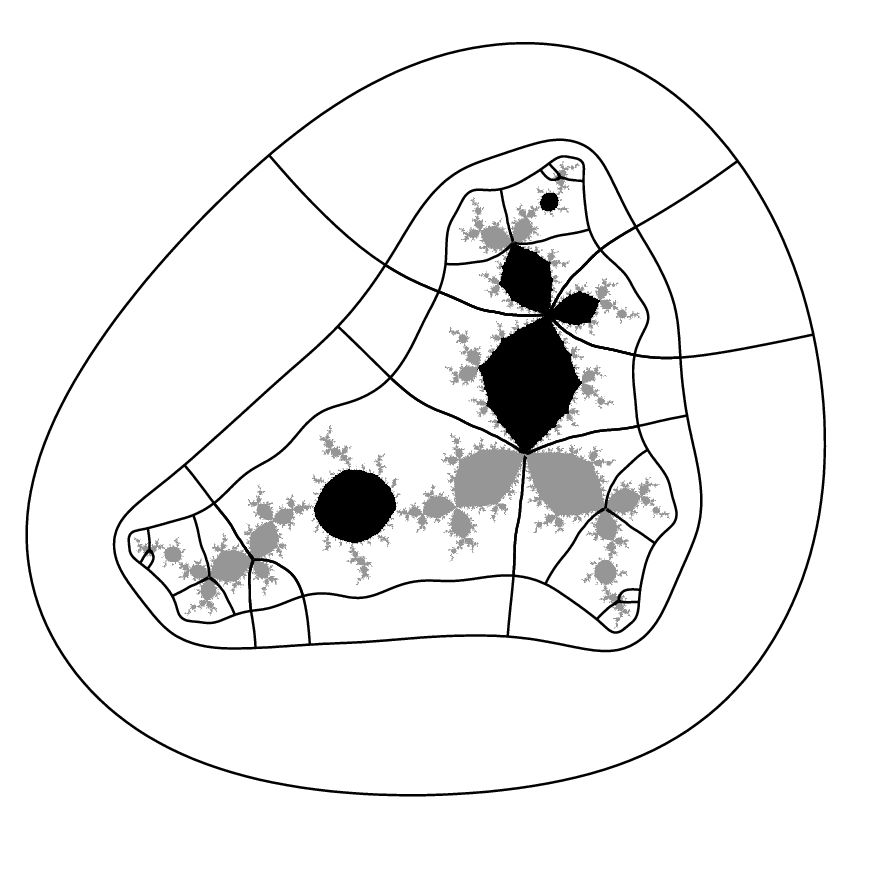}}
  \caption{Yoccoz puzzles of a separation depth.
  There are two periodic cycles of Fatou
  components, one of period two and the other of period three,  which are colored black.
  The period two cycle corresponds to $\CO^0$, and the period three
  cycle corresponds to $\CO^1$.}
  \label{fig-cubic2+3}
 \end{center}
\end{figure}

Denote by $\CO^{\rm per}$ the set formed by the periodic Fatou critical elements of $\lanot$.
Let $$\CO^{i,\rm per}=\CO^i \cap \CO^{\rm per},$$ for
$i=0,1$ (see Figure~\ref{fig-cubic2+3}). 

In the following, for a set $X \subset \C$, we denote by $\comp(X)$
 the set of
connected components of $X$ and let $\cN_{\epsilon}(X)=\{z \in \C;\
d(z,X)<\epsilon\}$ be the $\epsilon$-neighborhood of $X$.

Given a subset $S$ of $|\schema (\lanot)|$, we let
$$\cK_f (S) = \bigcup_{v \in S} K_f (v).$$

Fix $\epsilon_0>0$ sufficiently small,  so that
the neighborhoods $\cN_{\epsilon_0}(\cK)$  are pairwise disjoint for  
the components $\cK$ of  $\cK_f(\CO^{\per})$.
Recall that 
$$F = K(f) \setminus \bigcup_{n \ge 0} \bigcup_{v \in |\redschemanot|} f^{-n}( K_f(v)).$$

\begin{lem}
 \label{lem-a.e.z}
 The following hold for almost every $z \in F$\textup{:}
 \begin{enumerate}
  \item \label{item-Julia expansion}
	$\lim_{n\to \infty} \|(f^n)'(z)\|=\infty$ with respect to
	the hyperbolic metric on $\C \setminus PC(f)$\textup{;}
  \item \label{item-a.e.z-approach to K1,per}
	$\lim_{n \to \infty}
	d(f^n(z),\cK_f(\CO^{1,\per}))=0$\textup{;}
  \item \label{item-a.e.z-nbds}
	there exist $N=N(z)>0$ and $\cK \in
	\comp(\cK_f(\CO^{1,\per}))$ such that 
	$f^{N+n}(z) \in \cN_{\epsilon_0}(f^n(\cK))$, for all $n \geq 0$.
	Moreover, if $f^{N+n}(z) \in \cN_{\epsilon_0}(K_f(v))$
	for $v \in \CO^{1,\per}$, then $f^{N+n+1}(z) \in
	\cN_{\epsilon_0}(K_f(m_d(v)))$.
 \end{enumerate}
\end{lem}

\begin{proof}
 The assertion \ref{item-Julia expansion} is a stronger version of
 Ma\~{n}\'{e}'s Lemma due to  McMullen (see
 \cite[Theorem~3.6]{McMullen}).

 In view of \cite[Theorem~3.9]{McMullen}, for almost every point $z \in J(f)$ we have that
 \[
 \lim_{n \to \infty} d(f^n(z),PC(f)) =0.
 \]
 For all such $z \in F$, 
 \[
  \lim_{n \to \infty} d(f^n(z),\cK_f(\CO^{\per})) = 0
 \]
 because $PC(f) \setminus \cK_f(|\schema(\lanot)|)$ is a finite set
 consisting only of repelling periodic points and their backward
 images, and $f^n(\cK_f(|\schema(\lanot))|)=\cK_f(\CO^{\per})$
 for some $n \ge 0$.
 Since $\cK_f(\CO^{0,\per})$ and
 $\cK_f(\CO^{1,\per})$ are
 disjoint forward invariant compact sets, we have 
 \[
  \lim_{n \to \infty} d(f^n(z),\cK_f(\CO^{i,\per}))=0
 \]
 for some $i \in \{0,1\}$.

 Now assume  
 \begin{equation}
  \label{eqn-approaching to K0,per}
  \lim_{n \to \infty} d(f^n(z),\cK_f(\CO^{0,\per}))=0.
 \end{equation}
 For each $v \in \CO^0$, 
 $K_f(v)$ is contained in the interior of the puzzle piece
 $P(L_{k_0+1}(v))$. Therefore we have
 $K_f(v) \in \comp(\cK_f(\CO^{0,\per}))$ (i.e., $K_f(v)$
 is a component of $\cK_f(\CO^{0,\per})$).
 By \eqref{eqn-approaching to K0,per} and the continuity of $f$, it
 follows that there exist some $N \ge 0$ and $v \in \CO^0$ such
 that for any $n \ge 0$, $f^{N+n}(z) \in P(L_{k_0+1}(m_d^n(v)))$, which is a
 neighborhood of $K_f(m_d^n(v))$. Let $p$ be a period of $v$ under iterations of $m_d$.
 Then $f^{pn}(f^N(z))$ belongs to a neighborhood of $K_f(v)$ for all $n \ge 0$.
 Since $K_f(v)$ is the filled Julia set of a polynomial-like map, it
 follows that $f^N(z) \in K_f(v)$, which contradicts $z \in F$.
 Therefore, we have proved \ref{item-a.e.z-approach to K1,per}.

 The assertion \ref{item-a.e.z-nbds} easily follows from
 \ref{item-a.e.z-approach to K1,per} and the continuity of $f$.
\end{proof}

Now we define, for each $v \in
\CO^{1,\per}$,  an open set $\cN'(v)$, which is a slightly smaller
set than $\cN_{\epsilon_0}(K_f(v))$.
Each $\cK \in \comp(\cK_f(\CO^{1,\per}))$ can be
written as
\[
 \cK= \bigcup_{m=0}^M K_f(v_m), \quad v_m \in \CO^{1,\per}.
\]
Then the set 
\[
 I = \bigcup_{m \ne m'} (K_f(v_m) \cap K_f(v_m')) \subset
 \supp(f,\Lambda_{k_0+1})
\]
consists of repelling (pre)periodic points (see
\cite[Theorem~7.3]{McMullen} and \cite[Proposition~3.4]{Inou-ren}).
For each $x \in I$, take $\Theta(x) \subset \ang(x) \cap
\supp(\Lambda_{k_0})$
such that $d(\Theta(x))=\Theta(f(x))$ and each component of the
complement of 
\[
 \Gamma(x) = \bigcup_{\theta \in \Theta(x)} R_f(\theta) \cup \{x\}
\]
intersects exactly one of $K_f(v_m)$'s that contain $x$.
Let $\Gamma=\bigcup_{x \in I} \Gamma(x)$ 
and let $\Gamma_0$ denote the union of the components of $\Gamma \cap
\cN_{\epsilon_0}(\cK)$ intersecting $I$.
Then each component of $\cN_{\epsilon_0}(\cK) \setminus
\Gamma_0$ intersects exactly one of $K_f(v_1),\dots,K_f(v_M)$.
Let us denote by $\cN'(v_m)$ the component of
$\cN_{\epsilon_0}(\cK) \setminus
\Gamma_0$ containing $K_f(v_m) \setminus I$.

By construction, 
\begin{equation}
 \label{eqn-nbd partition}
  \overline{\bigcup_{m=1}^M \cN'(v_m)} =
  \overline{\cN_{\epsilon_0}(\cK)}
\end{equation}
and $\cN'(v)$ ($v \in \CO^{1,\per}$) are
pairwise disjoint.

\begin{lem}
 \label{lem-v(z,n)}

Let $F_0=\{z \in F;\ \mbox{Lemma~\ref{lem-a.e.z} holds}\}$. For $z \in F_0$, we have
 \begin{enumerate}
  \item \label{item-v(z,n)}
	For all $n \ge N(z)$, there exists $v(z,n) \in
	\CO^{1,\per}$ such that $f^n(z) \in \cN'(v(z,n))$\textup{;}
  \item \label{item-dv(z,n) ne v(z,n+1)}
	there exists arbitrarily large $n$ such that $m_d(v(z,n)) \ne
	v(z,n+1)$\textup{;}
  \item \label{item-epsilon1}
	there exists $\epsilon_1=\epsilon_1(\epsilon_0)>0$ independent
	of $z$ such that if
	$m_d(v(z,n)) \ne v(z,n+1)$, then $f^{n+1}(z) \in
	\cN_{\epsilon_1}(I)$ and $f^n(z) \in
	\cN_{\epsilon_1}((f^{-1}(I) \cap
	\cK_f(\CO^{1,\per}))\setminus I)$.
 \end{enumerate}
 Furthermore, $\lim_{\epsilon_0\to 0}\epsilon_1 =0$.
\end{lem}

\begin{proof}
 The assertion \ref{item-v(z,n)}  follows from the equation \eqref{eqn-nbd
 partition} above.
 If \ref{item-dv(z,n) ne v(z,n+1)} does not hold, then 
 we would have that 
 $f^n(z) \in K_f(v(z,n))$, for sufficiently large $n$ (since $\Lambda$ is a generator for
$\lanot$). This is a contradiction with $z \in F$.
 Finally, \ref{item-epsilon1} follows from the fact that
 $\cN_{\epsilon_0}(K_f(v)) \setminus \cN'(v)$ is
 contained in a small neighborhood of $I \cap K_f(v)$.
\end{proof}

Take a small neighborhood $\cO$ of $I$ such that $f(\cO)
\supset \overline{\cO}$ and $f|_\cO$ is injective.
Since $CO(f) \setminus \cK_f(\CO^{\per})$ is finite, we
may also assume the following:
\begin{enumerate}
 \item \label{item-no cv in O*}
       No critical value lies in $\cO^*= \cO \setminus I$, 
 \item \label{item-f^{-1}O and CO}
       If a component $\cO_1$ of $f^{-1}(\cO)$
       intersects $\cK_f(\CO^{1,\per})$, then 
       $\cO_1 \cap f^{-1}(I)$ and $\cO_1 \cap CO(f)$ are
       contained in $\cK_f(\CO^{1,\per})$.
\item $\tilde{\cO}=\cO^*/f$
is a union of tori.
\end{enumerate}

\begin{figure}
 \begin{center}
  \fbox{\includegraphics[width=10cm]{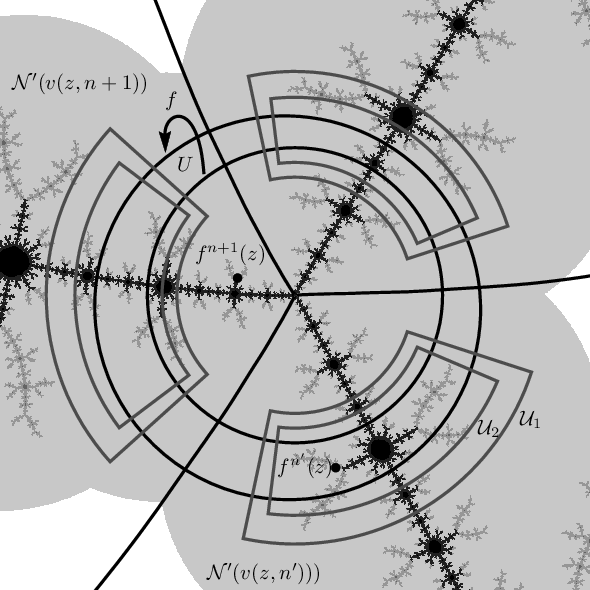}}

  \fbox{\includegraphics[width=10cm]{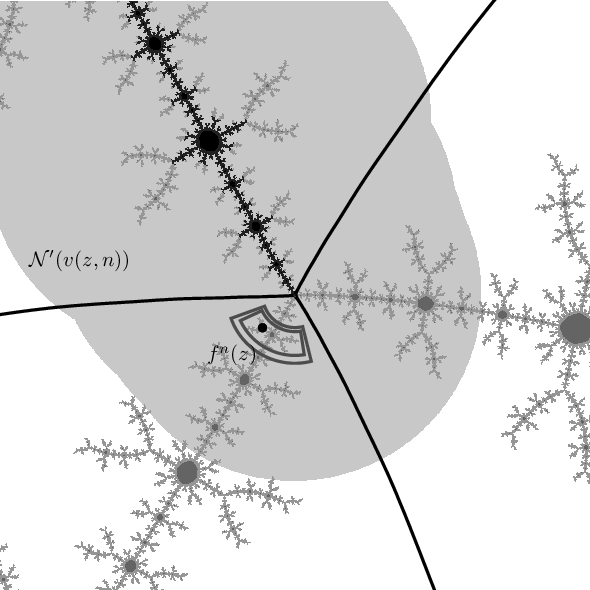}}
  \caption{Near $I$ and its inverse image.
  If $m_d(v(z,n)) \ne v(z,n+1)$, $f^n(z)$ must be close to $f^{-1}(I)
  \setminus I$.}
  \label{fig-fp}
 \end{center}
\end{figure}

Figure~\ref{fig-fp} illustrates the statement of the following lemma.

\begin{lem}
 \label{lem-O1O2}
 There exist open subsets $\cU_1$ and $\cU_2$ of $\cO$ such that the following statements hold:
 \begin{enumerate}
 \item $\cU_2 \Subset \cU_1 \Subset \cO$ and
       $\cU_1$ does not intersect $\Gamma$.
 \item Each one of the open sets $\cU_1$ and $\cU_2$ is a disjoint union of finitely many
       topological disks.
 \item Each component $A$ of $\cU_1$ contains exactly one
       component $B$ of $\cU_2$ and $A \setminus \overline{B}$
       is an annulus.
 \item $\cU_2/f \supset (\cO^* \cap K(f))/f$.
       In other words, there exists $\epsilon_2>0$ such that for any $z
       \in K(f) \cap \cN_{\epsilon_2}(I)$, there exists
       some $N>0$ such that $f^n(z) \in \cO$ for any $n$ with
       $0\le n <N$, $f^N(z) \in \cU_2$, and the
       branch of $f^{-N}$ sending $f^N(z)$ to $z$ can be univalently
       defined on the component of $\cU_1$ containing $f^N(z)$ 
 \end{enumerate}
\end{lem}

\begin{proof}
 (See Figure~\ref{fig-fp}.)
 Let $\tilde{K} = (K(f) \cap \cO^*)/f$ and 
 $\tilde{\Gamma}=(\Gamma \cap \cO^*)/f$.
 Note that $\tilde{\Gamma}$ and $\tilde{K}$ 
 are disjoint compact sets in $\tilde{\cO}$.
 Take a fundamental domain $U$ for the covering projection
 $\cO \setminus \Gamma \to \tilde{\cO} \setminus
 \tilde{\Gamma}$. 
 Let $\cU_2 \Subset \cU_1$ be neighborhoods of $U
 \cap K(f)$ in $\cO \setminus \Gamma$.
 It is easy to see that we can take $\cU_1$ and
 $\cU_2$ so that they satisfy the assertion.
\end{proof}

\begin{proof}[Proof of Key Lemma~\ref{keylem-zero-meas}]
 It is enough to show that any $z \in F_0$ is not a Lebesgue density point
 of $K(f) \supset F_0$.

 Let $\epsilon_1$ and $\epsilon_2$ be as in Lemma~\ref{lem-v(z,n)} and
 Lemma~\ref{lem-O1O2} respectively. 
 We may assume $\epsilon_1 \le \epsilon_2$.
 For $z \in F_0$, let 
 \[
 t(z) = \{n \ge N(z);\ m_d(v(z,n)) \ne v(z,n+1)\}.
 \] 
 By Lemma~\ref{lem-v(z,n)}, we have $t(z)$ contains infinitely many elements.
 For $n \in t(z)$, let $\cO_1(z,n)$ be the component of
 $f^{-1}(\cO \setminus \Gamma)$ containing
 $f^n(z)$.
 By the condition \ref{item-no cv in O*} on $\cO$ above,
 $f:\cO_1(z,n) \to \cO \setminus \Gamma$ is injective.
 (See Figure~\ref{fig-fp}.)

 For $z \in F_0$ and $n \in t(z)$, 
 $f^{n+1}(z)$ lies in a $\epsilon_1$-neighborhood of $I$ by
 Lemma~\ref{lem-v(z,n)}. Hence by Lemma~\ref{lem-O1O2}, there
 exists some $n'>n$ such that $f^{n'}(z) \in \cU_2$ and
 the branch of $f^{-n'+(n+1)}$ sending $f^{n'}(z)$ to $f^{n+1}(z)$
 is univalently defined on $\cU_1'(z,n)$,
 where $\cU_i'(z,n)$ is the component of $\cU_i$
 containing $f^{n'}(z)$ for $i=1,2$.
 Let $\tilde{\cU}_i(z,n)$ be the component of
 $f^{-n'+n}(\cU_i'(z,n))$ containing $f^n(z)$.
 Then $\tilde{\cU}_i(z,n) \subset \cO_1(z,n)$ and
 $f^{n'+n}:\tilde{\cU}_i(z,n) \to \cU_i'(z,n)$ is a
 conformal isomorphism.
 Since $\cO_1(z,n) \cap PC(f) \subset \cO_1(z,n) \cap
 CO(f) = \emptyset$, there exists a univalent inverse branch
 $f^{-n}:\tilde{\cU}_i(z,n) \to
 \cU_i(z,n)$ sending $f^n(z)$ to $z$.
 Therefore, $f^{n}:\cU_1(z,n) \to \tilde{\cU}_1(z,n)$ 
 and on $\cO$,
 $f^{n'}:\cU_1(z,n) \to
 \cU'_1(z,n)$ are conformal isomorphisms (condition \ref{item-f^{-1}O and CO}).

Denote by $\iota_{z,n}: \cU_1(z,n) \hookrightarrow  \C \setminus PC(f)$ the inclusion and
consider the following diagram:
 \[
 \xymatrix{
  \cU_1(z,n) \ar[r]^{f^n}_\cong \ar@{^{(}->}[d]^{\iota_{z,n}}
  & \tilde{\cU}_1(z,n) \ar[r]^{f^{n'-n}}_\cong \ar@{^{(}->}[d]
  & \cU_1'(z,n) \\
  \C \setminus PC(f) \ar@{.>}[r]^{f^n} & \C \setminus PC(f).
 }
 \]
 Since $\lim_{n\to \infty} \|(f^n)'(z)\|=\infty$ with respect to
 the hyperbolic metric on $\C \setminus PC(f)$ and inclusion does not
 expand hyperbolic metric, we have $\lim_{n \to \infty}
 \|\iota_{z,n}'(z)\| = 0$ with respect to the corresponding hyperbolic metric.
 By Koebe distortion theorem, this implies that the diameter
 $\cU_2(z,n)$ shrinks to zero with bounded distortion.
 Furthermore, since $f^{n'}:\cU_1(z,n) \to \cU_1'(z,n)$
 is a conformal isomorphism sending $\cU_2(z,n)$ to
 $\cU_2'(z,n)$, by the Koebe distortion theorem applied to the inverse of the conformal isomorphism
 $f^{n'}|_{\cU_1(z,n)}$, 
 there exist some $C,C'>0$ such that
 \[
 \frac{\area(\cU_2(z,n) \setminus
 K(f))}{\area(\cU_2(z,n))} 
 \ge C  \frac{\area(\cU_2'(z,n) \setminus
 K(f))}{\area(\cU_2'(z,n))} \ge C'.
 \]
 Here, we have only finitely many choices for $\cU_i'(z,n)$
 since each choice is a component of $\cU_i$. 
 This implies that $C$ and $C'$ are constants independent of $z
 \in F_0$ and $n \in t(z)$.
 
 Therefore, 
 \[
 \lim_{n \in t(z) \to \infty} \frac{\area(\cU_2(z,n) \cap
 K(f))}{\area(\cU_2(z,n))} \le 1-C' < 1
 \]
 and $z$ is not a Lebesgue density point of $K(f)$.
\end{proof}


\section{Onto hyperbolic maps}
\label{sec-onto hyperbolic}

The aim of this section is to prove Theorem~\ref{introthm-onto hyperbolic}.
First we will proof that hyperbolic maps are contained in the image of
straightening. 

\subsection{Onto hyperbolic maps}
Let $\lanot$ be a post-critically finite rational lamination with an internal angle system $\alpha_0$ such that
$|T(\lanot)| \neq \emptyset$ and $\cR(\lanot) \neq \emptyset$. Denote by $\chi$ the associated straightening map.
Consider a component $\cH$ of $\hyp \cC(\redschemanot)$.
Our aim now is to show that the image of $\chi$ contains $\cH$.

The image of $\chi$ contains at least one element of $\cH$. Indeed, 
consider the unique post-critically finite polynomial map $\hat{P} \in \cH$~\cite[Corollary~5.2]{Milnor-hyp}.
Since $\hat{P}$ is hyperbolic, by Theorem~\ref{thm-combinatorial onto},
there exists a post-critically finite polynomial  $\hat{f} \in \cR(\lanot)$ such that
$\chi (\hat{f})=\hat{P}$.

In order to prove that the image of $\chi$ contains $\cH$ we will
employ Milnor's parameterization of $\cH$~\cite{Milnor-hyp}.  Let $\hat{T} =
(|\hat{T}|,\hat\sigma, \hat\delta)$ be the reduced mapping schema of
$\hat{P}$. Namely, $|\hat{T}|$ is identified with the set formed by the
critical points of $\hat{P}$, $\hat{\sigma}: |\hat{T}| \to |\hat{T}|$ is the first
return map under iterations of $\hat{P}$, and $\hat{\delta} : |\hat{T}| \to \N$ is
the local degree of $\hat{P}$ at its critical points.  
Observe that the total degree of $\hat{T}$ is equal to
that of $\redschemanot$, because $\hat{P}$ is a hyperbolic polynomial map
over $\redschemanot$.
Consider the space $B(\hat{T})$ of all proper holomorphic maps
\[
 \beta:|\hat{T}| \times \Delta \to |\hat{T}| \times \Delta
\]
of the form $\beta(v,z)=(\hat{\sigma} (v),\beta_v(z))$
where $\beta_v:\Delta \to \Delta$ is of degree $\hat{\delta}(v)$ and the following
hold:
\begin{itemize}
 \item $\beta_v$ is 
{\itshape boundary-rooted}, i.e., the extension of $\beta_v$ to
       $S^1 = \partial \Delta$ satisfies $\beta_v(1)=1$;
 \item if $v$ is periodic under $\hat{\sigma}$,
       then $\beta_v$ is {\itshape fixed point centered}, i.e.,
       $\beta_v(0)=0$;
 \item if $v$ is not periodic, then $\beta_v$ is {\itshape critically
       centered}, i.e., the sum of its $\hat{\delta}(v)-1$ critical points
       (counted with multiplicity) is equal to zero. \changetwo
\end{itemize}
According to Milnor~\cite[Lemma~4.9]{Milnor-hyp} the space $B(\hat{T})$,
 endowed with the compact-open topology, is a topological
 cell of dimension $2n-2$ where $n$ is the total degree of $\hat{T}$. \changetwo
 Moreover, $B(\hat{T})$ 
 diffeomorphic to $\cH$~\cite[Theorem~5.1]{Milnor-hyp}. 
We proceed to describe a diffeomorphism  $\Phi: \cH \to B(\hat{T})$,
using the fact that $\cH$ is simply connected.

For all $v \in |\hat{T}|$ and all $P \in \cH$ we let $U_P (v)$ be the critical bounded Fatou component of $P$
such that  $U_P (v)$ 
depends continuously on $P$ (in the Caratheodory topology) and  $v \in U_{\hat{P}}(v)$. Such a (unique) labeling 
$U_P(v)$ exists, since $\cH$ is simply connected and the Julia set depends continuously on $P \in \cH$. 
It follows that $U_P(v') \neq U_P(v)$ if $v \neq v'$, and $P^{m_v} (U_P(v))=U_P(\hat{\sigma}(v))$ where
$m_v$ is the first return time of $v$ to $|\hat{T}|$ under iterations of $\hat{P}$.
Moreover, we may also consider a ''boundary marking''. That is, for all  
$v \in |\hat{T}|$ and all $P \in \cH$, we choose $z_P (v) \in \partial U_P(v)$  such that $z_P (v)$ varies continuously (in fact, holomorphically) with $P \in \cH$
and $P^{m_v}(z_P(v))=z_P(\hat{\sigma}(v))$. Again the existence of $z_P (v)$ relies on the fact that $\cH$ is simply connected.
Now, given $P \in \cH$, for all $v \in |\hat{T}|$ 
let $h_{P,v} : U_P(v) \to \Delta$ be the unique conformal isomorphism such that
the continuous extension of $h_{P,v}$ to $\partial U_P(v)$ maps $z_P(v)$ to $1 \in \partial \Delta$ and,
$\beta= \Phi(P): |\hat{T}| \times \Delta \to |\hat{T}| \times \Delta$ belongs to $B(\hat{T})$ where
$\beta(v,z)=(\hat{\sigma} (v),\beta_v(z))$ and 
$$\beta_v (z)= h_{P, \hat{\sigma}(v)} \circ P^{m_v} \circ h^{-1}_{P,v} (z).$$
(The existence of $h_{P,v}$ is guaranteed by~\cite[Corollary~4.16]{Milnor-hyp})
Milnor's Theorem (~\cite[Theorem~5.1]{Milnor-hyp}) states that $\Phi$ is a (real analytic) diffeomorphism.

\begin{lem}
  \label{lem-qc-surgery}
 There exists a differentiable map $\Psi:B(\hat{T}) \to \poly(d)$ such
 that $\Psi(B(\hat{T})) \subset \cR(\lanot)$ and $\Phi \circ \chi \circ
 \Psi:B(\hat{T}) \to B(\hat{T})$ is the identity.
 
 In particular, $\Phi'= \Phi \circ \chi: \Psi(B(\hat{T})) \to B(\hat{T})$ is onto.
\end{lem}

\begin{proof}
  We follow a  quasiconformal surgery construction which is a minor variation of
the one described in~\cite[Proof~of~5.7]{Milnor-hyp}.
  Given $\beta \in B(\hat{T})$ we will introduce a neighborhood $U_\beta$ which will be either contained in, or disjoint from, the image of $\Phi'$.
  Since $B(\hat{T})$ is a (connected) topological cell, the lemma will follow.

  Let $U_\beta$ be an open connected neighborhood of $\beta$ such that there exist $r_1 < r_2 <1$ with the properties that
all the critical points of $\beta_1$ are contained in $|\hat{T}| \times \Delta(r_1)$ and that the closure of $\beta_1(|\hat{T}| \times \Delta(r_2))$
is contained in $|\hat{T}| \times \Delta(r_1)$,  for all $\beta_1 \in U_\beta$. 

Now given $\beta_1, \beta_2 \in U_\beta$,  consider
the map $b =b(\beta_1,\beta_2): |\hat{T}| \times \Delta \to |\hat{T}| \times \Delta$  defined as $\beta_1$ outside $|\hat{T}| \times \Delta(r_2)$,
as $\beta_2$ in $|\hat{T}| \times \Delta(r_1)$, and by the linear interpolation 
$$\dfrac{|z| -r_1}{r_2-r_1} \beta_1(v, z) + \dfrac{r_2 - |z|}{r_2-r_1} \beta_2 (v,z)$$
for $r_1 \leq |z| \leq r_2$. Shrinking $U_\beta$, if necessary,  we may assume that for all $\beta_1, \beta_2 \in U_\beta$ the map $b$
is a local quasi-conformal homeomorphism.
Note that $b$ depends smoothly on $\beta_1$ and $\beta_2$ (in the appropriate smooth structure, see~\cite[Section~5]{Milnor-hyp}).

Assume that $\beta_1 \in U_\beta$ is in the image of $\Phi'$, say $\beta_1 = \Phi \circ \chi(f_1)$ and $P_1 = \chi (f_1)$. We proceed with a surgery 
which shows that 
$U_\beta$ is contained in the image of $\chi$. 

In fact, denote by $\psi = (\psi_w)_{w \in |T(\lanot)|}$ the hybrid conjugacy between the polynomial like map over $T(\lanot)$ extracted from
$f_1$ and $P_1$. For each critical Fatou component $U_{P_1}(v) \subset \{ w \} \times \C$ of $P_1$ denote by 
$V_{f_1}(v) = \psi_w^{-1}(U_{P_1}(v))$ the corresponding Fatou component of $f_1$.
 
For all  $\beta_2 \in U_\beta$ we may consider the local quasiconformal map
$b$, as above. 
Let $f_b : \C \to \C$ be the map which coincides with $f_1$ off the components $V_{f_1}(v)$ (i.e., $f_b(z) = f_1(z)$ provided 
$ z \notin V_{f_1}(v)$ for all $v \in |\hat{T}|$) and if $ z \in V_{f_1}(v)$ for some $v \in |\hat{T}|$, then
$$f_b(z) = \left( f_1^{n_v-1} |_{V'_{f_1}(v)} \right)^{-1} \circ \psi^{-1}_{w'} \circ h_{P_1, \hat{\sigma}(v)}^{-1} \circ b \circ h_{P_1,v} \circ \psi_w (z),$$
where $V'_{f_1}(v) = f_1(V_{f_1}(v))$,  the first return time of $V_{f_1}(v)$ to a critical Fatou component is $n_v \geq 1$, and $w' \in |T(\lanot)|$ 
is such that $U_{P_1}(\hat{\sigma}(v)) \subset \{ w' \} \times \C$.
It follows that $f_b$ is holomorphic off a finite union of annuli where the map is locally quasiconformal. 
Moreover, each forward orbit under $f_b$ passes through these annuli at most one
time. Thus there exists an invariant Beltrami differential invariant under $f_b$ which depends smoothly on the choice of $\beta_2$ and
vanishes when $\beta_2 = \beta_1$. Conjugating $f_b$ by a quasiconformal homeomorphism doubly tangent to the identity at infinity we
obtain a monic centered polynomial $g_b=g_{b(\beta_1,\beta_2)}$ which depends continuously on $\beta_2$ and $g_{b(\beta_1,\beta_1)}=f_1$.
It is not difficult to check that $g_b \in \cR(\lanot)$.

Now, given $\beta_2 \in U_\beta$ let $b=b(\beta_1,\beta_2)$ and observe that, by construction,
the element of $B(\hat{T})$ given by $\Phi \circ \chi(g_b)$
is conformally conjugate to $\beta_2$ in $|\hat{T}| \times \Delta(r_1)$. This conformal conjugacy extends, after taking successive 
preimages, to $|\hat{T}| \times \Delta$. Moreover, the conjugacy depends differentiably on $\beta_2$ and it is the identity when $\beta_2=\beta_1$.
Since the automorphism group of $B(\hat{T})$ is discrete, it follows that $\beta_2 = \Phi \circ \chi(g_b)$ for all $\beta_2 \in U_\beta$. 
Thus we have proven that $U_\beta$ has the desired property and we
 conclude that $\Phi'$ is onto.

 By the injectivity of $\Phi'=\Phi \circ \chi$, the map
 $g_{b(\beta_1,\beta_2)}$ does not depend on the choice of
 $U_\beta$ and $\beta_1$. Therefore we have obtained a well-defined map
 $\Psi:\beta \mapsto g_{b(\beta_1,\beta)}$ having the desired property.
\end{proof}

Let $\cH' = \Psi(B(\hat{T})) = \chi^{-1}(\cH)$.
\begin{lem}
 \label{lem-analytic-straightening}
For all $f \in \cH'$ there exists an open set $\cU_f \subset \poly(d)$ and an analytic family
 of polynomial-like mappings $((g^{\ell_v}:U_{g,v}' \to U_{g,\sigma(v)})_{v \in
 |T(\lanot)|})_{g \in \cU_f}$ over $T(\lanot)$, parameterized by
 $\cU_f$ such that the straightening map $\tilde{\chi}:\cU_f \to
 \cH$ is an extension of $\chi: \cH'\cap \cU_f \to \cH$ and $\tilde{\chi}$ is (complex) analytic.
\end{lem}

\begin{proof}
 By Lemma~\ref{lem-sufficient renor}, there exists a neighborhood $\cU$
 of $f$ in $\poly(d)$ and an analytic family of polynomial-like mappings 
 $F_g=((g^{\ell_v}:U_{g,v}' \to U_{g,\sigma(v)})_{v \in |T(\lanot)|})_{g \in
 \cU}$ over $T(\lanot)$.
Shrinking $\cU$, if necessary, we may assume that $F_g$ is hyperbolic, for all $g \in \cU$.
In particular, $K(F_g)$ is fiberwise connected for all $g \in \cU$.
Since the rays involved in an  external marking of $f$ (determined by an a priori chosen internal marking system $\alpha$ for $\lanot$) move continuously,
we may externally mark the whole analytic family of polynomial-like maps  so that the external marking of $F_g$  coincides
with that determined by $\alpha$ for all $g \in \cU \cap \cR(\lanot)$. 
Consider the associated straightening map $\tilde{\chi}:\cU \to \poly(T(\lanot))$.
 By the uniqueness part of the straightening theorem, $\tilde{\chi}$ is
 an extension of $\chi:\cR(\lanot)\cap \cU \to \cC(T(\lanot))$.
 Since the straightening map of an analytic family of hyperbolic polynomial-like map is complex
 analytic by \cite[page~313,
 Corollary~1]{Douady-Hubbard-poly-like}, it follows that $\tilde{\chi}$ is complex analytic.
Moreover, the straightening of a hyperbolic polynomial-like map is hyperbolic. Therefore,
$\tilde{\chi}(U) \subset \cH$, by the established continuity of $\tilde{\chi}$.
\end{proof}

Combining  Lemma~\ref{lem-qc-surgery} and
Lemma~\ref{lem-analytic-straightening}, we have the following:
\begin{cor}
 \label{cor-diff-submfd}
 $\cH' = \Psi(B(\hat{T}))$ is a (real) differentiable submanifold of
 $\poly(d)$ of (real) dimension $2d'-2$ where $d'$ is the total degree of $\hat{T}$.
\end{cor}

Thus, in order to prove  Theorem~\ref{introthm-onto hyperbolic} it remains 
to establish that $\cH'$ is a complex submanifold of $\poly(d)$.

\subsection{Proof of Theorem~\ref{introthm-onto hyperbolic}}

By the above lemma, the image of $\chi=\chi_{\lanot}$ contains $\hyp(\cC(T(\lanot)))$.
Consider a hyperbolic component $\cH$ of $\cC(T(\lanot))$ and let 
$\cH' = \chi^{-1}(\cH)$.

 Let $C_1, \dots, C_k$ be the critical classes of $\lanot$ and $C=C_1 \cup \cdots \cup C_k$.
Consider $f \in \cH'$ and an open connected neighborhood $\cU$ of $f$ in $\poly(d)$.
Take $\cU$ sufficiently small, so that for all $g \in \cU$, every cycle in the forward orbit of the landing point of the external rays with arguments in
$C$ are repelling. 
We may also assume that  attracting cycles of $f$ can be continued analytically in $\cU$ and that, for all $g \in \cU$, 
the  critical points of $g$ close to a critical point $\omega$ in the Fatou set of $f$,
belong to the basin of one of these attracting cycles. 
Note that the above implies that the landing point of the rays with arguments in $C$ depend continuously on $g \in \cU$.
Denote by $c_1, \dots, c_k$ the landing points, in the dynamical plane of $f$, of the external rays with arguments in
$C_1, \dots, C_k$, respectively.
Consider the analytic subset $X$ of $\poly(d)$  formed by all 
polynomials $g \in \cU$ with critical points 
$c_1(g), \dots, c_k(g)$, close to $c_1, \dots, c_k$,  so that the local degree of $g$ at  $c_j(g)$ is at least $\delta(C_j)$ and, for a 
generic $g$ (Zariski dense in $X$)  the local degree is exactly $\delta(C_j)$, for all $j$.
Let $S \subset X$ be the smaller analytic subset formed by all $g \in X$ such that
$g^n(c_i(g))=g^m(c_j(g))$ if  $d^n C_i = d^m C_j$ for some $n,m \geq 0$ and $1 \leq i \le j \le k$.
Hence, for all $g \in S$, each critical point of $g$ is either eventually periodic or contained in
the basin of an attracting cycle of $g$. Therefore, for all $g \in S$, every Julia  cycle of $g$ is repelling.
Note that $\cC(\lanot) \cap \cU \subset S$. Shrinking $\cU$, if necessary, we have that every irreducible component
of $S$ contains $f$.

We will show that a neighborhood $V$ of $f$ in $S$ is contained in $\cH'$. In particular, we will obtain that
$V$ is a neighborhood of $f$ in $\cC(\lanot)$ and
therefore conclude that $\cC(\lanot)$ is an analytic subset of $\poly(d)$ near $f$. 
Consider a  complex manifold
$S'$ and  a surjective analytic map $\pi:S' \to S$ (i.e. a resolution of the  singularities of $S$) such that
every component of $S'$ contains at least one element of $\pi^{-1}(f)$.
Via $S' \ni x \mapsto \pi(x) \in S$, we may regard $S'$ as a holomorphic family 
of polynomials. 
According to Ma\~{n}\'e-Sad-Sullivan~\cite{MSS}, for all $x \in S'$
close to $\pi^{-1}(f)$, the corresponding
polynomial $\pi(x)$ is structurally stable in a neighborhood of its Julia set, since all the Julia cycles of $g$ are repelling for 
all $g$ in a small neighborhood $V$ of $f$ in $S$.
By continuity of the landing point of (eventually) periodic external rays for structurally stable maps,
the rational lamination of all $g \in V$ coincide. Moreover, 
since $f \in \cR(\lanot)$, by Lemma~\ref{lem-sufficient renor},  we conclude that $V \subset \cR(\lanot)$.
Therefore, applying \cite[page~313, Corollary~1]{Douady-Hubbard-poly-like} we obtain that $\chi \circ \pi$ is complex analytic in $\pi^{-1}(V)$.
In particular, $\chi$ is continuous in $V$. It follows that $\chi(V) \subset \cH$.

Since $V$ is a neighborhood of $f$ both in $\cR(\lanot)$ and in $\cH'$,
by Corollary~\ref{cor-diff-submfd}, we conclude that $V$ (hence $S$) is smooth at $f$.
It follows that $\cH'$ is a complex submanifold of complex dimension $d'-1$ and
$\chi:\cH' \to \cH$ is a biholomorphism, where $d'$ is the total
degree of $\hat{T}$, which is equal to that of $\redschemanot$.
\qed


\section{Compactness}
\label{sec-compact}

The aim of  this section is to establish the following
compactness result which generalizes \cite[Theorem~1.3]{Inou-lim}.

\begin{thm}
\label{thm-compact}
  If $\lanot$ is a primitive invariant rational lamination, then $\cC(\lanot) = \cR(\lanot)$ and this set is compact.
\end{thm}

%

The reader may find the proofs of  
Theorem~\ref{thm-compact} and Theorem~\ref{introthm-compactness} 
after the statements and proofs of the following two lemmas.

Primitive rational laminations are rich in non-trivial classes:

\begin{lem}
\label{lem-finite-support}
  Let $\lanot$ be a non-trivial invariant rational lamination such that the
  support of $\lanot$ is contained in the $m_d$-grand orbit of a
  finite set. Then $\lanot$ is not primitive.
\end{lem}

\begin{proof}
 By passing to an iterate, we may assume that the support of $\lanot$
 is contained in the grand orbit of the finite set $F= \{ j/(d-1); j=
 0, \dots, d-2 \}$. That is, $t \in F$ if and only if $m_d (t) = t$.
 We proceed by contradiction and assume that $\lanot$ is primitive.

 Let $A_1, \dots, A_n$ be a complete list (without repetitions) of the
 $\lanot$-classes contained in $F$.  Let $f$ be a polynomial without
 neutral fixed points such that $\lambda_f = \lanot$. 
Let  $m$ be the number of 
non-trivial classes in $F$ and observe that $n+m \leq d-1$.

The union of the closure of
 the external rays with arguments in the set $F$ (which has cardinality $d-1$) cuts the complex plane into
$d-n$ regions $U_1, \dots, U_{d-n}$, since  
$$d -n = 1 + (d-1) -n =  1 + \sum^n_{j=1} (|A_j|-1).$$
According to Goldberg and Milnor \cite{Goldberg-Milnor}, for all $j$, there exists exactly 
one fixed point $z_j$ in the region $ U_j$.  These fixed points are
not the landing point of rational rays since the lamination is
supported in the grand orbit of $F$.  Therefore, $z_j$ is an
attracting fixed point, for all $j$. For each $j =1, \dots, d-n,$ call
$L_j$ the fixed infinite $\lanot$-unlinked class such that $z_j \in
K_f (L_j)$.  By~Proposition~\ref{prop-unlinked
  class}~\ref{item-prop-unlinked class-complement}, it follows that
the boundary of $L_j$ intersects a non-trivial class contained in $F$,
say $A(L_j)$.  
Hence, for some $j \neq k$, we have that $A(L_j) =
A(L_k)$. For otherwise, $d-n \le m$, but $n +m \leq d -1$. Therefore, $\lanot$ is not primitive which gives us the desired contradiction.
\end{proof}

\begin{lem}
\label{lem-perfect-lamination}
  Consider a non-trivial primitive $d$-invariant rational lamination $\lanot$.
  Let $E$ be the $m_d$-grand orbit of a finite set and $\lambda = \lanot \cap ((\QS \setminus E) \times (\QS \setminus E))$.
Then $\lanot$ is the smallest closed equivalence relation in $\QS$ containing $\lambda$.
\end{lem}

\begin{proof}
  Let $\lambda'$ be the smallest closed equivalence relation in $\QS$
  containing $\lambda$. It follows that $\lambda'$ is a
  $d$-invariant rational lamination and $\lambda' \subset
  \lanot$. 

Let $L_0$ be an infinite periodic $\lambda'$-unlinked class, say of
period $p$. We claim that $L_0$ is a $\lanot$-unlinked class.  For
this, consider $\pi : \overline{L_0} \rightarrow \S$ such that $\pi
\circ m_d^p = m_{d'} \circ \pi$ for some $d' \ge 2$.  Let $\pi_*
\lanot$ be the equivalence relation in $\QS$ that identifies $\theta$
and $\thetap$ if and only if there exists $t \in \pi^{-1} (\theta)$
and $s \in \pi^{-1} (\thetap)$ such that $s$ and $t$ are 
$\lanot$-equivalent arguments. 
Observe that $\pi_* \lanot$ is a $d'$-invariant rational
lamination with  support contained in $\pi (E \cap \overline{L_0})$. By the previous
lemma, $\pi_* \lanot$ is trivial or it is not primitive. The latter
alternative is impossible, since taking the preimage of appropriate
$\pi_* \lanot$-unlinked classes we would have that $\lanot$ is not
primitive. Therefore, $\pi_* \lanot$ is trivial. Hence, $L_0$ is a
$\lanot$-unlinked class.

Now let $L_1$ be a strictly preperiodic  infinite $\lambda'$-unlinked class. We also claim that $L_1$ is a $\lanot$-class.
In fact, let $\ell \ge 1$ be such that $m_d^\ell (L_1)$ is a periodic $\lanot$-unlinked class $L_0$ as in the previous paragraph. 
Let $\pi : \overline{L_0} \rightarrow \S$ be as above and consider $\hat\pi: \overline{L_1} \rightarrow \S$
such that $\pi \circ m_{d}^\ell = m_{d'} \circ \hat\pi$ for some $d' \ge 1$. Then $m_{d' *} (\hat\pi_* \lanot)$ is the trivial rational lamination.
In particular, every non-trivial $\hat\pi_* \lanot$-class is mapped, under $m_{d'}$ onto a singleton. Denote by $\Theta_1, \dots, \Theta_j$
these classes, note that the elements of each $\Theta_i$ differ by a multiple of $1/d'$. 
It follows unlinked classes with respect to the relation in $\R/\Z$ whose only non-trivial classes are $\Theta_1, \dots, \Theta_j$
have total length a multiple of $1/d'$ and there are at least $j+1$ of these unlinked classes (compare with Critical Portraits in~\cite{Poirier}).
Hence, $j \leq d'-1$.  Therefore, either $\hat\pi_* \lanot$ 
is trivial and 
$L_1$ is a $\lanot$-unlinked class, or $\lanot$ is not primitive.
Thus, $L_1$ is a $\lanot$-unlinked class.

Therefore, we have $\lambda'=\lanot$ by
 Lemma~\ref{lem-strict-inclusion}.
\end{proof}

With the previous lemma we may now establish our compactness result.

\begin{proof}[Proof of Theorem~\ref{thm-compact}]
Assume that $f \in \cC(\lanot)$.
In order to show that $\cC(\lanot) = \cR(\lanot)$ we must prove that $f$ is $\lanot$-renormalizable (i.e., $f \in \cR(\lanot)$).
More precisely, we must extract a polynomial-like map $g$ over $\redschema(\lanot)$ (see Definition~\ref{defn-renormalization}).
To extract a polynomial-like map from $f$ we will apply Lemma~\ref{lem-sufficient renor} to an appropriate Yoccoz puzzle.

We start by finding an appropriate combinatorial Yoccoz puzzle.
Let $F \subset \QS$ be the set formed by the arguments of external
rays of $f$ landing at parabolic or critical points. Let $m$ be the
minimum common multiple of the periods, under $m_d$, of the arguments
in the forward orbit of $F$. Denote by  $F' \subset \QS$ the set of all $m_d$-periodic arguments of period
dividing $m$ and, by $E'$ the grand orbit of $F'$.  From the previous lemma (Lemma~\ref{lem-perfect-lamination}),
there exist $\lanot$-classes $A_1, \dots, A_j$ contained in $\QS
\setminus E'$ that separate the critical elements of $\lanot$.  
The support of the depth $0$ puzzle will be 
$E_0 \subset \QS$, the set formed by all periodic arguments in the
forward orbit of $A_1 \cup \cdots \cup A_j$. 
More precisely, let $\lambda$ be the restriction of $\lanot$ to $E_0 \times E_0$.
For all $\ell \ge 1$, let $E_\ell = m^{-\ell}_d (E_0)$ and define $\Lambda_\ell$ as  the restriction of $\lanot$
to $E_\ell \times E_\ell$. It follows that $\Lambda = (\Lambda_\ell)$ is a combinatorial Yoccoz puzzle admissible for $\lanot$.
Note that, for  some $k \ge 0$, we have that $E_k \supset  A_1 \cup \cdots \cup A_j$.
In particular, from Definition~\ref{def-generator}, we have that $\Lambda$ is a generator for $\lanot$ with $k$ as separation depth.

Since $E_\ell \cap E' = \emptyset$,  for all
 $\theta \in E_\ell$, the landing point of the external
 ray of angle $\theta$ is a non-critical eventually repelling periodic point. 
By  Lemma~\ref{lem-sufficient renor},
we have  that $f \in \cR(\lanot)$.

\medskip
Now we show that $\cC(\lanot)$ is compact. Due to
  the compactness of $\cC(d)$, we only have to show that $\cC(\lanot)$ is
  closed in $\cC(d)$.  For this consider a sequence $(f_n)$ of polynomials in
  $\cC(\lanot)$ which converges to $f \in \cC(d)$.  Let $F \subset
  \QS$ be the arguments of the rays landing at parabolic periodic
  points of $f$ or at critical points of $f$.  Denote by $E$ the
  $m_d$-grand orbit of $F$.

We claim that 
 \[
 \lanot \cap ((\QS \setminus E) \times (\QS \setminus E))
 \subset 
 \lambda_f \cap ((\QS \setminus E) \times (\QS \setminus E)).
 \]
Assume that $\theta$ and $\thetap$ are two periodic elements in $\QS \setminus E$
which are $\lanot$-equivalent.
Since the landing points of the external rays
with arguments $\theta$ and $\thetap$ are repelling for $f$ and $\lambda_{f_n} \supset \lanot$,
 we have that $\theta$ and
$\thetap$ are $\lambda_f$-equivalent. A similar reasoning works for preperiodic arguments.

From the previous lemma and the fact that $\lambda_f$ is closed in
$\QS \times \QS$ we conclude that $\lambda_f \supset \lanot$. Hence,
$f \in \cC(\lanot)$.
\end{proof}

\begin{proof}[Proof of Theorem~\ref{introthm-compactness}]
Assume that
$\lanot$ is  hyperbolic post-critically finite but not primitive.
From the previously proven theorem, it is sufficient to show that: $\cR(\lanot)$ is not compact and $\cC(\lanot) \neq \cR(\lanot)$.

 Consider $f_0 \in \cR(\lanot)$ without neutral cycles and such that  
 $\lambda_{f_0}=\lanot$. 
 Since $\lanot$ is not primitive and hyperbolic, there exist two periodic Fatou
 components $\Omega_1$, $\Omega_2$ such that $\partial \Omega_1 \cap
 \partial \Omega_2 = \{x_0\}$ where $x_0$ is a repelling periodic point
 of period dividing the period of at least one of these  Fatou components, 
say $\Omega_1$.
 We may assume that $f_0$ has an attracting (not superattracting)
 periodic point $y_0$ in $\Omega_1$.
 By Ha\"{i}ssinsky's pinching theorem~\cite{Haissinsky-surgery}
 \cite{Haissinsky-pinching},
 there exist a continuous path of quasiconformal deformations
 of polynomials $(f_t)_{t \in [0,1)}$ such that
 \begin{enumerate}
  \item $f_t$ converges uniformly to a polynomial $f_1$, shrinking
	progressively some path $\gamma$ connecting $x_0$ and $y_0$.
	There exists a semiconjugacy $\varphi:\C \to \C$ from
	$f_0$ to $f_1$ such that $\varphi$ is a homeomorphism outside
	$\gamma$ and its preimages,
	and $\varphi$ sends $\gamma$ to a point, as well as 
	each component of the preimages of $\gamma$.
  \item The point $x_1=\varphi(\gamma)$ is a parabolic periodic point $f_1$ 
and its
	immediate basin is equal to $\varphi(\Omega_1 \setminus \Gamma)$
	where $\Gamma$ is the union of $\gamma$ and all the preimages.
 \end{enumerate}
 In particular, $f_0$ and $f_1$ are topologically conjugate on their
 Julia sets, by $\varphi$. 
 Therefore, $f_t \in \cR(\lanot)$ for $t \in [0,1)$ 
 but $f_1 \in \cC(\lanot)$ and $f_1 \not \in
 \cR(\lanot).$ 
  Otherwise, a $\lanot$-renormalization $g$ of $f_1$ would
 satisfy $\{x_1\} = K(g,v_1) \cap K(g,v_2)$
 where $v_i$ is the infinite $\lanot$-unlinked class such that
 $K(f_0,v_i)=\overline{\Omega_i}$.
 Therefore, $x_1$ would be repelling (\cite[Theorem~7.3]{McMullen} and
 \cite[Proposition~3.4]{Inou-ren}) which is a contradiction with (ii).
\end{proof}

\begin{rem}
 For hyperbolic post-critically finite invariant rational laminations $\lanot$, 
 the remaining problem is to characterize when $\cC(\lanot)$ is compact.
 Note that  $\cC(\lanot)$ is known to be compact in some cases and
 non-compact in some other cases.
 For example, 
 if $\lanot$ is of degree two, then $\cC(\lanot)$
 is always compact (it is a baby Mandelbrot set). 

 On the other hand, the following example shows that 
 $\cC(\lanot)$ is not compact in general.
\end{rem}

\begin{exam}
 Let $\lanot$ be the    rational lamination of
 $f_0(z)=z^3+3z/2$.
 The critical points $\pm i/\sqrt{2}$ of $f_0$ are fixed and
 the origin lies in the boundaries of their immediate basin of
 attraction. Since $f_0$ is real and $K(f_0) \cap \R=\{0\}$,
 the external rays of angles $0$ and $1/2$ are the positive and negative
 real line respectively, so $0$ and $\frac{1}{2}$ are $\lanot$-equivalent.

 \begin{figure}[tbh]
  \begin{center}
   \fbox{\includegraphics[width=6cm]{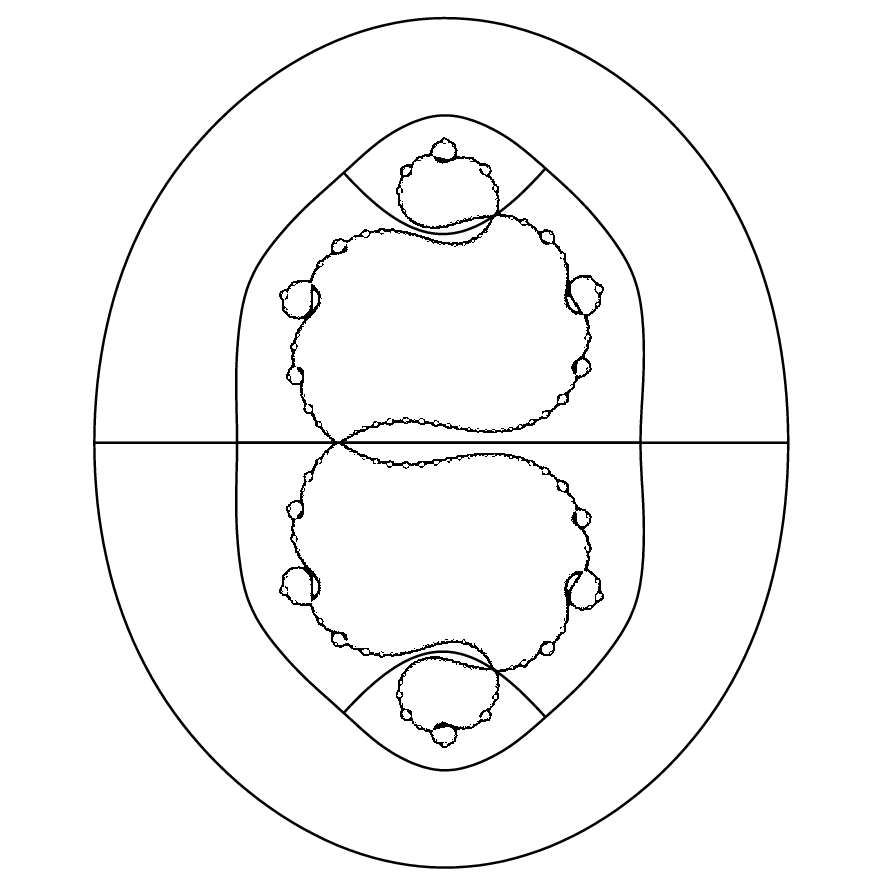}}
   \fbox{\includegraphics[width=6cm]{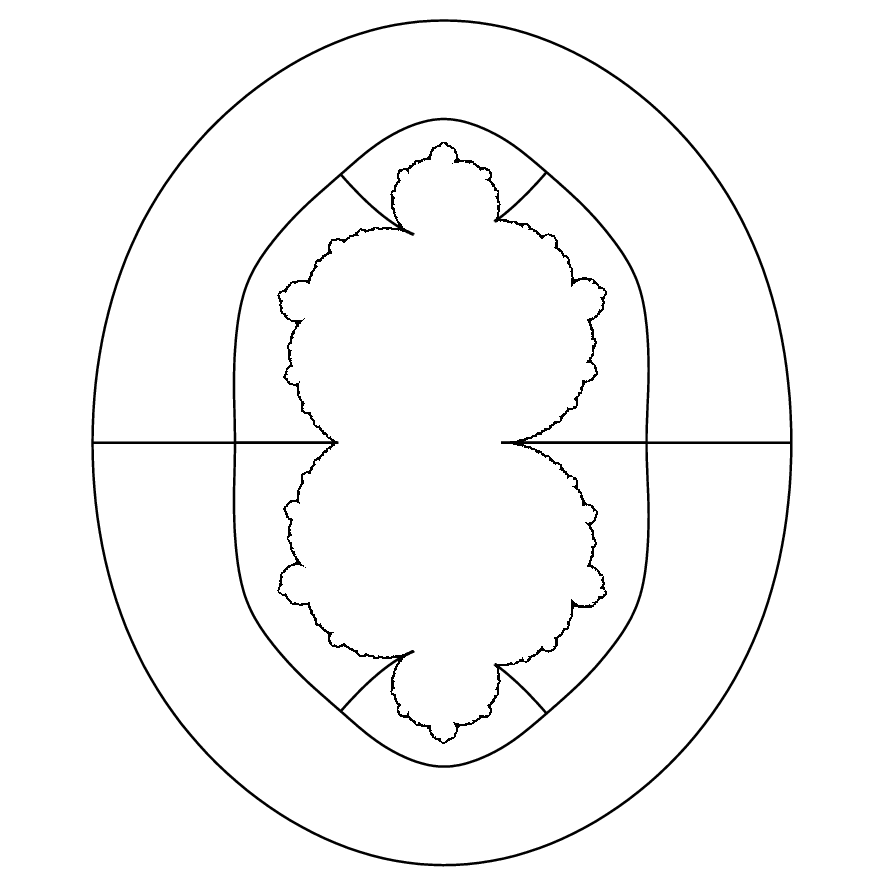}}
   \caption{The Julia sets of $f_{\mu_n}$ and its limit.}
   \label{fig-cubic1+1-collapse}
  \end{center}
 \end{figure}
For $\mu=\zeta+i\xi \ne 1$, let 
 \[
  f_\mu(z) = z^3 - \frac{2\xi}{\sqrt{2(1-\zeta)}}z^2 -
 \frac{1}{4}\left(2\zeta-6+\frac{2\xi^2}{\zeta-1}\right)z.
 \]
 (Note that for $\mu=0$, $f_\mu$ is equal to the previously
 defined $f_0$.)
Also, let $$\alpha=\frac{\xi}{\sqrt{2(1-\zeta)}}+\sqrt{\frac{1-\zeta}{2}}i.$$
 Then the fixed points of $f_\mu$ are $0$,
 $\alpha$ and 
 $\overline{\alpha}$. The multipliers of $\alpha$ and $\overline{\alpha}$ are 
 $\mu$ and $\overline{\mu}$ respectively.
 Hence if $|\mu|<1$, then $f_\mu \in \cC(\lanot)$.
 Fix $k>1/2$ and let $\mu_n=1-k/n^2+i/n$. Then
 $|\mu_n|<1$ for sufficiently large $n$, $\mu_n \to 1$, and
 \[
   f_{\mu_n}(z) \to g_k(z)=z^3 - \sqrt{\frac{2}{k}}z^2 +
 \left(1+\frac{1}{2k}\right)z.
 \]
 Clearly, $g_k$ is real and has a repelling fixed point at $0$ and a
 parabolic fixed point at $1/\sqrt{2k}>0$, hence $K(g_k) \cap \R$ is a
 closed interval and the external rays of angles $0$ and $1/2$ do not
 land at the same point, therefore $g_k \not \in \cC(\lanot)$ (see
 Figure~\ref{fig-cubic1+1-collapse}).
 In particular, $\cC(\lanot)$ is not compact.

 Another example is illustrated in Figure~\ref{fig-cubic2+2-collapse},
 which was proposed by Goldberg-Milnor \cite{Goldberg-Milnor}.
 \begin{figure}[tbh]
  \begin{center}
   \fbox{\includegraphics[width=6cm]{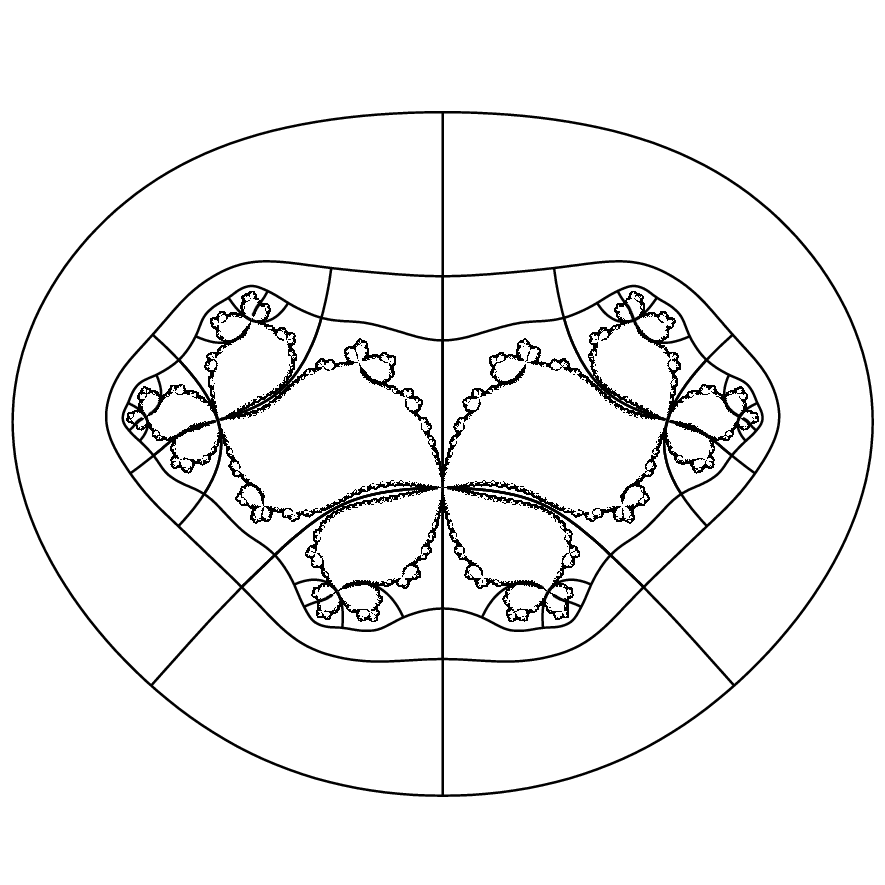}}
   \fbox{\includegraphics[width=6cm]{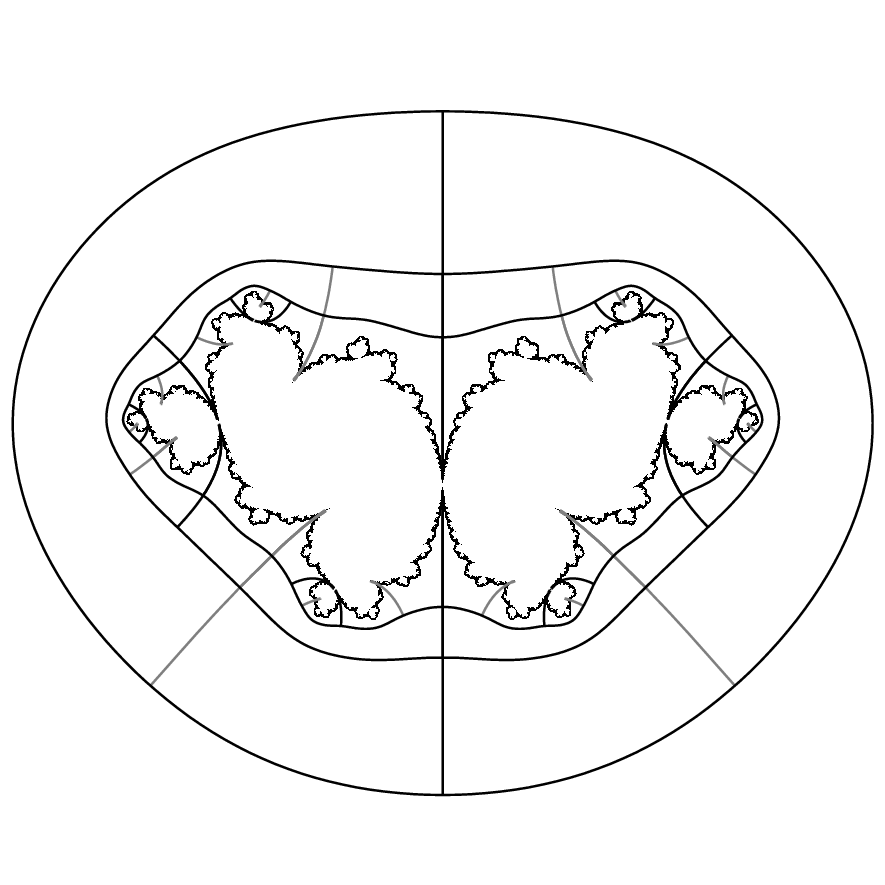}}
   \caption{Butterfly collapses to a period two parabolic cycle. See
   Goldberg-Milnor \cite{Goldberg-Milnor}.}
   \label{fig-cubic2+2-collapse}
  \end{center}
 \end{figure}
\end{exam}

\begin{exam}
\label{exam:last}
  A related example of a post-critically finite polynomial of degree $4$ such that  $\cC(f_0)$ is not compact arises from the presence of a strictly preperiodic critical point in $J(f_0)$. The polynomial $f_0$ has two critical points $\omega$ and $\omega'$, where $\omega$ has multiplicity $2$ . While $\omega$ is fixed, the critical point $\omega'$ maps in two iterations onto a fixed point $z_0$ in the boundary of the basin of $\omega$. The point $z_0$ is the landing point of the $0$ ray and $\omega'$ is the landing point of the rays $3/8$ and $5/8$. 
This is illustrated in Figure~\ref{fig:last} (left). Now $\cC(f_0)$ contains
a sequence of polynomials $f_n$ with two attracting fixed points converging
to a polynomial $g$ which has a parabolic fixed point as landing point of the $0$ ray. 
Figure~\ref{fig:last} (middle) illustrates the Julia set of $f_n$, for some value of $n$, and it suggests that
$\lambda(f_n) \supset \lambda (f_0)$ (i.e., $f_n \in \cC(f_0)$). 
For example, a Julia fixed point of $f_n$ is the landing point
of the rays corresponding to the $\lambda(f_n)$-class $\{0,1/3,2/3 \}$ which (strictly) contains
each of the $\lambda(f_0)$-fixed classes $\{0\}$ and $\{1/3,2/3\}$.  
As the illustration in Figure~\ref{fig:last} (right) of $J(g)$ suggests  we have that $g \notin \cC(f_0)$, since
the $3/8$ and $5/8$ rays land at distinct points. Thus, the $\lambda(f_0)$-class $\{3/8,5/8\}$ is not contained in
a $\lambda(g)$-class.
Note that the pre-fixed critical point of $g$ maps onto the repelling fixed point in 
the boundary of the parabolic basin which is the landing point of the $1/3$ and $2/3$ rays.

 \begin{figure}[tbh]
  \begin{center}
   \fbox{\includegraphics[width=3.9cm]{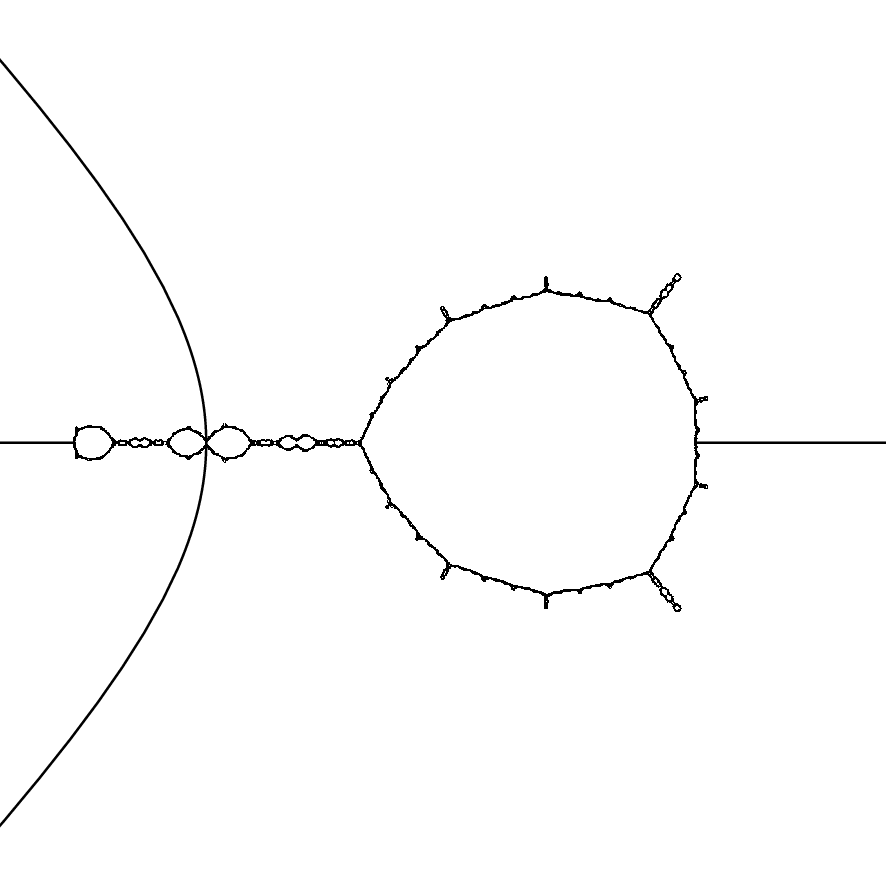}}
   \fbox{\includegraphics[width=3.9cm]{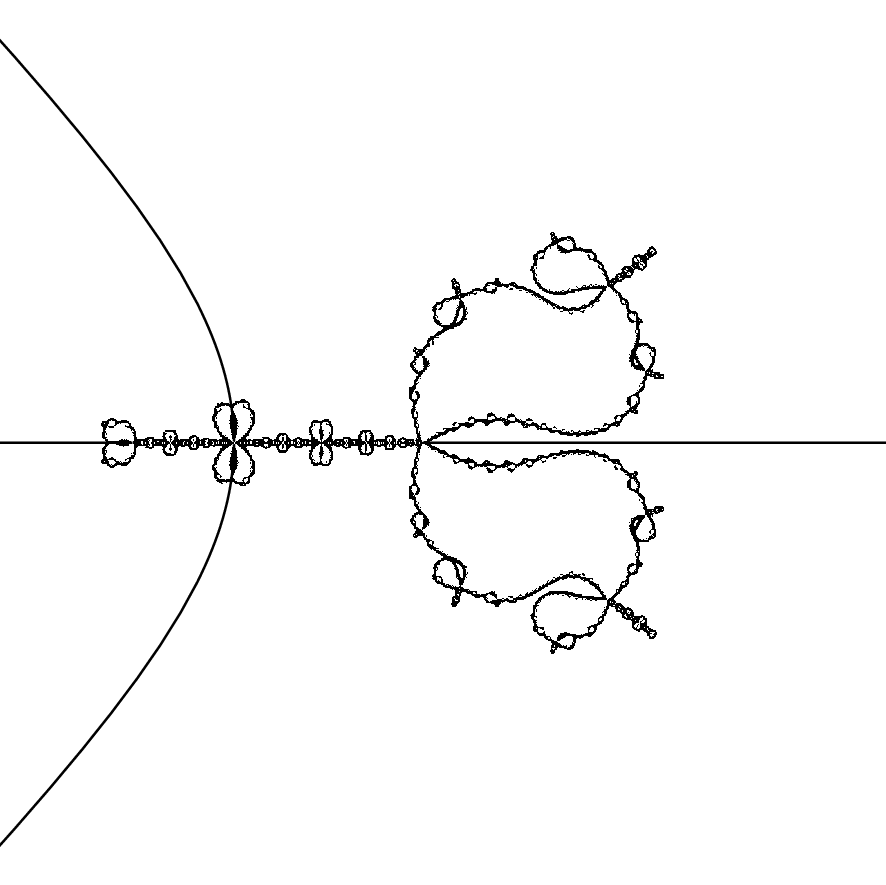}}
   \fbox{\includegraphics[width=3.9cm]{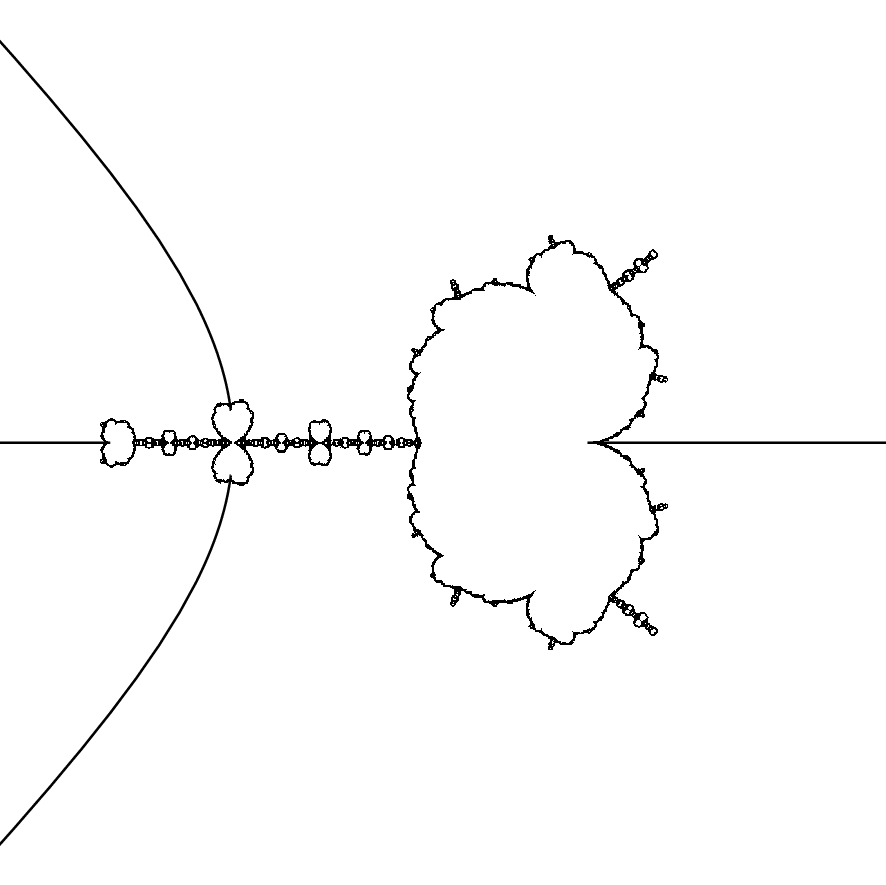}}
   \caption{Illustration of Example~\ref{exam:last}.}
   \label{fig:last}
  \end{center}
 \end{figure}
\end{exam}


\section{Surjectivity of straightening maps}
\label{sec-surjectivity}

\subsection{Preliminaries}
\label{subsec-surjetivity preliminaries} 
To prove our surjectivity results we will employ the compactness of primitive renormalizable polynomials and the injectivity of straightening proven
in previous sections together with the nice properties of straightening
of quadratic-like families, established by Douady and Hubbard
in~\cite{Douady-Hubbard-poly-like}. We summarize some of these
properties as follows (compare with ~\cite[Corollary 2 of Proposition~13 and Chapter IV]{Douady-Hubbard-poly-like}).

\begin{thm}
\label{thm-DH-cont}
  Let $X$ be a complex manifold and $\{f_\mu \}_{\mu \in X}$ be an analytic family of quadratic-like maps. Let
$$\cM_X = \{ \mu \in X; K(f_\mu) \mbox{ is connected } \},$$
and denote by 
$$\chi : \cM_X \to \cM$$
the corresponding straightening map. That is, $\chi(\mu) =c$ where $c \in \cM$ is such that $z \mapsto z^2 +c$ is hybrid equivalent to $f_\mu$.
Then the following statements hold:
\begin{itemize}
\item
$\chi$ is continuous.
\item
For all $c \in \cM$, we have that $\chi^{-1} (c)$ is a complex analytic subspace of $X$.
\item 
If $X$ has complex dimension one and $\mu_0 \in \cM_X$,
then there exists a neighborhood $V$ of $\mu_0$ in $\cM_X$ such that $\chi$ is constant on $V$ or
$\chi(V)$ contains a neighborhood of $\chi(\mu_0)$ in $\cM$.
\end{itemize}
\end{thm}

We will also need the following lemma on quasiconformal deformations:
\begin{lem}
 \label{lem-qc deform}
 Let $\lanot$ be a $d$-invariant rational lamination.
 Let $f \in \cR(\lanot)$ and $\straightening_{\lanot}(f)=P \in
 \cC(\redschemanot)$.
 Consider a quasiconformal deformation of $P$ on $K(P)$,
 i.e., there exists a quasiconformal map $\phi:\C \to \C$ whose complex
 dilatation is supported on $K(P)$ such that
 $\phi \circ P = P_1 \circ \phi$ for some $P_1 \in \cC(\redschemanot)$.

 Then there exists a quasiconformal deformation $f_1 \in \cR(\lanot)$
 of $f$ on $K(f)$ such that $\straightening_{\lanot}(f_1)=P_1$.
\end{lem}

\begin{proof}
 Let us denote by $\sigma_0$ the standard complex structure on $\C$ 
 and let $\sigma = \phi^*\sigma_0$.
 Define an almost complex structure $\sigma_1$ on $\C$ as follows:
 \[
  \sigma_1 = 
  \begin{cases}
   (\psi_v \circ f^n)^* \sigma &
   \mbox{on $f^{-n}(K_f(v))$ for some $v \in |\redschemanot|$ and }n \ge
   0, \\
   \sigma_0 & \mbox{otherwise,}
  \end{cases}
 \]
 where $\psi=(\psi_v)$ is a hybrid conjugacy from a
 $\lambda$-renormalization of $f$ to $P$.
 Then $\sigma_1$ is well-defined because $\sigma$ is $P$-invariant.
 Its dilatation is uniformly bounded
 because $\sigma_1$ is defined via  pullbacks of $\psi_v^*\sigma=(\phi \circ
 \psi_v)^*\sigma_0$ by holomorphic maps which preserve the maximal dilatation.
 Furthermore, $\sigma_1$ is $f$-invariant by construction.
 Therefore, by the measurable Riemann mapping theorem,
 $f:(\C,\sigma_1) \to (\C,\sigma_1)$ is conformally conjugate to a
 polynomial $f_1:\C \to \C$.
 Since this deformation does not change the complex structure of the
 basin of infinity, we have $\lambda_f=\lambda_{f_1}$, in particular
 we have $f_1 \in \cC(\lanot)$ and a
 $\lanot$-renormalization of $f$ gives rise to a $\lanot$-renormalization of
 $f_1$.
 More precisely, let $\phi_1:\C \to \C$ be the quasiconformal
 homeomorphism such that $\phi_1^*\sigma_0 = \sigma_1$ and $f_1 = \phi_1
 \circ f \circ \phi_1^{-1}$.
 Take a $\lanot$-renormalization $g=(f:U_v' \to U_{\sigma_{\lanot}(v)})$
 of $f$.
 Then $g_1=(f_1:\phi_1(U_v') \to \phi_1(U_{\sigma_{\lanot}(v)}))$
 is a $\lanot$-renormalization of $f_1$
 and it is hybrid equivalent to $P_1$ by $\phi \circ \psi \circ
 \phi_1^{-1}$.
\end{proof}

\subsection{Surjectivity of straightening: primitive disjoint type}
\begin{thm}
 \label{thm-disjoint-surjectivity}
 Let $\fnot \in \cC(d)$ be an internally angled primitive
hyperbolic post-critically finite  polynomial with rational lamination $\lanot$.
If $\fnot$ has exactly $d-1$ superattracting periodic orbits,
 then the corresponding straightening map $\chi_\lanot : \cR(\lanot) \to \cM^{d-1}$ is a homeomorphism. 
\end{thm}

\begin{proof}
Recall that the interior of the Mandelbrot set $\cM$ is dense in $\cM$.
Therefore, the interior of $\cM^{d-1}$ is also dense in $\cM^{d-1}$. 
The components of the interior of $\cM$ are either hyperbolic (maps with an attracting cycle) 
or queer. The closure of the union of the hyperbolic components is the complement of the queer components.
That is, $\partial \cM$ is contained in the closure of the union of the hyperbolic components.

Let $\lanot$ be an internally angled hyperbolic primitive rational lamination of degree $d$ with reduced mapping schema of disjoint type
(i.e. $|T(\lanot)|$ consists of $d-1$ vertices, each fixed by the schema map).
 By Lemma~\ref{lem-sufficient renor}, for every $f \in \cR(\lanot)$ there exist a
 neighborhood $\cV$ of $f$ in $\poly(d)$ and an analytic family of
 polynomial-like maps over $\redschemanot$ whose straightening map,
 which is continuous by Theorem~\ref{thm-DH-cont}, is an extension of
 $\chi_{\lanot}$. 
 Therefore, $\chi_{\lanot}$ is continuous.

Since $\chi_{\lanot}$ is continuous and $\cR(\lanot)$ is compact (Theorem~\ref{introthm-compactness}), the
image of $\chi_{\lanot}$ is closed in $\cM^{d-1}$. Taking into account that the interior of $\cM^{d-1}$ is dense in
$\cM^{d-1}$, to prove that $\chi_\lanot$ is surjective it is sufficient to prove that its image contains
every connected component $S_0$ of the interior of $\cM^{d-1}$.
Given such a connected component $S_0$, there exist $d-1$ components $W_1, \dots, W_{d-1}$ of the interior of $\cM$ 
such that 
$S_0= W_1 \times \cdots \times W_{d-1}$.
Since hyperbolic maps are in the image of $\chi_\lanot$ by
 Theorem~\ref{introthm-onto hyperbolic}, we have that $\left( \partial
 \cM \right)^{d-1}$ is also contained in the image.
In particular, $\partial W_1 \times \cdots \times \partial W_{d-1}  \subset \chi_\lanot (\cR(\lanot))$.
Now, for  $n \in \{ 1, \dots, d-1 \}$, let
$$S_n = \partial W_1 \times \cdots \partial W_n \times W_{n+1} \times \cdots \times W_{d-1}.$$
Thus, $S_{d-1} \subset \chi_\lanot (\cR(\lanot))$ and
to prove that $S_0$ is contained in  $\chi_\lanot (\cR(\lanot))$  it is enough to establish the following.

\noindent
{\bf Claim.} For all $1 \le n \le d-1$, if $S_n \subset \chi_\lanot (\cR(\lanot))$, then $S_{n-1} \subset \chi_\lanot (\cR(\lanot))$.

\noindent
{\it Proof of the claim.}  Let $\mathbf{c}=(c_1, \dots, c_{d-1}) \in
S_n$. Denote by $\pi : \C^{d-1} \to \C^{d-2}$ the projection which
forgets the $n$th coordinate and by $\pi_n :\C^{d-1} \to \C$ the projection onto the $n$th coordinate. 

Consider the set $Y=(\pi \circ \chi_\lanot)^{-1}(\pi(\mathbf{c}))$. Recall that $\lambda_0$-renormalizable polynomials 
$g$ are such that its renormalization  consists of $d-1$ quadratic like maps $g_i$ for $i=1, \dots, d-1$. 
Then, $Y$ consists of
all $\lanot$-renormalizable polynomials $g$
such that, for all $i \neq n$,
the quadratic like map $g_i$ is hybrid equivalent to the prescribed quadratic polynomial $z^2 + c_i$.   
We claim that $Y$ is a complex analytic space $Y$.
In fact, take any $g \in Y$, then by~Lemma~\ref{lem-sufficient renor} 
and the proof Theorem~\ref{thm-compact}, there exists a neighborhood $\cV$ of
$g$ in $\poly(d)$ that parametrizes an analytic family of polynomial like maps over $T(\lanot)$.
That is, $d-1$ families of quadratic like maps. Since $Y$ corresponds, locally around $g$, to 
the intersection of $d-2$ complex analytic sets, according to Theorem~\ref{thm-DH-cont}, it follows that
$Y$ is a complex analytic set.

Since $\cR(\lanot)$ is compact (Theorem~\ref{introthm-compactness}) 
and the uncountable set $\partial W_n$ is contained in $\pi_n \circ \chi_\lanot (\cR(\lanot) \cap Y)$, there
exists $f \in \cR(\lanot)$ contained in a component of dimension at
 least one of $Y$ such that  $\pi_n \circ \chi_\lanot(f) \in \partial W_n$.
 Moreover, we may assume that $f$ is a smooth point of $Y$ and not isolated in $Y \cap \cR(\lanot)$.
Hence there exists a one dimensional submanifold of $Y$ containing $f$ and some other point of $Y \cap \cR(\lanot)$.
Theorem~\ref{introthm-injectivity} implies that
 $\pi_n \circ \chi_\lanot$ is not constant along this one dimensional submanifold of $Y$. Therefore
$\pi_n \circ \chi_\lanot (\cR(\lanot) \cap Y)$ contains a neighborhood of $c_n$ in $\cM$, by Theorem~\ref{thm-DH-cont}. 
From Lemma~\ref{lem-qc deform}, we conclude that  $\pi_n \circ \chi_\lanot (\cR(\lanot) \cap Y)$ contains $W_n$.
Hence, for all $\mathbf{c}=(c_1, \dots, c_d) \in S_n$
we have that $\{c_1 \} \times \cdots \times \{ c_{n-1} \} \times W_n
\times \{ c_{n+1} \} \times \cdots \times \{c_{d-1} \} \subset
\chi_\lanot (\cR(\lanot))$ and the claim follows.
\end{proof}

\subsection{Surjectivity of straightening: cubic primitive capture type}

\begin{thm}
\label{thm-capture-surjectivity}
Let $\fnot \in \cC(3)$ be an internally angled primitive
hyperbolic post-critically finite  polynomial with rational lamination $\lanot$ and reduced mapping schema of capture type.
 Then the associated straightening map
$\straightening_\lanot:\cR(\lanot) \to \cMK$ is a
      bijection. 
\end{thm}

The proof of Theorem~\ref{thm-capture-surjectivity} is similar to the
disjoint case. However we have to be careful since the straightening maps 
involved are discontinuous
\cite{Inou-discont}. Nevertheless, we use
the fact that  these maps are continuous along carefully chosen 
sequences.

\begin{lem}
 \label{lem-deformation}
 Let $f_n:U_n' \to U_n$ be a sequence of quadratic-like maps with
 connected Julia sets. Assume that the following hold:
 \begin{enumerate}
  \item $f_n$ converges to a quadratic-like map $f:U' \to U$.
	More precisely, $f_n \to f$ uniformly on some neighborhood of
	$K(f;U',U)$ as $n \to \infty$;
  \item $f_n$ is hybrid equivalent to $g_n \in \cM$.
  \item $g_n \to g\in \cM$;
  \item \label{item-approx fatou comp}
	For any $z \in \Int K(g)$, $z \in \Int K(g_n)$ for all sufficiently
	large $n$.
 \end{enumerate}
 Then, for each $n$, we can choose a hybrid conjugacy $\psi_n$, between $f_n$ and $g_n$,
 such that $\psi_n$ converges to a hybrid conjugacy $\psi$ between $f$
 and $g$.

 In particular, let $\lanot$ be a hyperbolic 3-invariant rational
 lamination of capture type.
 Let $f_n \to f$ be a convergent sequence in $\cR(\lanot)$.
 If the quadratic renormalization $f_n^p:U_n' \to U_n$ of $f_n$
 satisfies the above assumption, then $\straightening_\lanot(f_n) \to
 \straightening_\lanot(f)$.
\end{lem}

\begin{rem}
 \label{rem-approx fatou comp}
 The assumption \ref{item-approx fatou comp} is equivalent to  $J(g_n)
 \to J(g)$ in the Hausdorff topology. Furthermore,  \ref{item-approx fatou comp} holds if one of
 the following hold:
 \begin{enumerate}
  \item $J(g)=K(g)$, i.e., $\Int K(g)$ is empty;
  \item $g$ lies in the boundary of a hyperbolic component $\cW$
	of $\cM$ and $g_n \in \cW \to g$
	non-tangentially (see \cite{McMullen-HD2} for the case $g$
	parabolic and see \cite{Yoccoz} for the case $g$ has a Siegel
	disk).
 \end{enumerate}
 In particular, for any quadratic polynomial $g(z)=z^2+c$ with $c \in
 \partial M$, there exists a sequence $g_n \to g$ such that $g_n$ is
 hyperbolic and satisfies the assumption \ref{item-approx fatou comp}.
\end{rem}

\begin{proof}
After affine conjugacy, we may assume that $f_n$ and $f$ are normalized to have Taylor series of the form $c+z^2 + \cdots$, around the origin. 
Passing to the space of quadratic-like germs $\mathcal{QG}$, introduced by Lyubich~\cite{Lyubich-ql-space}, we
have that  the classes $[f_n] \in \mathcal{QG}$ converge to $[f] \in \mathcal{QG}$.  
Hence, $f$ is hybrid equivalent to $g$, by continuity of straightening in the space of quadratic-like germs with connected filled Julia sets (e.g. see the stronger result~\cite[Theorem~4.13]{Lyubich-ql-space}).

 Since $f_n$ converges to $f$,
 we may assume diffeomorphisms between fundamental annuli in the proof of
 the straightening theorem for $f_n$ also converge to that of $f$.
 In particular, they are uniformly quasiconformal.
 Then it follows by construction that the hybrid
 conjugacies $\psi_n$ are also uniformly quasiconformal.
 (compare the tubing construction of hybrid conjugacies
 \cite{Douady-Hubbard-poly-like}). 
 The limits of subsequences of $\psi_n$ agree on the complement of the filled Julia set of $f$. Thus to show that $\psi_n$ converges, it is sufficient to show that
if $\psi$ is a limit of a subsequence of $\psi_n$, then $\psi$ is a hybrid conjugacy.
 Therefore, by passing to a subsequence, we assume that $\psi_n$ converges
 to some quasiconformal map $\psi$, which conjugate $f$ to $g$ and proceed to show that in fact is a hybrid conjugacy.

 If $g$ is structurally stable (equivalently, $g \in \Int
 \cM$), then $f$ also is structurally stable in the space of
 quadratic-like germs \cite{Lyubich-ql-space}.
 Thus the proof that hybrid conjugacies that arise from the
 tubing construction vary continuously in the stability locus
 \cite[Proposition~12]{Douady-Hubbard-poly-like} can be applied to the
 convergent sequence $f_n \to f$.
 Therefore $\psi$ is a hybrid conjugacy in this case.

 Now we assume $g \in \partial \cM$.
 Let $\mu=\frac{\bar{\partial}\psi^{-1}}{\partial\psi^{-1}}$ be the
 complex dilatation of $\psi^{-1}$.
 Then, by the equation $\psi \circ f = g \circ
 \psi$, we have $g^*\mu=\mu$.
 Since $g$ does not carry an invariant line field on its Julia set
 (otherwise, $g$ must lie in a queer component), we may assume that
 $\mu$ vanishes  on $J(g)$.
 Furthermore, since $\psi_n$ is holomorphic in the interior of $K(f_n)$, 
 $\mu$ also vanishes on  $\Int K(g)$ by \ref{item-approx fatou comp}.
 Therefore, $\psi$ is a hybrid conjugacy between $f$ and $g$.

 For $f_n \to f \in \cR(\lanot)$, let $\omega_n$ be the
 captured critical point for $f_n$ and $k$ be the capture time.
 Then we have $\straightening_\lanot(f_n)=(g_n,x_n)$ where
 $x_n=\psi_n(f_n^k(\omega_n))$. 
 If $f_n:U_n'\to U_n$ satisfy the assumption,
 then, after passing to the limit, we have $x = \psi(f^k(\omega))$, where
 $\omega$ is the captured critical point for $f$.
 Therefore, $\straightening_\lanot(f)=(g,x)=\lim\straightening_\lanot(f_n)$.
\end{proof}

\begin{proof}[Proof of Theorem~\ref{thm-capture-surjectivity}] 
To simplify notation, we identify the map $h(z) = z^2 + c$ with $c 
\in \C$.
That is, straightening takes values in the set:
 \[
  \cMK= \{(h,x);\ h(z)=z^2+c \in \cM,\ x \in K(h)\}
 \]
 
Observe that if $h$ is hyperbolic and $x \in \Int K(h)$, then
the pair $(h, x)$ represents a hyperbolic dynamical system over the reduced 
schema of $\lanot$. By Theorem~\ref{introthm-onto hyperbolic}, 
$(h,x)$ is in the image of $\chi_\lanot$.
 
Now we consider $(h_0, x_0) \in \cMK$.
To prove that  $(h_0, x_0) \in \chi_{\lanot}(\cR(\lanot))$ we consider two cases
according to whether $h_0$ is in a queer component or not.

 \smallskip
 \noindent
 {\bfseries Case I:} $h_0$ does not lie in any queer component.

 Take sequences $h_n \to h_0$ and $x_n \to x_0$ satisfying
 \ref{item-approx fatou comp} in Lemma~\ref{lem-deformation}
 such that $h_n$ is hyperbolic and $x_n \in \Int K(h_n)$.
 Since $(h_n,x_n)$ corresponds to a hyperbolic polynomial over
 $T(\lanot)$,
 there exists $f_n \in \cR(\lanot)$ satisfying
 $\straightening_\lanot(f_n)=(h_n,x_n)$. Since $\cR(\lanot)$ is compact,
 we may assume $f=\lim f_n \in \cR(\lanot)$ exists.
 Then, by
 Lemma~\ref{lem-deformation}, we have $\straightening_\lanot(f)=(h_0,x_0)$.

\smallskip
When $h_0$ is in a queer component we further subdivide into two cases according
to whether the captured critical forward orbit is finite or not.

 \smallskip
 \noindent
 {\bfseries Case II.a:} 
$h_0 \in \cW$ for some queer component
 $\cW$ and  $x_0$ is periodic for $h_0$.
 
 In this case, we have $x_0$ is a repelling period point of period, say $j$.
Consider the continuous map $x: \overline{\cW} \to \C$ such that $x(h_0) = x_0$ and 
$x(h)$ is a period $j$ periodic point of $h$ (this map is unique because for all $h \in \overline{\cW}$ we have that
all cycles of $h$ are repelling). 

For $f \in \cR(\lanot)$ we denote by $(g_{f,i} : U'_{f,i} \to  U_{f,0})_{i=0,1}$ its renormalization where
$g_{f,0}$ is a quadratic-like map obtained from an appropriate restriction of $f^p$, for some $p$ independent of $f$
and $g_{f,1}$ is a degree $2$ branched covering obtained from a restriction of $f^\ell$, for some $\ell$ independent of $f$.
 
Now consider the one dimensional complex analytic set $V \subset \poly(3)$ formed by all cubic polynomials $f$ that posses a critical point $\omega$ such that $f^{pj + \ell}(\omega) = f^\ell (\omega)$. 
By Case I and the definition of $V$, for all $h \in \partial \cW$ we have that 
$(h, x(h)) \in \chi_\lanot (\cR(\lanot)) \cap V$. 

We choose a smooth point $f_1$ of $V$ such that $\chi_\lanot(f_1) = (h_1 , x(h_1))$ with $h_1 \in \partial \cW$.
By Lemma~\ref{lem-sufficient renor} and Lemma~\ref{lem-hyperbolic thickening} there exists a neighborhood $\cU$ of $f_1$ in $\poly(3)$ 
that parametrizes an analytic family
of polynomial like maps $g_f= (g_{f,i})$ over $T(\lanot)$ such that if $g_f$ has fiberwise connected filled Julia set $K(g_f)$, then $f$ lies in $\cR(\lanot)$ and $\chi_\lanot(f)$ is the straightening of  $g_f$. 
After shrinking
$\cU$, if necessary, for all  $f \in \cU \cap V$, the critical value of $g_{f,1}$ is a periodic point
 in the filled Julia set $K(g_{f,0})$ of the quadratic-like polynomial $g_{f,0}$.
Hence, $K(g_f)$ is fiberwise connected if and only if $K(g_{f,0})$ is connected, for all  $f \in \cU \cap V$.

The straightening of the quadratic-like map $g_{f,0}$ is not constant in $\cU \cap V$, otherwise, 
for all $f \in \cU \cap V$, we would have $\chi_\lanot(f) = (h_1,x)$ where $x$ is a period $j$ point of $h_1$,  
which would be a contradiction with Theorem~\ref{introthm-injectivity}.
Thus, for all $h$ in a neighborhood of $h_1$ in $\overline{\cW}$ there exists 
 $f \in \cU \cap V$, such that $\chi_\lanot (f) = (h, x)$ where $x$ has period $j$ under $h$.
By Theorem~\ref{thm-DH-cont}, $h \in \overline{\cW}$ depends continuously on $f$ and by Lemma~\ref{lem-deformation}, $x$ is also continuous on $f$,
it follows that $x=x(h)$ since the repelling periodic point $x_1$ of $h_1$ has a unique (analytic) continuous continuation.
We conclude that there exists $h_2 \in \cW$ such that $(h_2,x(h_2)) = \chi_\lanot(f_2)$ for some $f_2 \in \cR(\lanot)$.
Since $(h_0,x_0)$ is a quasiconformal deformation of $(h_2, x(h_2))$,  by Lemma~\ref{lem-qc deform} we have
that  $(h_0,x_0)$ lies in the image of $\chi_\lanot$.

 \smallskip
 \noindent
 {\bfseries Case II.b:} Assume that 
$h_0 \in \cW$ for some queer component
 $\cW$ and $x_0 \in J(h_0)(=K(h_0))$.

 Since periodic points are dense in $J(h_0)$,
 we can take a sequence $x_n \to x_0$ with $x_n$ periodic.
 Take $f_n \in \cR(\lanot)$ which satisfy
 $\straightening_\lanot(f_n)=(h_0,x_n)$. Since $\cR(\lanot)$ is
 compact, we may assume $f_n$ converges to some $f \in
 \cR(\lanot)$. 
 With the notation of case II.a, we may choose the quadratic-like maps $g_{f_n,0}: U_{f_n,0}' \to U_{f_n,0}$
converging to $g_{f,0}: U_{f,0}' \to U_{f,0}$. Since all  $g_{f_n,0}$ are hybrid equivalent to $h_0$ and its filled Julia set
$K(h_0)$ has empty interior, we may apply Lemma~\ref{lem-deformation} to conclude that
 $\straightening_\lanot(f)=(h_0,x_0)$.
\end{proof}

\subsection{Complex submanifolds and quasiconformality of straightening}
\label{sec-qc}
Here we prove some regularity properties of straightening maps of
primitive disjoint type and primitive cubic capture type. That is, we
finish the proof of theorems~\ref{introthm-disjoint-surjectivity} and
~\ref{introthm-capture-surjectivity}.

\subsubsection{Proof of Theorem~\ref{introthm-disjoint-surjectivity}}
In view of Theorem~\ref{thm-disjoint-surjectivity} the straightening map $\chi_\lanot$ is surjective.

To prove the existence of a complex submanifold,
we work in the space $\cQG$ of quadratic-like germs introduced by
\cite{Lyubich-ql-space}. 
Given $f \in \cR(\lanot)$, let $\cU$ be a neighborhood of $f$ in $\poly(d)$ such that
for all $g \in \cU$  and $1 \le i \le d-1$, we may extract a 
quadratic-like map  $$\iota_i(g) = [g^{\ell_i}:U_i' \to U_i] \in \cQG,$$
as in Lemma~\ref{lem-sufficient renor},  which is hybrid conjugate
to $z^2 + \chi_i(g)$, for all $g \in  \cU \cap \cR(\lanot)$. Recall that the analytic structure of $\cQG$ is
induced from complex Banach spaces of holomorphic functions defined on
domains equipped with the sup-norm
\cite[Section~4]{Lyubich-ql-space}.
Hence it follows that $\iota_i$ is analytic. 

For each $i=1,\dots,d-1$, consider analytic sets
\begin{align*}
 \cS_i=\cS_i(f) &= \{g \in \cU;\ \chi_i(g)=\chi_i(f)\},\\
 \cS_i^*=\cS_i^*(f) &= \{g \in \cU;\ 
 \chi_j(g)=\chi_j(f)\mbox{ for all }j \ne i\}.
\end{align*}
Then $\chi_i:\cS_i^* \to \poly(2)$ is continuous and maps $\cS_i^*
\cap \cR(\lanot)$ homeomorphically onto a neighborhood of $\chi_i(f)$ in $\cM$. 
In particular, $f$ is not isolated in $\cS_i^*$.

Consider a one-dimensional irreducible analytic subset $V \subset \cS_i^*$
containing $f$ and contained in a small neighborhood of $f$.
Since $\chi_i|_{\cS_i^* \cap \cR(\lanot)}$ is injective, 
$\iota_i(V)$ is transverse to the codimension one foliation $\cF$ of
hybrid classes at $f$ by \cite[Lemma~4.25]{Lyubich-ql-space}.
Moreover, $\iota_i(\cS_i)$ is contained in a leaf of $\cF$,
so it follows that there exists a tangent vector $v_i \in T_f\cS_i^*
\setminus T_f\cS_i$.

Given $i$, for all $j \neq i$, we have that $v_j \in T_f \cS^*_j \subset  T_f\cS_i$
and $v_i \notin T_f\cS_i$. Thus, for all $i$, we have that 
$v_i$ is not contained in the linear span of $\{ v_j; j \neq i\}$.
Hence, $\{ v_1, \dots, v_{d-1} \}$ are linearly independent vector in the $(d-1)$-dimensional
vector space $T_f\poly(d)$. That is, $v_1,\dots,v_{d-1}$ form a basis of $T_f\poly(d)$.

By repeating the same argument at a smooth point $g$ of $V$ with $g \in
\cR(\lanot)$, it follows that each component of $\cS_i^*$ containing
$f$ is one-dimensional. In particular, those components are smooth at $f$.
By Theorem~\ref{thm-DH-cont} and the injectivity of $\chi_i|_{\cS_i^*\cap
\cR(\lanot)}$, we conclude that there is only one such a component,
i.e., $\cS_i^*$ is a smooth curve near $f$.

Since $\iota_i$ is an analytic map transverse to the leaf $L$ of $\cF$
containing $\iota_i(f)$, it follows that $\cS_i=\iota_i^{-1}(L)$ is a
codimension one smooth complex manifold near $f$.
Moreover we have $T_f\cS_i=\bigoplus_{j\ne i} T_f\cS_j^*$,
so it follows that for any given $0 \le k \le d-1$ and
$1 \le i_1 < \cdots < i_k \le d-1$, 
the transversal intersection
\[
 \cS_{i_1,\dots,i_k}(f)=\cS_{i_1}(f) \pitchfork \cS_{i_2}(f) \pitchfork
 \dots \pitchfork \cS_{i_k}(f)
\]
is a complex submanifold of codimension $k$,
by shrinking $\cU$ if necessary. 

Now for any $0 \le k \le d-1$, given 
$1 \le i_1 < \cdots < i_k \le d-1$ and  $c_1, \dots, c_k \in \cM$, 
let 
\[
 \cC = \{ g \in \cR(\lanot);\ \chi_{i_j}(g) = c_j \mbox{ for all }
 j=1, \dots, k \}.
\]
If $f' \in \cC \cap \cS_{i_1,\dots,i_k}(f)$,
then $\cS_{i_1,\dots,i_k}(f')$ agrees with $\cS_{i_1,\dots,i_k}(f')$ in 
a small neighborhood of $f'$.
Therefore, $\tilde{\cS} = \bigcup_{f \in \cC} \cS_{i_1,\dots,i_k}(f)$ is a
codimension $k$ complex submanifold of $\poly(d)$ containing $\cC$ (by
shrinking $\tilde{\cS}$ if necessary).

Now consider the case $k=d-2$.
We may assume that $\tilde{S}$ is a one complex dimensional disk
$\tilde{S}$ embedded in $\poly(d)$ containing 
$$\{ g \in \cR(\lanot);\ \chi_j(g) = \chi_j(f) \mbox{ for all }j \ne i\}.$$ 
By the above construction, there exists a finite covering $\{ \cU_k \}$
of $\tilde{S}$ such that each $\cU_k$ parameterizes an analytic family of
quadratic like maps $$g_{(h,k)}= h^{\ell_i}: U'_i \to U_i,$$
where $h \in \cU_k$. 
Choosing $\tilde{S}$ sufficiently small, 
if $h \in \cU_k \cap \cU_{k'}$, then
the a priori distinct
 quadratic like germs $\iota_i(h) =[g_{(h,k)}]$ and  
$\iota_i(h)=[g_{(h,k')}]$ agree  as elements of $\cQG$.
Therefore, $\iota_i : \tilde{S} \to \cQG$ is well defined, analytic and
the image is transverse to hybrid classes.
Moreover, $\iota_i(h)$ belongs to the connectedness locus in $\cQG$ if and only if  the quadratic like map $h^{\ell_i}: U'_i \to U_i$ has connected Julia set, which is equivalent to $h \in \cR(\lanot) \cap \tilde{S}$, after applying Lemma~\ref{lem-hyperbolic thickening}. Thus the connectedness locus of   $\iota_i(\tilde{S})$ is compactly contained in  $\iota_i(\tilde{S})$.
That is, the family $\iota_i(\tilde{S})$ is ``full''~\cite[Section 4.11]{Lyubich-ql-space}.
The quasiconformality of $\chi_i$ in a neighborhood
of $\tilde{S} \cap \cR(\lanot)$ now follows 
directly from~\cite[Theorem~4.26]{Lyubich-ql-space}, since
$\chi:\tilde{S} \cap \cR(\lanot) \to \cM$ is injective (i.e. $\iota_i(\tilde{S})$ is an unfolded family in $\cQG$~\cite[Section 4.11]{Lyubich-ql-space}).
\qed

\subsubsection{Capture case}
In order to complete the proof of Theorem~\ref{introthm-capture-surjectivity}
we must establish that the straightening map involved is quasiconformal along appropriate
slices. 
The proof is similar than that of Buff and Henriksen in~\cite{Buff-Henriksen}: We construct a holomorphic motion on the filled
Julia set and use its quasiconformal extension, obtained via the $\lambda$-lemma, to
show that the straightening map has the desired quasiconformal regularity.

\begin{proof}[Proof of Theorem~\ref{introthm-capture-surjectivity}]
Let $c \in \cM$ and $f_0 \in \cR(\lanot)$ such that $\chi_1(f_0)=c$. Extend $\chi_1$  continuously to a neighborhood $\cU$ of
 $f_0$ in $\poly(3)$  (Lemma~\ref{lem-sufficient renor}).
Let $\cS'_c=\{f \in \cU; \chi_1(f)=c\}$.
Then from Douady-Hubbard's Theorem~\ref{thm-DH-cont} it follows that $\cS'_c$ is 
a complex analytic subspace of $\poly(3)$. 
By Theorem~\ref{thm-capture-surjectivity},
$\cS'_c$ is uncountable and not of full dimension. Hence it has dimension
 one.
Denote by $\cS$ an irreducible component (i.e. branch)  of $\cS'_c$ at $f_0$.
Note that  $\cS$ has dimension one.
We will show that $\chi_2$ is quasiconformal on $\cS$.
Then, since a quasiconformal map is an open map, from the injectivity of $\chi$
we will  conclude that there is at most one branch
of $\cS_c$ at $f_0$.
 Hence, it will immediately follow that $\cS'_c$ is locally irreducible.
Note that, without loss of generality, we may assume that there exists a holomorphic parameterization
$\pi : \Delta \rightarrow \cS \subset \poly(3)$. Since $\pi$ is a bijection, for simplicity, we identify $f \in \cS$ with
$\pi^{-1}(f)$.

Consider a family $(f^{\ell_{v_0}}:U_{v_0}' \to U_{v_0})_{f \in \Delta}$
 of quadratic-like mappings.
 Since the hybrid class of the elements of this family is constant,
 it is stable on the whole parameter space $\Delta$, in the sense of
 Ma\~{n}\'{e}-Sad-Sullivan.
 It follows that the Julia set $J_f(v_0)=\partial K_f(v_0)$ admits a
 holomorphic motion on $\Delta$ \cite{McMullen}.
 Furthermore, we have a natural conformal conjugacy
 $\psi_f^{-1} \circ \psi_{f_0}: \Int K_{f_0}(v_0) \to \Int K_f(v_0)$
 between $f_0$ and $f$, which is also a holomorphic motion on $\Delta$,
 where $\psi_f$ is a hybrid conjugacy between $f^{\ell_{v_0}}:U_{v_0}'
 \to U_{v_0}$ and $Q(z)=z^2+c$. In fact, on a periodic Fatou component,
 $\psi_f^{-1} \circ \psi_{f_0}$ can be also written locally as a
 composition $\phi_f^{-1} \circ \phi_{f_0}$ where $\phi_f$ is the
 (properly normalized) linearizing
 coordinate (if the component has a periodic point) or the
 Fatou coordinate (if the component is a parabolic basin).
 Since these coordinates depend holomorphically on $f \in \Delta$,
 $\psi_{f}^{-1} \circ \psi_{f_0}$ is also holomorphic on periodic
 components. On each preperiodic component, it is a pullback, so
 it is also analytic.

 Gluing both holomorphic motions together, we obtain a holomorphic motion
 $h_f: K_{f_0}(v_0) \to K_f(v_0)$.
 For $f \in \Delta \cap \cR(\lanot)$, we have
 \begin{equation*}
  \straightening (f) = (c,\psi_f (f^{\ell_{v_1}}(\omega_f)))
  = (c,\psi_{f_0} \circ h^{-1}_f (f^{\ell_{v_1}}(\omega_f))),
 \end{equation*}
 where $\omega_f$ is the captured critical point.

 By applying the $\lambda$-lemma \cite{MSS}
 \cite{Slodkowski}, $h_f$ extends to a holomorphic motion 
 \[
  h: \Delta \times \C \to \C,
 \]
 such that $h_f:\C \to \C$ is quasiconformal.
 By the same argument as \cite[Lemma~13]{Buff-Henriksen},
 the map
 \[
  f \mapsto h^{-1}_f(f^{\ell_{v_1}}(\omega_f))
 \]
 is locally quasiconformal on ${\Delta}$.
Hence, $\chi_2$ extends to a  quasiconformal map on a neighborhood of $f_0$ in $\cS$.
As described above, we conclude that $\cS$ is the unique branch of $\cS'_c$ at $f_0$.

 Passing to a finite cover of $\chi_1^{-1}(c)$, by open sets $\cU$ as above, we conclude that
there exists a one dimensional, locally irreducible, complex analytic space $\cS_c$ such that $\chi_1^{-1}(c) \subset \cS_c$
and $\chi_1 : \cR(\lanot) \cap \cS_c \to K(z^2+c)$ is a homeomorphism which extends locally to a quasiconformal map.

To finish the proof we have to show that $\cR(\lanot)$ is connected.
Let $C$ be both open and closed in 
$\cR(\lanot)$ and for $c \in \cM$ let $\cS_c$ be as above.
Then $\cS_c \cap \cR(\lanot)$, which is homeomorphic to the connected
set $K(z^2+c)$, is either contained in $C$ or disjoint
from $C$.
Therefore, $\chi_1(C)=\{c \in \cM;\ \cS_c \cap \cR(\lanot) \subset C\}$
and $\chi_1(\cR(\lanot) \setminus C)=\{c \in \cM;\ (\cS_c \cap
\cR(\lanot)) \cap C = \emptyset\}$ are closed
and disjoint subsets of $\cM$ because $\chi_1$ is continuous.
Then one of these sets must be empty, since $\cM$ is connected and,
 $\chi_1(C) \sqcup \chi_1(\cR(\lanot) \setminus C) = \cM$, by surjectivity.
 Hence, either $C$ is empty or $C= \cR(\lanot)$.
\end{proof}




\end{document}